\documentclass[10pt]{amsart}
 
\usepackage{amscd,a4,graphics,amssymb,amsbsy}
\usepackage[all]{xy}
\usepackage{color}

\newtheorem{thm}{Theorem}[section] \newtheorem{lemma}[thm]{Lemma}
\newtheorem{prop}[thm]{Proposition}
\newtheorem{cor}[thm]{Corollary}

\theoremstyle{definition}
 \newtheorem{dfn}[thm]{Definition}
\newtheorem{ntn}[thm]{Notation} \newtheorem{rmk}[thm]{Remark}
\newenvironment{pf}{\medskip\noindent{{\em Proof: }}}{\qed}
\newtheorem{ex}[thm]{Example}
\newtheorem{ass}[thm]{Assumption}

\newcommand{\Spec}{\mathrm{Spec}}
\newcommand{\Frac}{\mathrm{Frac}}
 
\newcommand{\Char}{\mathrm{char}}

\newcommand{\Gcd}{\mathrm{gcd}}
\newcommand{\Deg}{\mathrm{deg}}
\newcommand{\Mult}{\mathrm{mult}}
\newcommand{\Lcm}{\mathrm{lcm}}

\newcommand{\Sing}{\mathrm{Sing}}
\newcommand{\Aut}{\mathrm{Aut}}
\newcommand{\Id}{\mathrm{id}}
\newcommand{\Reg}{\mathrm{Reg}}
\newcommand{\Dim}{\mathrm{dim}}
\newcommand{\Ker}{\mathrm{Ker}}
\newcommand{\Rk}{\mathrm{rk}}
\newcommand{\Hom}{\mathrm{Hom}}
\newcommand{\Card}{\mathrm{Card}}
\newcommand{\Pic}{\mathrm{Pic}}
\newcommand{\Lim}{\mathrm{lim}}
\newcommand{\Tr}{\mathrm{Tr}}
\newcommand{\Sch}{\mathrm{Sch}}
\newcommand{\Sets}{\mathrm{Sets}}

\numberwithin{equation}{section}

\begin{document}


\title{Galois actions on models of curves}


  \author{Lars Halvard Halle}\email{halle@math.kth.se}
  \address{Department of Mathematics, KTH, S--100 44 Stockholm,
    Sweden}


\begin{abstract}
We study group actions on regular models of curves. If $X$ is a smooth curve defined over the fraction field $K$ of a complete discrete valuation ring $R$, every tamely ramified extension $K'/K$ with Galois group $G$ induces a $G$-action on the extension $X_{K'}$ of $X$ to $K'$. In this paper we study the extension of this $G$-action to certain regular models of $X_{K'}$. In particular, we are interested in the induced action on the cohomology groups of the structure sheaf of the special fiber of such a regular model. We obtain a formula for the Brauer trace of the endomorphism induced by a group element on the alternating sum of the cohomology groups. Inspired by this global study, we also consider similar group actions on the cohomology of the structure sheaf of the exceptional locus of a tame cyclic quotient singularity, and obtain an explicit polynomial formula for the Brauer trace of the endomorphism induced by a group element on the alternating sum of the cohomology groups. 

We apply these results to study a natural filtration of the special fiber of the N\'eron model of the Jacobian of $X$ by closed, unipotent subgroup schemes. We show that the jumps in this filtration only depend on the fiber type of the special fiber of the minimal regular model with strict normal crossings for $X$ over $\Spec(R)$, and in particular are independent of the residue characteristic. Furthermore, we obtain information about where these jumps occur. We also compute the jumps for each of the finitely many possible fiber type for curves of genus $1$ and $2$.   
\end{abstract}


\maketitle

\section{Introduction}

\subsection{Stable reduction of curves and Jacobians}
Let $X$ be a smooth, projective and geometrically connected curve of genus $ g(X) > 0 $, defined over the fraction field $K$ of a complete discrete valuation ring $R$, with algebraically closed residue field $k$. By a \emph{model} for $X$ over $R$, we mean an integral and normal scheme $ \mathcal{X} $ that is flat and projective over $ S = \Spec(R) $, and with generic fiber $ \mathcal{X}_K \cong X $. 

Recall the semi-stable reduction theorem, due to Deligne and Mumford (\cite{DelMum}, Corollary 2.7), which states that there exists a finite, separable field extension $ L/K$ such that $X_L$ admits a semi-stable model over the integral closure $R_L$ of $R$ in $L$.

It can often be useful to work with the Jacobian $J/K$ of $X$. The question whether $X$ has semi-stable reduction over $S = \Spec(R)$ is reflected in the structure of the \emph{N\'eron model} $ \mathcal{J}/S$ (cf. ~\cite{Ner}) of $J$. In fact, $X$ has semi-stable reduction over $S$ if and only if $ \mathcal{J}_k^0 $, the identity component of the special fiber, has no \emph{unipotent radical} (\cite{DelMum}, Proposition 2.3). 

In general, it is necessary to make ramified base extensions in order for $X$ to obtain semi-stable reduction. If the residue characteristic is positive, it can often be difficult to find explicit extensions over which $X$ obtains stable reduction. In the case where a tamely ramified extension suffices one can do this by considering the geometry of suitable regular models for $X$ over $ S $ (cf. ~\cite{Tame}). In this paper we study, among other things, how the geometry of the N\'eron model contains information that is relevant for obtaining semi-stable reduction for $X$.

\subsection{N\'eron models and tame base change}
Let $ K'/K $ be a finite, separable and tamely ramified extension of fields, and let $R'$ be the integral closure of $R$ in $ K' $. Then $R'$ is a complete discrete valuation ring, with residue field $k$. Furthermore, $ K'/K $ is Galois, with group $ G = \boldsymbol{\mu}_n $, where $ n = \Deg(K'/K) $.

Let $ \mathcal{J}'/S' $ be the N\'eron model of the Jacobian of $ X_{K'} $, where $ S' = \Spec(R') $. Due to a result by B. Edixhoven (\cite{Edix}, Theorem 4.2), it is possible to describe $ \mathcal{J}/S $ in terms of $ \mathcal{J}'/S' $, together with the induced $G$-action on $ \mathcal{J}' $. Namely, if $W$ denotes the \emph{Weil restriction} of $ \mathcal{J}'/S' $ to $S$ (cf. ~\cite{Ner}, Chapter 7), one can let $G$ act on $W$ in such a way that $ \mathcal{J} \cong W^G $, where $ W^G $ denotes the scheme of invariant points. In particular, one gets an isomorphism $ \mathcal{J}_k \cong W_k^G $. By \cite{Edix}, Theorem 5.3, this decription of $ \mathcal{J}_k $ induces a descending filtration 
$$ \mathcal{J}_k = F_n^0 \supseteq \ldots \supseteq F_n^i \supseteq \ldots \supseteq F_n^n = 0 $$
of $ \mathcal{J}_k $ by closed subgroup schemes.

In \cite{Edix}, Remark 5.4.5, a generalization of this setup is suggested. If we define $ \mathcal{F}^{i/n} = F_n^i $, where $ F_n^i $ is the $i$-th step in the filtration induced by the extension of degree $n$, one can consider the filtration 
$$ \mathcal{J}_k = \mathcal{F}^0 \supseteq \ldots \supseteq \mathcal{F}^a \supseteq \ldots \supseteq \mathcal{F}^1 = 0, $$
with indices in $ \mathbb{Z}_{(p)} \cap [0,1] $. In order for this to make sense, it is necessary that the construction of $ \mathcal{F}^a $ is independent of the choice of representatives $n$ and $i$ for $a$, and that $ \mathcal{F}^a $ descends with increasing $a$. We give a proof for these properties in Section \ref{Neron}. 

The filtration $ \{ \mathcal{F}^a \} $ contains significant information about $\mathcal{J}$. For instance, one can show that the subgroup schemes $ \mathcal{F}^a $ are \emph{unipotent} for $ a > 0 $. So in a natural way, this filtration gives a measure on how far $ \mathcal{J}/S $ is from being \emph{semi-abelian}. 

One way to study the filtration $ \{ \mathcal{F}^a \} $ is to determine where it \emph{jumps}. This will occupy a considerable part of this paper. The jumps in the filtration often give explicit numerical information about $X$. For instance, if $X$ obtains stable reduction after a tamely ramified extension, we show that the jumps occur at indices of the form $i/\tilde{n}$, where $\tilde{n}$ is the degree of the minimal extension that realizes stable reduction for $X$. 

It follows from Edixhoven's theory that to determine the jumps in the filtration $ \{ F^i_n \} $ induced by an extension of degree $n$, one needs to compute the irreducible characters for the representation of $ \boldsymbol{\mu}_n $ on the tangent space $ T_{ \mathcal{J}'_k, 0 } $. We shall use such computations for infinitely many integers $n$ to describe the jumps of the filtration $ \{ \mathcal{F}^a \} $ with rational indices.

\subsection{}
Contrary to the case of general abelian varieties, N\'eron models for Jacobians can be constructed in a fairly concrete way, using the theory of the relative Picard functor (cf. ~\cite{Ner}, Chapter 9). The following property will be of particular importance to us: If $ \mathcal{Z}/S' $ is a regular model for $ X_{K'}/K' $, then there is a canonical isomorphism
$$ \Pic_{\mathcal{Z}/S'}^0 \cong (\mathcal{J}')^0, $$  
where $ \Pic_{\mathcal{Z}/S'}^0 $ (resp. $(\mathcal{J}')^0$) is the identity component of $ \Pic_{\mathcal{Z}/S'} $ (resp. $\mathcal{J}'$). It follows that there is a canonical isomorphism
$$ H^1(\mathcal{Z}_k, \mathcal{O}_{\mathcal{Z}_k}) \cong T_{\mathcal{J}'_k, 0}. $$
 
We shall work with regular models $ \mathcal{Z} $ of $ X_{K'} $ that admit $G$-actions that are compatible with the $G$-action on $ \mathcal{J}' $. It will then follow that the representation of $G$ on $ T_{\mathcal{J}'_k, 0} $ can be described in terms of the representation of $G$ on $ H^1(\mathcal{Z}_k, \mathcal{O}_{\mathcal{Z}_k}) $.

\subsection{}
In order to find an $S'$-model for $X_{K'} $ with a compatible $G$-action, we take a model $ \mathcal{X} $ of $X$ over $S$, and consider its pullback $ \mathcal{X}_{S'} $ to $ S' $. This is in general not a model of $ X_{K'} $, but its normalization $ \mathcal{X}' $ will be a model. Furthermore, there exists a \emph{minimal} desingularization $ \mathcal{Y} \rightarrow \mathcal{X}' $ which is an isomorphism on generic fibers, thus producing in a natural way a regular model for $ X_{K'} $ over $S'$. 

There is a natural $G$-action on $ \mathcal{X}_{S'} $ via the action on the second factor. This action lifts uniquely to the normalization $ \mathcal{X}' $, and to the minimal desingularization $ \mathcal{Y} $. The $G$-action restricts to the special fiber $ \mathcal{Y}_k $, and in particular, $G$ will act on the cohomology groups $ H^i(\mathcal{Y}_k, \mathcal{O}_{\mathcal{Y}_k}) $, for $ i = 0, 1 $. 

In order to understand the $G$-action on $ H^i(\mathcal{Y}_k, \mathcal{O}_{\mathcal{Y}_k}) $, it is important that we have a good description of the geometry of $\mathcal{Y}$ and of the $G$-action on $ \mathcal{Y} $, and this is studied in Section \ref{extensions and actions} and Section \ref{lifting the action}. For this purpose, we demand that the model $ \mathcal{X} $ has good properties. To begin with, we shall require that $ \mathcal{X} $ is regular, and that the special fiber is a divisor with strict normal crossings, i.e., that $ \mathcal{X} $ is an SNC-model. However, we have to impose some restrictions on the geometry of $ \mathcal{X}_k $. In fact, we shall always require that any two irreducible components of $ \mathcal{X}_k $, whose multiplicities are both divisible by the residue characteristic, have empty intersection. This condition is automatically fulfilled if $X$ obtains stable reduction after a tamely ramified extension, but holds also for a larger class of curves.  

Under these assumptions, it turns out that the normalization $ \mathcal{X}' $ of $ \mathcal{X}_{S'} $ has rather well behaved singularities, known as \emph{tame cyclic quotient singularities} (cf. ~\cite{CED}, Definition 2.3.6 and \cite{Tame}, Proposition 4.3). Furthermore, these singularities can be resolved explicitly, and it can be seen that the minimal desingularization $ \mathcal{Y} $ is an SNC-model for $ X_{K'} $. 

We shall also make the assumption that $ n = \Deg(K'/K) $ is relatively prime to the multiplicities of all the irreducible components of $ \mathcal{X}_k $. With this additional hypothesis, it turns out that we can describe the combinatorial structure of the special fiber $ \mathcal{Y}_k $ (i.e., the intersection graph of the irreducible components, their genera and multiplicities), in terms of the corresponding data for $ \mathcal{X}_k $. 

If all the assumptions above are satisfied, it follows that all irreducible components of $ \mathcal{Y}_k $ are stable under the $G$-action on $ \mathcal{Y} $, and that all intersection points in $ \mathcal{Y}_k $ are fixed points. We can explicitly describe the action on the cotangent space of $ \mathcal{Y} $ at these intersection points, and the restriction of the $G$-action to each irreducible component of $ \mathcal{Y}_k $. 

\subsection{Action on cohomology}
Next, we study the representation of $ G = \boldsymbol{\mu}_n $ on $ H^1(\mathcal{Y}_k, \mathcal{O}_{\mathcal{Y}_k}) $. In particular, we would like to compute the irreducible characters for this representation. So for every $ g \in G $, we want to compute the trace of the endomorphism of $ H^1(\mathcal{Y}_k, \mathcal{O}_{\mathcal{Y}_k}) $ induced by $g$, and then use this information to find the characters.

There are some technical problems that need to be overcome in order to do this. First, since we allow the residue characteristic to be positive, just knowing the trace for each $g \in G$ may not give sufficient information to compute the characters. Instead, we have to compute the so called \emph{Brauer trace} for every $g \in G $ (cf. \cite{SerreLin}, Chapter 18). This means that we have to lift the eigenvalues and traces from characteristic $p$ to characteristic $0$. From knowing the Brauer trace for every $ g \in G $ we can compute the irreducible Brauer characters, and then the ordinary characters are obtained by reducing the Brauer characters modulo $p$. Second, the special fiber $ \mathcal{Y}_k $ will in general be singular, and even non-reduced. This complicates trace computations considerably.

To deal with these problems, we introduce in Section \ref{section 6} a certain filtration of the special fiber $ \mathcal{Y}_k $ by effective subdivisors, where the difference at the $i$-th step is an irreducible component $ C_i $ of $ \mathcal{Y}_k $. Since $ \mathcal{Y} $ is an SNC-model, each $ C_i $ is a smooth and projective curve, and with our assumption on $n$, the $G$-action restricts to each $ C_i $. Furthermore, to each step in this filtration, one can in a natural way associate an invertible $ G $-sheaf $ \mathcal{L}_i $, supported on $C_i$. 

We apply the so called Lefschetz-Riemann-Roch formula (\cite{Don}, Corollary 5.5), in order to get a formula for the Brauer trace of the endomorphism induced by each $ g \in G $ on the formal difference $ H^0(C_i,\mathcal{L}_i) - H^1(C_i,\mathcal{L}_i) $. An important step is to show that our description of the action on $ \mathcal{Y} $ is precisely the data that is needed to obtain these formulas. Then we show that these traces add up to give the Brauer trace for the endomorphism induced by each $ g \in G $ on the formal difference $ H^0(\mathcal{Y}_k, \mathcal{O}_{\mathcal{Y}_k}) - H^1(\mathcal{Y}_k, \mathcal{O}_{\mathcal{Y}_k}) $. In particular, we give in Theorem \ref{thm. 9.13}, which is the first main result in this work, a formula for this Brauer trace, and show that it only depends on the combinatorial structure of $ \mathcal{X}_k $. 

Let us also remark that in our applications, we already know the character for $ H^0(\mathcal{Y}_k, \mathcal{O}_{\mathcal{Y}_k}) $, and hence we will be able to compute the irreducible characters for $ H^1(\mathcal{Y}_k, \mathcal{O}_{\mathcal{Y}_k}) $ in this way.

\subsection{Trace formulas for singularities}
If $ x' \in \mathcal{X}' $ is a singular point, the $G$-action on $ \mathcal{Y} $ restricts to the exceptional fiber $ \mathcal{E}_{x'} := \rho^{-1}(x') \subset \mathcal{Y}_k $, where $ \rho : \mathcal{Y} \rightarrow \mathcal{X}' $ is the minimal desingularization. Hence $ G = \boldsymbol{\mu}_n $ acts on the cohomology groups $ H^i(\mathcal{E}_{x'}, \mathcal{O}_{\mathcal{E}_{x'}}) $, for $ i = 0,1 $. This situation is studied in Section \ref{trace formula}. We observe that the methods developed earlier in the paper also apply to this situation, and enable us to compute the Brauer trace of the endomorphism induced by $ g \in G $ on the formal difference $ H^0(\mathcal{E}_{x'}, \mathcal{O}_{\mathcal{E}_{x'}}) - H^1(\mathcal{E}_{x'}, \mathcal{O}_{\mathcal{E}_{x'}}) $. 

The singularity $ x' \in \mathcal{X}' $ is determined by \emph{parameters} $ n $, $m_1$ and $m_2$, where $ n $ is the order of $G$, and $ m_1$ and $ m_2 $ are the multiplicities of the components of $ \mathcal{X}_k' $ intersecting at $x'$. If $n$ is large enough compared to $m_1$ and $m_2$, we obtain in Theorem \ref{Formula} an explicit closed polynomial formula for the Brauer trace in terms of the parameters of the singularity. It turns out that there is precisely one polynomial for each element in $ (\mathbb{Z}/M)^* $, where $M = \Lcm(m_1,m_2) $.

Apart from the fact that we find this to be an interesting problem in its own right, these formulas are used later on in the paper, both for theoretical issues as well as for the explicit computations in Section \ref{computations and jumps}. In particular, by combining Theorem \ref{Formula} and Theorem \ref{thm. 9.13}, we obtain our main result Theorem \ref{improved formula}, which gives an explicit effective formula for the Brauer trace of the automorphism induced by $\xi \in \boldsymbol{\mu}_n $ on the alternating sum $ H^0(\mathcal{Y}_k, \mathcal{O}_{\mathcal{Y}_k}) - H^1(\mathcal{Y}_k, \mathcal{O}_{\mathcal{Y}_k}) $.

\subsection{}
If now $ \mathcal{X}/S $ is the minimal SNC-model for $ X/K $, we prove in Theorem \ref{main character theorem} that the irreducible characters for the representation of $ G = \boldsymbol{\mu}_n $ on $ H^1(\mathcal{Y}_k, \mathcal{O}_{\mathcal{Y}_k}) $ only depend on the combinatorial structure of the special fiber $ \mathcal{X}_k $, as long as $n$ is relatively prime to $l$, where $l$ is the least common multiple of the multiplicities of the irreducible components of $ \mathcal{X}_k $. 

If $ \{ \mathcal{F}^a \} $ denotes the filtration of $ \mathcal{J}_k $, where $ \mathcal{J} $ is the N\'eron model of the Jacobian of $X$, we prove as a corollary that the jumps in the filtration $ \{ \mathcal{F}^a \} $ only depend on the combinatorial structure of $ \mathcal{X}_k $ (Corollary \ref{main jump corollary}). This is due to the fact that $ \mathbb{Z}_{ (p l) } \cap [0,1] $ is ``dense'' in $ \mathbb{Z}_{ (p) } \cap [0,1] $. Furthermore, in Corollary \ref{specific jump corollary}, we draw the conclusion that the jumps are actually independent of ~$p$, and that the jumps can only occur at finitely many rational numbers of a certain kind, depending on the combinatorial structure of $ \mathcal{X}_k $. 

It is known that for a fixed genus $ g \geq 1 $, there are only a finite number of possible combinatorial structures for $ \mathcal{X}_k $, modulo a certain equivalence relation. Furthermore, in case $ g = 1 $ or $ g = 2 $, one has complete classifications (cf. ~\cite{Kod} for $g=1$ and \cite{Ueno}, \cite{Ogg} for $g=2$). In Section \ref{computations and jumps} we compute the jumps for each possible fiber type for $g=1$ (which were also computed by Schoof in \cite{Edix}) and for $g=2$.

\subsection{Acknowledgements}
I would like to thank Bas Edixhoven for suggesting this subject to me, and generously sharing his ideas. I would also like to thank my thesis advisor Carel Faber for discussing the material in this paper with me.

\section{N\'eron models and tamely ramified extensions}\label{Neron}

\subsection{N\'eron models} 
Let $ R $ be a discrete valuation ring, with fraction field $K$ and residue field $k$, and let $A$ be an abelian variety over $K$. There exists a canonical extension of $A$ to a smooth group scheme $ \mathcal{A} $ over $ S = \Spec(R) $, known as the \emph{N\'eron model} (\cite{Ner}, Theorem 1.4/3). The N\'eron model is characterized by the following universal property: for every smooth morphism $ T \rightarrow S $, the induced map $ \mathcal{A}(T) \rightarrow A(T_K) $ is \emph{bijective}.

\subsection{N\'eron models and base change}\label{neronbase}
We assume from now on that $R$ is strictly henselian. Let $ K'/K $ be a finite, separable extension of fields, and let $R'$ be the integral closure of $R$ in $K'$. Let $ \mathcal{A}'/S' $ denote the N\'eron model of the abelian variety $ A_{K'}/K' $, where $ S' = \Spec(R') $. In general, it is not so easy to describe how N\'eron models change under ramified base changes. However, in the case where $K'/K$ is tamely ramified, one can relate $ \mathcal{A}'/S' $ to $ \mathcal{A}/S $ in a nice way, due to a result by B. Edixhoven (\cite{Edix}, Theorem 4.2). We will in this section explain this relation, following the treatment in \cite{Edix}. We refer to this paper for further details. 

From now on, we will assume that $ K'/K $ is tamely ramified. Then $ K'/K $ is Galois with group $ G = \boldsymbol{\mu}_n $, where $n$ is the degree of the extension. Let $G$ act on $ A_{K'} = A \times_{\Spec(K)} \Spec(K') $ (from the right), via the action on the right factor. By the universal property of $ \mathcal{A}' $, this $G$-action on $ A_{K'} $ extends uniquely to a right action on $ \mathcal{A}' $, such that the morphism $ \mathcal{A}' \rightarrow S' $ is equivariant. The idea is now to reconstruct $ \mathcal{A} $ as an invariant scheme for this action. However, since $ \mathcal{A} $ is an $S$-scheme, one first needs to ``push down'' from $S'$ to $S$. 

The \emph{Weil restriction} of $ \mathcal{A}' $ to $S$ is a contravariant functor 
$$ \Pi_{S'/S}(\mathcal{A}'/S') : (\Sch/S)^0 \rightarrow (\Sets), $$ 
defined by assigning, for any $S$-scheme $T$, $ \Pi_{S'/S}(\mathcal{A}'/S')(T) = \mathcal{A}'(T') $, where $ T' = T \times_S S' $ (cf. ~\cite{Ner}, Chapter 7). This functor is representable by an $S$-scheme, which we will denote by $X$ (\cite{Edix}, Remark 2.1).   

In \cite{Edix}, an equivariant $G$-action on $ X \rightarrow S $ is defined in the following way: Let $ T $ be an $S$-scheme, and let $ P \in X(T) = \mathcal{A}'(T') $ be a $T$-point, where $ T' = T \times_S S' $. For any $ g \in G $, let $ \rho_{\mathcal{A}'}(g) $ and $ \rho_{S'}(g) $ be the automorphisms of $ \mathcal{A}' $ and $S'$ induced by $g$. Finally, let $ \rho_{T'}(g) = 1_T \times \rho_{S'}(g) $. The action is then given by
$$ P \cdot g = \rho_{\mathcal{A}'}(g) \circ P \circ \rho_{T'}(g)^{-1}. $$

\subsection{Relating $\mathcal{A}'$ and $\mathcal{A}$}
In \cite{Edix}, Proposition 4.1, it is shown that the generic fiber of $ \Pi_{S'/S}(\mathcal{A}'/S') $ is an abelian variety, having $ \Pi_{S'/S}(\mathcal{A}'/S') $ as its N\'eron model. Furthermore, Construction 2.3 in \cite{Edix} gives a closed immersion
$$ A \hookrightarrow (\Pi_{S'/S}(\mathcal{A}'/S'))_K, $$
which extends uniquely to a morphism
$$ \mathcal{A} \rightarrow \Pi_{S'/S}(\mathcal{A}'/S'), $$
since $ \mathcal{A}/S $ is smooth and $ \Pi_{S'/S}(\mathcal{A}'/S') $ is a N\'eron model.

According to \cite{Edix}, Theorem 4.2, this morphism is a closed immersion, and induces an isomorphism
\begin{equation}\label{Bastheorem}
\mathcal{A} \cong (\Pi_{S'/S}(\mathcal{A}'/S'))^G, 
\end{equation}
where $(\Pi_{S'/S}(\mathcal{A}'/S'))^G$ is the \emph{scheme of invariant points} for the $ G $-action defined above (\cite{Edix}, Chapter 3). 

\subsection{Filtration of $\mathcal{A}_k$} 
One can use the isomorphism in \ref{Bastheorem} to study the special fiber $ \mathcal{A}_k $ in terms of $ \mathcal{A}'_k $, together with the $G$-action. Indeed, let $ R \subset R' = R[\pi']/(\pi'^n - \pi) $ be a tame extension, where $ \pi $ is a uniformizing parameter for $R$. Then we have that $ R'/\pi R' = k[\pi']/(\pi'^n) $. For any $ k $-algebra $C$, it follows that 
$$ \mathcal{A}_k(C) \cong X_k^G(C) \cong X_k(C)^G \cong \mathcal{A}'(C[\pi']/(\pi'^n))^G. $$

In \cite{Edix}, Chapter 5, this observation is used to construct a filtration of $\mathcal{A}_k$. To do this, let us first consider an $R$-algebra $ C $. The Weil restriction induces a map
$$ \mathcal{A}(C) \rightarrow \mathcal{A}'(C \otimes_R R'), $$
which gives a map 
$$ \mathcal{A}(C) \rightarrow \mathcal{A}'(C \otimes_R R') \rightarrow \mathcal{A}'(C \otimes_R R'/(\pi'^i)), $$
for any integer $i$ such that $ 0 \leq i \leq n $. Define functors $ F^i\mathcal{A}_k $ by
$$ F^i\mathcal{A}_k(C) = \Ker(\mathcal{A}(C) \rightarrow \mathcal{A}'(C \otimes_R R'/(\pi'^i))), $$
for any $ k = R/(\pi)$-algebra $C$. The functors $ F^i\mathcal{A}_k $ are represented by closed subgroup schemes of $ \mathcal{A}_k $, and give rise to a descending filtration

$$ \mathcal{A}_k = F^0\mathcal{A}_k \supseteq F^1\mathcal{A}_k \supseteq \ldots \supseteq F^n\mathcal{A}_k = 0. $$
Let us also remark that the group schemes $ F^i\mathcal{A}_k $ are \emph{unipotent} for $ i > 0 $.

One can describe the successive quotients of this filtration quite accurately: Let $ Gr^i \mathcal{A}_k $ denote the quotient $ F^i\mathcal{A}_k/F^{i+1}\mathcal{A}_k $, for $ i \in \{ 0, \ldots, n-1 \} $. Then, by Theorem 5.3 in \cite{Edix}, we have that $ Gr^0(\mathcal{A}_k) = (\mathcal{A}'_k)^{\boldsymbol{\mu}_n} $, and for $ 0 < i < n $, we have that
$$ Gr^i \mathcal{A}_k \cong T_{\mathcal{A}'_k,0}[i] \otimes_k (m/m^2)^{\otimes i}, $$
where $ m \subset R'$ is the maximal ideal, and where $ T_{\mathcal{A}'_k,0}[i] $ denotes the subspace of $ T_{\mathcal{A}'_k,0} $ where $ \xi \in \boldsymbol{\mu}_n $ acts by multiplication by $ \xi^i $.

The filtration \emph{jumps} at the index $ i \in \{ 0, \ldots, n-1 \} $ if $ Gr^i \mathcal{A}_k \neq 0 $. Since 
$$ T_{\mathcal{A}'_k,0}[0] = (T_{\mathcal{A}'_k,0})^{\boldsymbol{\mu}_n} = T_{(\mathcal{A}'_k)^{\boldsymbol{\mu}_n},0} $$
(use \cite{Edix}, Proposition 3.2), it follows that the jumps are completely determined by the representation of $ \boldsymbol{\mu}_n $ on $ T_{\mathcal{A}'_k,0} $. In particular, it follows that there are at most $ \Dim(A) $ jumps, since $ \Dim_k T_{\mathcal{A}'_k,0} = \Dim(A) $.

\subsection{Compositions of tame extensions}
In Remark 5.4.5 \cite{Edix}, a generalization of the filtration of $ \mathcal{A}_k $ discussed above is suggested. Let $ \{ F^i_n \mathcal{A}_k \} $ denote the filtration induced by the tame extension of degree $n$. The idea is to put all the $ F^i_n \mathcal{A}_k $, for all positive integers not divisible by the residue characteristic $p$, in a common filtration of $ \mathcal{A}_k $. We shall later see that this filtration gives interesting information about $ \mathcal{A}_k $.

In order to set this up, it is necessary to understand how two filtrations $ \{ F^i_{n_1} \mathcal{A}_k \} $ and $ \{ F^i_{n_2} \mathcal{A}_k \} $ are related in the case where $ n_1 $ divides $n_2$. We will now explain this in some detail, since this is not done in \cite{Edix}.

Let $ R_1 $ and $ R_2 $ be tame extensions of $R$, where $ R_1 = R[t_1]/(t_1^{n_1} - \pi) $, $ R_2 = R[t_2]/(t_2^{n_2} - \pi) $ and where $ R_2 = R_1[t_2]/(t_2^{m} - t_1) $. Let $ \mathcal{A}_i $ denote the N\'eron model of $ A_{K_i} $ over $ R_i $ for $ i \in {1,2} $, where $ K_i $ is the fraction field of $R_i$. By Remark 2.1 in \cite{Edix}, we have that
$$ \Pi_{S_i/S} (\mathcal{A}_i/S_i) \times_S \Spec(K) = \Pi_{K_i/K} (A_{K_i}/K_i). $$

Furthermore, Construction 2.3 in \cite{Edix} gives, for every $ T \rightarrow \Spec(K) $, a commutative diagram
$$ 
\xymatrix{
\Hom_K(T,A) \ar[r] \ar[dr] & \Hom_{K_1}(T_{K_1}, A_{K_1}) \ar[d] \\
   &  \Hom_{K_2}(T_{K_2}, A_{K_2}), }
$$
where all arrows are injections. So we get a commutative diagram
$$ \xymatrix{
A \ar[rr]^-{\phi_1} \ar[drr]_-{\phi_2} & & \Pi_{K_1/K} (A_{K_1}/K_1) \ar[d]^{\phi_{1,2}} \\
 & & \Pi_{K_2/K} (A_{K_2}/K_2). }
$$
We now claim that all maps in this diagram are closed immersions. To see this, let $ G $ denote the Galois group of $ K_2/K $, and let $ H \subseteq G $ be the subgroup that fixes $K_1 $. Then $ H $ is the Galois group of $ K_2/K_1 $, and $ G/H $ is the Galois group of $K_1/K $. 

We have that $ G $ acts on $ A_{K_2} = A \times_{\Spec(K)} \Spec(K_2) $ by its action on the right factor. Therefore, $G$ acts on $ \Pi_{K_2/K} (A_{K_2}/K_2) $ as defined in Section \ref{neronbase}, with invariant scheme $A$ (see the proof of \cite{Edix}, Theorem 4.2), and $ \phi_2 $ can be identified with the inclusion. By Proposition 3.1 in \cite{Edix}, it follows that $\phi_2$ is a closed immersion. The same argument for the $G/H$-action on $ \Pi_{K_1/K} (A_{K_1}/K_1) $ shows that $\phi_1$ is a closed immersion.

\begin{lemma}
The morphism $ \phi_{1,2} $ above is a closed immersion. 
\end{lemma}
\begin{pf}
We need to understand how $ H $ acts on $ \Pi_{K_2/K} (A_{K_2}/K_2) $. For any $ g \in G $, let $ \rho_{K_2}(g) : \Spec(K_2) \rightarrow \Spec(K_2) $ be the corresponding automorphism. For any $ K $-scheme $T$, there is an induced automorphism given as $ \rho_{T_{K_2}}(g) = \Id_T \times \rho_{K_2}(g) $ on $ T_{K_2} $, and similarly, there is an induced automorphism $ \rho_{A_{K_2}}(g) = \Id_{A} \times \rho_{K_2}(g) $ on $ A_{K_2} $. Then $G$ acts, for any $ P \in \Pi_{K_2/K} (A_{K_2}/K_2)(T) $, by $ P \cdot g = \rho_{A_{K_2}}(g) \circ P \circ \rho_{T_{K_2}}(g)^{-1} $. 

However, if $ h \in H $, the morphism $ \rho_{K_2}(h) $ is relative to $ \Spec(K_1) $. So we can in fact write $ \rho_{T_{K_2}}(h) = \Id_{T_{K_1}} \times \rho_{K_2}(h) $ and $ \rho_{A_{K_2}}(g) = \Id_{A_{K_1}} \times \rho_{K_2}(g) $. 

Observe that we have the identity
$$ \Pi_{K_2/K} (A_{K_2}/K_2)(T) = \Hom_{K_2}(T_{K_2}, A_{K_2}) = \Pi_{K_2/K_1} (A_{K_2}/K_2)(T_{K_1}), $$
and an $H$-action on $ \Pi_{K_2/K_1} (A_{K_2}/K_2) $ as defined in Section \ref{neronbase}. It follows immediately that the $H$-action on $ \Pi_{K_2/K} (A_{K_2}/K_2)(T) $ and $ \Pi_{K_2/K_1} (A_{K_2}/K_2)(T_{K_1}) $ can be identified. 

But since $A_{K_1} $ is the invariant scheme of $ \Pi_{K_2/K_1} (A_{K_2}/K_2) $ under this $H$-action, we have that 
$$ \Pi_{K_2/K_1} (A_{K_2}/K_2)(T_{K_1})^H = \Hom_{K_1}(T_{K_1}, A_{K_1}) = \Pi_{K_1/K} (A_{K_1}/K_1)(T), $$
and it therefore follows that
$$ \Pi_{K_2/K} (A_{K_2}/K_2)(T)^H = \Hom_{K_1}(T_{K_1}, A_{K_1}) = \Pi_{K_1/K} (A_{K_1}/K_1)(T). $$ 
Therefore, since $ \Pi_{K_1/K} (A_{K_1}/K_1) $ is the invariant scheme of $ \Pi_{K_2/K} (A_{K_2}/K_2) $ under the action of $H$, and $ \phi_{1,2} $ can be identified with the inclusion, it follows that $ \phi_{1,2} $ is a closed immersion.
\end{pf}

\vspace{0.5cm}

By the N\'eronian property, the maps $ \phi_1 $, $ \phi_2 $ and $ \phi_{1,2} $ lift uniquely so that we have a commutative diagram
\begin{equation}\label{GH diagram}
 \xymatrix{
\mathcal{A} \ar[r]^-{\Phi_1} \ar[dr]_-{\Phi_2} &  \Pi_{S_1/S} (\mathcal{A}_1/S_1) \ar[d]^{\Phi_{1,2}} \\
& \Pi_{S_2/S} (\mathcal{A}_2/S_2). }
\end{equation}

\begin{lemma}
The morphisms $ \Phi_1 $, $ \Phi_2 $ and $ \Phi_{1,2} $ are closed immersions. Furthermore, via these maps, we have that $ \mathcal{A} = \Pi_{S_2/S} (\mathcal{A}_2/S_2)^G $, $ \mathcal{A} = \Pi_{S_1/S} (\mathcal{A}_1/S_1)^{G/H} $ and $ \Pi_{S_1/S} (\mathcal{A}_1/S_1) = \Pi_{S_2/S} (\mathcal{A}_2/S_2)^H $.
\end{lemma}
\begin{pf}
Except for the statement regarding $ \Phi_{1,2} $, this is Theorem 4.2 in \cite{Edix}. But the same proof can be applied to $ \Phi_{1,2} $.
\end{pf}

\vspace{0.5cm}

Having established this relationship between the Weil restrictions, we get the following result for the filtrations:

\begin{lemma}\label{12extensions}
Let $ \{ F_{n_j}^i \mathcal{A}_k \} $ be the filtration of $ \mathcal{A}_k $ induced by the extension $ R_j/R $, for $ j \in \{ 1,2 \} $. Then we have that $ F_{n_1}^i \mathcal{A}_k = F_{n_2}^{im} \mathcal{A}_k $, for any $ i \in \{ 0, \ldots, n_1 \} $.
\end{lemma}
\begin{pf}
Let $C$ be an $R$-algebra. As $ C \otimes_R R_2 = (C \otimes_R R_1) \otimes_{R_1} R_2 $, Diagram \ref{GH diagram} above induces a commutative diagram
\begin{equation}\label{diagram 14.4}
 \xymatrix{
\mathcal{A}(C) \ar[r]^-{\Phi_1} \ar[dr]_-{\Phi_2} & \mathcal{A}_1(C \otimes_R R_1) \ar[d]^{\Phi_{1,2}} \\
                                     & \mathcal{A}_2((C \otimes_R R_1) \otimes_{R_1} R_2). } 
\end{equation}
Note that $ C \otimes_R R_j = C [t_j]/(t_j^{n_j} - \pi) $, and in the case where $ C $ is also an $ R/(\pi) $-algebra, we get that $ C \otimes_R R_j = C [t_j]/(t_j^{n_j}) $. In the following, we shall only consider the latter case. Let $i$ be an integer such that $ 0 \leq i \leq n_1 $. Then we have that $ C[t_1]/(t_1^i) \otimes_{R_1} R_2 = C[t_2]/(t_2^{im}) $, and that there is a commutative diagram
\begin{equation}\label{diagram 14.5}
\xymatrix{
\mathcal{A}_1(C[t_1]/(t_1^{n_1})) \ar[r]^{\beta_1} \ar[d]_{\Phi_{1,2}} & \mathcal{A}_1(C[t_1]/(t_1^i)) \ar[d]^{\Phi_{1,2}} \\
\mathcal{A}_2(C[t_2]/(t_2^{n_2})) \ar[r]_{\beta_2}  & \mathcal{A}_2(C[t_2]/(t_2^{im})). } 
\end{equation}
 
Combining Diagram \ref{diagram 14.4} with Diagram \ref{diagram 14.5} gives a commutative diagram
\begin{equation}\label{diagram 14.6}
\xymatrix{
\mathcal{A}_k(C) \ar[r]^-{\alpha_1} \ar[dr]_{\alpha_2} & \mathcal{A}_1(C[t_1]/(t_1^i)) \ar[d]^-{\Phi_{1,2}} \\
                          & \mathcal{A}_2(C[t_2]/(t_2^{im})). } 
\end{equation}
Since $\Phi_{1,2}$ is \emph{injective}, it follows that $ \Ker(\alpha_1) = \Ker(\alpha_2) $ for any $k$-algebra $C$, and hence $ F_1^i\mathcal{A}_k = F_2^{im}\mathcal{A}_k $ for all $i$ such that $ 0 \leq i \leq n_1 $.
\end{pf}

\subsection{Filtration with rational indices}\label{ratfil}

Let $ a \in \mathbb{Z}_{(p)} \cap [0,1] $. If $ a = i/n $, then we define $ \mathcal{F}^a\mathcal{A}_k = F^i_n \mathcal{A}_k $, where $ F^i_n \mathcal{A}_k $ denotes the $i$-th step in the filtration induced by the tame extension of degree $n$.

\begin{prop}\label{ratfil}
The construction above gives a descending filtration
$$ \mathcal{A}_k = \mathcal{F}^0 \mathcal{A}_k \supseteq \ldots \supseteq \mathcal{F}^a\mathcal{A}_k \supseteq \ldots \supseteq \mathcal{F}^1 \mathcal{A}_k = 0 $$
of $ \mathcal{A}_k $ by closed subgroup schemes, where $ a \in \mathbb{Z}_{(p)} \cap [0,1] $.
\end{prop}
\begin{pf}
It follows from Lemma \ref{12extensions} that the definition of $ \mathcal{F}^a\mathcal{A}_k $ does not depend on the choice of representatives for $a$. Hence the $ \mathcal{F}^a\mathcal{A}_k $ are well defined. 

In order to show that the filtration is descending, take $a_1, a_2 \in \mathbb{Z}_{(p)} \cap [0,1] $ such that $ a_1 \leq a_2 $. Let $ a_1 = i_1/n_1 $ and $ a_2 = i_2/n_2 $. Since $ \mathcal{F}^a \mathcal{A}_k $ does not depend on the choice of representatives for $a$, we write $ a_1 = i_1 n_2/n_1 n_2 $ and $ a_2 = i_2 n_1/n_1 n_2 $. In particular, $ i_1 n_2 \leq i_2 n_1 $. But then we get that
$$ \mathcal{F}^{a_1}\mathcal{A}_k = F^{i_1 n_2}_{n_1 n_2} \mathcal{A}_k \supseteq F^{i_2 n_1}_{n_1 n_2} \mathcal{A}_k = \mathcal{F}^{a_2}\mathcal{A}_k. $$
\end{pf}

\vspace{0.3cm}

 Let $ x \in [0,1] $ be a real number, and let $ (x^j)_j $ (resp. $ (x_k)_k $) be a sequence of numbers in $ \mathbb{Z}_{(p)} \cap [0,1] $ converging to $x$ from above (resp. from below). We will say that $ \{ \mathcal{F}^a \mathcal{A}_k \} $ jumps at $x$ if $ \mathcal{F}^{x_k} \mathcal{A}_k \supsetneq \mathcal{F}^{x^j} \mathcal{A}_k $ for all $j$ and $k$. It is natural to ask how many jumps there are, and \emph{where} they occur. It is easily seen that since every filtration $ \{ F^i_n \mathcal{A}_k \} $ jumps at most $ g = \Dim(A) $ times, it follows that the filtration $ \{ \mathcal{F}^a \mathcal{A}_k \} $ can have at most $g$ jumps. 

Consider a positive integer $n$ that is not divisible by $p$, and let $ \{ F^i_n \mathcal{A}_k \} $ be the filtration induced by the extension of degree $n$. Let us assume that this filtration has a jump at $ i \in \{ 0, \ldots, n-1 \} $. Then we can say that $ \{ \mathcal{F}^a \mathcal{A}_k \} $ has a jump in the interval $ [i/n, (i+1)/n] $. By computing jumps in this way for increasing $n$, we get finer partitions of the interval $ [0,1] $, and increasingly better approximations of the jumps in $ \{ \mathcal{F}^a \mathcal{A}_k \} $. 

It follows that one can compute the jumps of $ \{ \mathcal{F}^a \mathcal{A}_k \} $ by computing the jumps for the filtrations $ \{ F^i_n \mathcal{A}_k \} $ for ``sufficiently'' many $n$ that are not divisible by $p$. This would for instance be the case for a multiplicatively closed subset $ \mathcal{U} \subset \mathbb{N} $ such that $ \mathbb{Z}[\mathcal{U}^{-1}] \cap [0,1] $ is dense in $ \mathbb{Z}_{(p)} \cap [0,1] $. 

\subsection{Jumps in the tamely ramified case}
In the case where $ A/K $ obtains semi-abelian reduction over a tamely ramified extension $K'$ of $K$, the jumps of $ \{ \mathcal{F}^a \mathcal{A}_k \} $ have an interesting interpretation, which we will now explain. Let $ \widetilde{K} $ be the minimal extension over which $A$ has semi-abelian reduction (cf. ~\cite{Deschamps}, Th\'eor\`eme 5.15), and let $ \tilde{n} = \Deg(\widetilde{K}/K) $. Then the jumps occur at rational numbers of the form $ k/\tilde{n} $, where $ k \in \{0, \ldots, \tilde{n} - 1 \} $. This is essentially due to the following observation:


\begin{lemma}\label{tame jumps}
Let $ \widetilde{K}/K $ be the minimal extension over which $A/K$ obtains semi-abelian reduction, and let $ \tilde{n} = \Deg(\widetilde{K}/K) $. Consider a tame extension $ K'/K $ of degree $ n $, factoring via $ \widetilde{K} $, and let $ m = n/\tilde{n} $. Let $\mathcal{A}'/S' $ be the N\'eron model of $A_{K'}$. 

Then we have that the jumps in the filtration $ \{ F^i_n \mathcal{A}_k \} $ induced by $ S'/S $ occur at indices $ i = k n / \tilde{n} $, where $ 0 \leq k \leq \tilde{n} - 1$.
\end{lemma}
\begin{pf}
Let $ \widetilde{\mathcal{A}}/\widetilde{S} $ be the N\'eron model of $ A_{\widetilde{K}} $. By assumption, we have that both $\mathcal{A}'$ and $ \widetilde{\mathcal{A}} $ are semi-abelian. Since $ \widetilde{\mathcal{A}}_{S'} $ is smooth, and $\mathcal{A}'$ has the N\'eronian property, we get a canonical morphism $ \widetilde{\mathcal{A}}_{S'} \rightarrow \mathcal{A}' $, lifting the identity map on the generic fibers. Since $ \widetilde{\mathcal{A}}_{S'} $ is semi-abelian, it follows from Proposition 7.4/3 in \cite{Ner} that this morphism induces an isomorphism $ (\widetilde{\mathcal{A}}_k)^0 \cong (\mathcal{A}'_k)^0 $. In particular, we get that $ T_{\widetilde{\mathcal{A}}_k,0} = T_{\mathcal{A}'_k, 0} $.

Consider now the filtration $ \{ F^i_m \widetilde{\mathcal{A}}_k \} $ of $ \widetilde{\mathcal{A}}_k $ induced by the extension $ S'/\widetilde{S} $. Since $ \widetilde{\mathcal{A}} $ is semi-abelian, we have that $ F^i_m \widetilde{\mathcal{A}}_k = 0 $ for all $ i > 0 $. Therefore, we get that
$$ \widetilde{\mathcal{A}}_k = F^0_m \widetilde{\mathcal{A}}_k = Gr^0_m \widetilde{\mathcal{A}}_k = (\mathcal{A}'_k)^{\boldsymbol{\mu}_m}. $$
But now 
$$ (T_{\mathcal{A}'_k, 0})^{\boldsymbol{\mu}_m} = T_{(\mathcal{A}'_k)^{\boldsymbol{\mu}_m}, 0} = T_{\mathcal{A}'_k, 0}, $$ 
and so it follows that $ \boldsymbol{\mu}_m $ acts trivially on $ T_{\mathcal{A}'_k, 0} $. 

Let us now consider the filtration $ \{ F^i_n \mathcal{A}_k \} $ induced by the extension $ S'/S $. The jumps in this filtration are determined by the $ \boldsymbol{\mu}_n $-action on $ T_{\mathcal{A}'_k, 0} $. Assume that $ T_{\mathcal{A}'_k, 0}[i] \neq 0 $, for some $ i \in \{0, \ldots, n -1 \} $. On this subspace, every $ \xi \in \boldsymbol{\mu}_n $ acts by multiplication by $ \xi^i $. We can identify $ \boldsymbol{\mu}_m $ with the $\tilde{n}$-th powers in $ \boldsymbol{\mu}_n $, and since we established above that $ \boldsymbol{\mu}_m $ acts trivially, it follows that $ \xi^{\tilde{n} i} = 1 $. So therefore $ \tilde{n} i = k n $ for some $ k \in \{0, \ldots, \tilde{n} - 1 \} $, and we get that $ i = k n/\tilde{n} $.
\end{pf}

\vspace{0.5cm}

We can now formulate the following result:

\begin{prop}\label{tamejumpprop}
If $ A/K $ obtains semi-abelian reduction over a tamely ramified extension of $K$, then the jumps in the filtration $ \{ \mathcal{F}^a \mathcal{A}_k \} $ occur at indices $ k/\tilde{n} $, where $ k \in \{0, \ldots, \tilde{n} - 1 \} $, and where $ \tilde{n} $ is the degree of the minimal extension $ \widetilde{K}/K $ that realizes semi-abelian reduction for $ A $.
\end{prop}
\begin{pf}
Let us consider the sequence of integers $ (\tilde{n} m)_m $, where $m$ runs over the positive integers that are not divisibe by $p$. For the extension of degree $n = \tilde{n} m$, Lemma \ref{tame jumps} gives that the jumps of $ \{ F^i_n \mathcal{A}_k \} $ occur at indices $ i = k n / \tilde{n} $, where $ 0 \leq k \leq \tilde{n} - 1$. It follows that the jumps of $ \{ \mathcal{F}^a \mathcal{A}_k \} $ will be among the limits of the expressions $ i/n = k/ \tilde{n} $, as $m$ goes to infinity, and the result follows.
\end{pf}

\subsection{The case of Jacobians}\label{jacobiancase}
Let $X/K$ be a smooth, projective and geometrically connected curve of genus $ g > 0 $. We shall also make the assumption that $ X(K) \neq \emptyset $. If $K'/K $ is a tame extension, it follows that $ X_{K'} $ hase the same properties. Let $ J' = J_{K'} $ denote the Jacobian of $X_{K'}$, and let $ \mathcal{J}'/S' $ be the N\'eron model of $J'$ over $S'$. 

We can let $G$ act on $ X_{K'} $ via the action on the second factor. Let $ \mathcal{Y}/S' $ be a regular model of $ X_{K'} $ such that the $G$-action on $ X_{K'} $ extends to $ \mathcal{Y} $. According to \cite{Ner}, Theorem 9.5/4, there is a canonical isomorphism
$$ \Pic^0_{\mathcal{Y}/S'} \cong \mathcal{J}'^0, $$
where $ \mathcal{J}'^0 $ is the identity component of $ \mathcal{J}' $, and where $ \Pic^0_{\mathcal{Y}/S'} $ is the identity component of the relative Picard functor $ \Pic_{\mathcal{Y}/S'} $. Hence, on the special fibers, we get an isomorphism
$$ \Pic^0_{\mathcal{Y}_k/k} \cong \mathcal{J}'^0_k. $$
By \cite{Ner}, Theorem 8.4/1, it follows that we can canonically identify
\begin{equation}\label{H^1=T}
H^1(\mathcal{Y}_k, \mathcal{O}_{\mathcal{Y}_k}) \cong T_{\mathcal{J}'_k,0}. 
\end{equation}

We are interested in computing the irreducible characters for the representation of $\boldsymbol{\mu}_n$ on $T_{\mathcal{J}'_k,0}$. With the identification in \ref{H^1=T} above, we see that this can be done by computing the irreducible characters for the representation of $\boldsymbol{\mu}_n$ on $ H^1(\mathcal{Y}_k, \mathcal{O}_{\mathcal{Y}_k}) $.

By combining the discussion in this section with properties of the representation of $\boldsymbol{\mu}_n$ on $ H^1(\mathcal{Y}_k, \mathcal{O}_{\mathcal{Y}_k}) $, we obtain in Corollary \ref{specific jump corollary} a quite precise description of the jumps of the filtration $ \{ \mathcal{F}^a \mathcal{J}_k \} $.

\section{Tame extensions and Galois actions}\label{extensions and actions}

\subsection{Some definitions.}
Throughout this paper, $R$ will denote a complete discrete valuation ring, with fraction field $K$, and algebraically closed residue field $k$. 

$X/K$ will be a smooth, projective, geometrically connected curve over $K$, of genus $g(X)>0$. For simplicity, we will say that $X$ is a \emph{curve} of genus $g$ over $K$.
\begin{dfn}
A scheme $ \mathcal{X} $ is called a \emph{model} of $X$ over $ S = \Spec(R) $ if $ \mathcal{X} $ is integral and normal, projective and flat over $S$, and with generic fibre $ \mathcal{X}_K \cong X $.
\end{dfn}

It is easily seen that one can always find a model for $X$ (cf. ~\cite{Liubook}, Proposition 10.1.8). Let us now mention a few types of models that will frequently occur in this paper.

Since a model $ \mathcal{X} $ is normal, with smooth generic fiber, it follows that the singular locus consists of a finite set of closed points in the special fiber. Furthermore, there exists a \emph{strong} desingularization $ \phi : \mathcal{Z} \rightarrow \mathcal{X} $, i.e., $ \phi $ is an isomorphism over the regular locus of $ \mathcal{X} $ (\cite{Liubook}, Corollary 8.3.51). In fact, there even exists a \emph{minimal} (strong) desingularization $ \rho : \mathcal{Y} \rightarrow \mathcal{X} $, characterized by the property that any other desingularization of $ \mathcal{X} $ factors via $ \rho $ (\cite{Liubook}, Proposition 9.3.32).

It follows that we can always find a regular model $ \mathcal{X}/S $ for $X/K$. By blowing up points in the special fiber, we can even ensure that the irreducible components of the special fiber are smooth, and intersect transversally. Such a model will be called a \emph{strict normal crossings} model for $X/K$, or for short, an SNC-model.
\subsection{Construction}\label{2.2}
Let $X/K$ be a curve of genus $g > 0 $, and let $ \mathcal{X}/S $ be an SNC-model for $X/K$. 
Let $ K \subset K' $ be a finite, separable field extension, and let $R'$ be the integral closure of $R$ in $K'$. Since $R$ is complete, we have that $R'$ is a complete discrete valuation ring (\cite{Serre}, Proposition II.3). Making the finite base extension $S' = \Spec(R') \rightarrow S = \Spec(R) $,  we obtain a commutative diagram
$$ \xymatrix{
 \mathcal{Y} \ar[d] \ar[r]^{\rho} & \mathcal{X}' \ar[d] \ar[r]^{f}  & \mathcal{X} \ar[d]\\ 
  S' \ar[r]^{\Id} & S' \ar[r] & S, } $$
where $ \mathcal{X}' $ is the normalization of the pullback $ \mathcal{X}_{S'} = \mathcal{X} \times_S S' $ ($ \mathcal{X}_{S'} $ is integral by Lemma \ref{lemma construction} below), and $ \rho : \mathcal{Y} \rightarrow \mathcal{X}' $ is the minimal desingularization. The map $ f : \mathcal{X}' \rightarrow \mathcal{X} $ is the composition of the projection $ \mathcal{X}_{S'} \rightarrow \mathcal{X} $ with the normalization $ \mathcal{X}' \rightarrow \mathcal{X}_{S'} $.

\begin{lemma}\label{lemma construction} With the hypotheses above, the following statements hold:
\begin{enumerate}
\item The pullback $ \mathcal{X}_{S'} $ is integral.
\item $ f : \mathcal{X}' \rightarrow \mathcal{X} $ is a finite morphism.
\end{enumerate}
\end{lemma}
\begin{pf}
(i) Let us first note that the generic fiber of $ \mathcal{X}_{S'} $ is the pullback $\mathcal{X}_K \otimes_K K' $, where $ \mathcal{X}_K $ is the generic fiber of $ \mathcal{X} $. By assumption $ \mathcal{X}_K $ is smooth and geometrically connected over $K$, hence the generic fiber of $ \mathcal{X}_{S'} $ is in particular integral. Now, since $ \mathcal{X}_{S'} \rightarrow S' $ is flat, it follows from \cite{Liubook}, Proposition 4.3.8, that $ \mathcal{X}_{S'} $ is integral as well.

(ii) Since $ R' $ is a \emph{complete} discrete valuation ring it is excellent (\cite{Liubook}, Theorem 8.2.39). As $ \mathcal{X}_{S'} $ is of finite type over $S'$, it follows that $ \mathcal{X}_{S'} $ is an excellent scheme, and hence the normalization morphism $ \mathcal{X}' \rightarrow \mathcal{X}_{S'} $ is finite (\cite{Liubook}, Theorem 8.2.39). The projection $ \mathcal{X}_{S'} \rightarrow \mathcal{X} $ is finite, since it is the pullback of the finite morphism $ S' \rightarrow S $. So the composition $f$ of these two morphisms is indeed finite.
\end{pf} 

\subsection{Galois actions} Let us now assume that the field extension $ K \subset K' $ is Galois with group $G$. Every $ \sigma \in G $ induces an automorphism of $R'$ that fixes $R$, and we have furthermore that $ R'^{G} = R $. So there is an injective group homomorphism $ G \rightarrow \Aut(S') $, and we may view $ S' \rightarrow S $ as the quotient map. 


We can lift the $G$-action to $ \mathcal{X}_{S'} $, via the action on the second factor. So there is a group homomorphism $ G \rightarrow \Aut(\mathcal{X}_{S'}) $. For any element $ \sigma \in G $, we shall still denote the image in $ \Aut(\mathcal{X}_{S'}) $ by $ \sigma $. Proposition \ref{prop. 2.3} below shows that this action can be lifted uniquely both to the normalization $ \mathcal{X}' $ and to the minimal desingularization $ \mathcal{Y} $ of $ \mathcal{X}' $.

Recall the universal property of the normalization, saying that if $ g : \mathcal{Z} \rightarrow \mathcal{X}_{S'} $ is a \emph{dominant} morphism where $ \mathcal{Z} $ is a normal scheme, then $g$ factors uniquely via $ \mathcal{X}' \rightarrow \mathcal{X}_{S'} $. 
 
\begin{prop}\label{prop. 2.3} With the hypotheses above, the following statements hold:
\begin{enumerate}
\item The $G$-action on $ \mathcal{X}_{S'} $ lifts uniquely to the normalization $ \mathcal{X}' $.

\item The $ G $-action on $ \mathcal{X}' $ lifts uniquely to the minimal desingularization $ \mathcal{Y} $.

\item For any $ \sigma \in G $, let $ \sigma $ denote the induced automorphism of $ \mathcal{X}' $, and let $ \tau $ be the unique lift of $ \sigma $ to $ \Aut(\mathcal{Y}) $. Then we have that $ \tau(\rho^{-1}(\Sing(\mathcal{X}'))) = \rho^{-1}(\Sing(\mathcal{X}')) $. That is, the exceptional locus is mapped into itself under the $G$-action on $ \mathcal{Y} $.
\end{enumerate}
\end{prop}
\begin{pf}
(i) The proofs of the liftings of the $G$-actions to $ \mathcal{X}' $ and $ \mathcal{Y} $ are similar, so we do not write out the details for $ \mathcal{X}' $.

(ii) Let $ \sigma \in G $. Then $ \sigma $ acts as an automorphism $ \sigma : \mathcal{X}' \rightarrow \mathcal{X}' $. Consider the diagram
$$ \xymatrix{
 \mathcal{Y} \ar[d]_{\rho}  & \mathcal{Y} \ar[d]^{\rho} \\ 
   \mathcal{X}' \ar[r]^{\sigma} &  \mathcal{X}'. } $$

The composition $ \sigma \circ \rho : \mathcal{Y} \rightarrow \mathcal{X}' $ is a desingularization, so by the minimality of $ \rho : \mathcal{Y} \rightarrow \mathcal{X}' $, there exists a unique morphism $ \tau : \mathcal{Y} \rightarrow \mathcal{Y} $ such that the diagram
$$ \xymatrix{
 \mathcal{Y} \ar[d]_{\rho}  \ar[r]^{\tau} & \mathcal{Y} \ar[d]^{\rho} \\ 
   \mathcal{X}' \ar[r]^{\sigma} &  \mathcal{X}' } $$
commutes. Let us first show that $ \tau $ is an automorphism. Let $ \sigma' = \sigma^{-1} $. Since $ \sigma $ is an automorphism, also $ \sigma' \circ \rho $ is a desingularization, hence there exists a unique $ \tau' : \mathcal{Y} \rightarrow \mathcal{Y} $ such that $ \rho \circ \tau' = \sigma' \circ \rho $. 

It now follows that
$$ \rho \circ \tau' \circ \tau = \sigma' \circ \rho \circ \tau = \sigma' \circ \sigma \circ \rho = \rho. $$
Hence the diagram
$$ \xymatrix{
 \mathcal{Y} \ar[dr]_{\rho}  \ar[rr]^{\tau' \circ \tau} & & \mathcal{Y} \ar[dl]^{\rho} \\ 
    & \mathcal{X}' & } $$
commutes. But also $ \Id_{\mathcal{Y}} $ makes this diagram commute, and by the universal property of $ \rho $, we get that $ \tau' \circ \tau = \Id_{\mathcal{Y}} $. By symmetry we also have that $ \tau \circ \tau' = \Id_{\mathcal{Y}} $, hence $ \tau' = \tau^{-1} $, so in particular $ \tau $ is an automorphism.

This shows the existence of the map $ G \rightarrow \Aut(\mathcal{Y}) $. The proof that this map is indeed a group homomorphism is straightforward, and is omitted here. 

(iii) The desingularization $ \rho : \mathcal{Y} \rightarrow \mathcal{X}' $ is \emph{strong}, i.e., it is an isomorphism over $ \Reg(\mathcal{X}') $. Since $ \sigma $ is an automorphism, we have that $ \sigma(\Sing(\mathcal{X}')) = \Sing(\mathcal{X}') $. But then $ \tau(\rho^{-1}(\Reg(\mathcal{X}'))) = \rho^{-1}(\Reg(\mathcal{X}')) $, so the claim follows immediately.
\end{pf}

\subsection{Tame extensions}\label{assumption on degree}
For the rest of this paper, we shall always make the assumption that the degree $ n = [K':K] $ is not divisible by the residue characteristic $p$. In other words, $ S' \rightarrow S $ is a tamely ramified extension. Since $k$ is algebraically closed, $k$ has a full set $ \boldsymbol{\mu}_n $ of $n$-th roots of unity. Furthermore, since $R$ is complete, and in particular henselian, we may lift all $n$-th roots of unity to $R$. We can choose a uniformizing parameter $ \pi \in R $ such that $ K' = K[\pi']/(\pi'^n - \pi) $. The extension $ K \subset K' $ is Galois, with group $ G = \boldsymbol{\mu}_n $. Then $ R' := R[\pi']/(\pi'^n - \pi) $ is the integral closure of $R$ in $K'$, and $ \pi'$ is a uniformizing parameter for $R'$.

\subsection{Assumptions on $\mathcal{X}$}\label{assumption on surface}
In the rest of this paper, we shall make two assumptions in the situation in Section \ref{2.2}:

\begin{ass}\label{ass. 2.4}
Let $ x \in \mathcal{X} $ be a closed point in the special fiber such that two irreducible components $C_1$ and $C_2$ of $ \mathcal{X}_k $ meet at $x$, and let $ m_i = \Mult(C_i) $. We will always assume that \emph{at least} one of the $m_i$ is not divisible by $p$. 
\end{ass}

With this assumption, we can find an isomorphism 
$$ \widehat{\mathcal{O}}_{\mathcal{X},x} \cong R[[u_1,u_2]]/(\pi - u_1^{m_1} u_2^{m_2}) $$
(cf. \cite{CED}, proof of Lemma 2.3.2).

\begin{ass}\label{ass. 2.5}
Let $l$ be the least common multiple of the multiplicities of the irreducible components of $\mathcal{X}_k$. Then we assume that $ \Gcd(l,n) = 1 $. 
\end{ass}

Let us now recall a few facts developed in \cite{Tame}:

$\bullet$ Let $ x \in \mathcal{X} $ be a closed point in the special fiber. Because of Assumption \ref{ass. 2.5}, there is a unique point $ x' \in \mathcal{X}'_k $ that maps to $x$. The local analytic structure of $ \mathcal{X}' $ at $x'$ depends only on $n$ and the local analytic structure of $ \mathcal{X} $ at $x$. If $x$ belongs to a unique component of $ \mathcal{X}_k $, then $x'$ belongs to a unique component of $ \mathcal{X}'_k $, and $ \mathcal{X}' $ is regular at $x'$. If $x$ is an intersection point of two distinct components, then the same is true for $x'$, and $ \mathcal{X}' $ will have a \emph{tame cyclic quotient singularity} at $x'$. 

$\bullet$ The minimal desingularization $ \mathcal{Y} $ of $ \mathcal{X}' $ is an SNC-model. Furthermore, the structure of $ \mathcal{Y} $ locally above a tame cyclic quotient singularity $x' \in \mathcal{X}'$ is completely determined by the structure locally at $ x = f(x') $ and the degree $n$ of the extension. The inverse image of $x'$ consists of a chain of smooth and rational curves whose multiplicities and self intersection numbers may be computed from the integers $n, m_1$ and $m_2$.

$\bullet$ For every irreducible component $C$ of $ \mathcal{X}_k $, there is precisely one component $C'$ of $ \mathcal{X}'_k $ that dominates $C$. The component $C'$ is isomorphic to $C$, and we have that $ \Mult_{\mathcal{X}'_k}(C') = \Mult_{\mathcal{X}_k}(C) $. It follows that the combinatorial structure of $ \mathcal{Y}_k $ is completely determined by the combinatorial structure of $ \mathcal{X}_k $ and the degree of $S'/S$.

\subsection{}
We will now begin to describe the $G$-action on $ \mathcal{X}' $ and $ \mathcal{Y} $ in more detail. Assumptions \ref{ass. 2.4} and \ref{ass. 2.5} will impose some restrictions on this action.


    
\begin{prop}\label{action on desing}
Let $ \rho : \mathcal{Y} \rightarrow \mathcal{X'} $ be the minimal desingularization. Then the following properties hold: 
\begin{enumerate}
\item Let $D$ be an irreducible component of $ \mathcal{Y}_k $ that dominates a component of $ \mathcal{X}_k $. Then $D$ is stable under the $G$-action, and $G$ acts \emph{trivially} on $D$.
\item Let $ x' \in \mathcal{X}' $ be a singular point, and let $E_1, \ldots, E_{l}$ be the exceptional components mapping to $x'$ under $ \rho $. Then every $E_i$ is stable under the $G$-action, and every node in the chain $ \rho^{-1}(x') $ is fixed under the $G$-action.
\end{enumerate}
\end{prop}
\begin{pf}
Let us first note that the map $ \mathcal{X}_{S'} \rightarrow \mathcal{X} $ is an isomorphism on the special fibers. Moreover, the action on the special fiber of $ \mathcal{X}_{S'} $ is easily seen to be trivial, so every closed point in the special fiber is fixed. Since the action on $ \mathcal{X}' $ commutes with the action on $ \mathcal{X}_{S'} $, it follows that every point in the special fiber of $ \mathcal{X}' $ is fixed. In particular, every irreducible component $C'$ of $ \mathcal{X}'_k $ is stable under the $G$-action, and the restriction of this action to $C'$ is trivial. Since the action on $ \mathcal{Y} $ commutes with the action on $ \mathcal{X}' $, it follows that the same is true for the strict transform $D$ of $C'$ in $ \mathcal{Y} $. This proves (i).

For (ii), we observe that since $x'$ is fixed, we have that $\rho^{-1}(x')$ is stable under the $G$-action. But also the two branches meeting at $x'$ are fixed. Let $D$ be the strict transform of any of these two branches. From part (i), it follows that the point where it meets the exceptional chain $\rho^{-1}(x')$ must be fixed. So if $E_1$ is the component in the chain meeting $\widetilde{D}$, then $E_1$ must be mapped into itself. Let $E_2$ be the next component in the chain. Then the point where $E_1$ and $E_2$ meet must also be fixed, so $E_2$ must also be mapped to itself. Continuing in this way, it is easy to see that all of the exceptional components are stable under the $G$-action, and that all nodes in $\rho^{-1}(x')$ are fixed points.  
\end{pf}

\begin{cor}\label{g^{-1}(Z) = Z}
Let $ 0 \leq Z \leq \mathcal{Y}_k $ be an effective divisor. Then we have that the $G$-action restricts to $Z$.
\end{cor}
\begin{pf}
 Since $Z$ is an effective Weil divisor, we can write $ Z = \sum_{C} r_C C $, where $C$ runs over the irreducible components of $ \mathcal{Y}_k $, and $r_C$ is a non-negative integer for all $C$. But Proposition \ref{action on desing} states that all irreducible components $C$ of $ \mathcal{Y}_k $ are stable under the $G$-action, and hence we get that the same holds for $ Z $. In other words, the action restricts to $Z$.
\end{pf}

\vspace{0.3cm}

From Proposition \ref{action on desing} above, it follows that every node $y$ in $ \mathcal{Y}_k $ is a fixed point for the $G$-action on $ \mathcal{Y} $. Hence there is an induced action on $ \mathcal{O}_{\mathcal{Y},y} $ and on the cotangent space $ m_y/m_y^2 $, where $ m_y \subset \mathcal{O}_{\mathcal{Y},y} $ is the maximal ideal. In order to get a precise description of the action on the cotangent space, we will first describe the action on the completion $ \widehat{\mathcal{O}}_{\mathcal{Y},y} $.
 
Since, by Proposition \ref{action on desing}, every irreducible component $D$ of $ \mathcal{Y}_k $ is mapped to itself under the $G$-action, it follows that the $G$-action restricts to $D$ and that the points where $D$ meets the rest of the special fiber are fixed. In the case where $G$ acts non-trivially on $D$, we will see in Proposition \ref{prop. 3.3} that the fixed points for the $G$-action on $D$ are precisely the points where $D$ meets the rest of the special fiber. In particular, we wish to describe the action on $D$ locally at the fixed points. 

\section{Desingularizations and actions}\label{lifting the action}
In this section, we study how one can explicitly describe the action on the minimal desingularization $ \rho : \mathcal{Y} \rightarrow \mathcal{X}' $. Since we are only interested in this action locally at fixed points or stable components in the exceptional locus of $\rho$, we will begin with showing that we can reduce to studying the minimal desingularization locally at a singular point $ x' \in \mathcal{X}' $. This is an important step, since we have a good description of the complete local ring $ \widehat{\mathcal{O}}_{\mathcal{X}',x'} $. In particular, we can find a nice algebraization of this ring, with a compatible $G$-action. It turns out that it suffices for our purposes to study the minimal desingularization of this ring, and the lifted $G$-action.

In the second part of this section, we study the desingularization of the algebraization of $ \widehat{\mathcal{O}}_{\mathcal{X}',x'} $. We will use the explicit blow up procedure given in \cite{CED} for the resolution of tame cyclic quotient singularities, which will allow us to describe precisely how the $G$-action lifts. In particular, we describe the action on the completion of the local rings at the nodes in the exceptional locus, and the action on the exceptional components. These results are gathered in Proposition \ref{prop. 3.3}.    

\subsection{Reduction to the complete local rings}
If $ x' \in \mathcal{X}' $ is a singular point, we need to understand how $G$ acts on $ \widehat{\mathcal{O}}_{\mathcal{X}',x'} $. In order to do this, we consider $ f(x') = x \in \mathcal{X} $, where $ f : \mathcal{X}' \rightarrow \mathcal{X} $. Then $x$ is a closed point in the special fiber, and we have that   
$$ \widehat{\mathcal{O}}_{\mathcal{X},x} \cong R[[v_1,v_2]]/(\pi - v_1^{m_1} v_2^{m_2}), $$
where $m_1$ and $m_2$ are positive integers. Let $n$ be the degree of $ R'/R $. By assumption, we have that $ n $ is relatively prime to $m_1$ and $m_2$. Before we proceed, we would like to remark that we will in the discussion that follows use some properties that were proved in Section 2 of \cite{Tame}.

Consider now the pullback $ \mathcal{X}_{S'} $. We let $ G = \boldsymbol{\mu}_n $ act on $ \mathcal{X}_{S'} $ via its action on the second factor. We point out that we here choose the action given by $ [\xi](\pi') = \xi \pi' $ for any $ \xi \in \boldsymbol{\mu}_n $. Choosing this action is notationally convenient when we work with rings. However, the natural right action on $ \mathcal{X}_{S'} $ is the inverse to the one we use here. In particular, the irreducible characters for the representation of $ \boldsymbol{\mu}_n $ on $ H^1(\mathcal{Y}_k, \mathcal{O}_{\mathcal{Y}_k}) $ induced by the action chosen here on $\mathcal{X}_{S'}$ will be the inverse characters to those induced by the right action.

By abuse of notation, let $ x \in \mathcal{X}_{S'} $ be the point mapping to $ x \in \mathcal{X} $. Then we have that the map $ \mathcal{O}_{\mathcal{X},x} \rightarrow \mathcal{O}_{\mathcal{X}_{S'},x} $ induced by the projection $ \mathcal{X}_{S'} \rightarrow \mathcal{X} $ can be described by the tensorization
$$ \mathcal{O}_{\mathcal{X},x} \rightarrow \mathcal{O}_{\mathcal{X},x} \otimes_R R', $$
and that the $G$-action on $ \mathcal{O}_{\mathcal{X}_{S'},x} = \mathcal{O}_{\mathcal{X},x} \otimes_R R' $ is given by the action on $R'$.

Since $ \mathcal{O}_{\mathcal{X},x} \rightarrow \mathcal{O}_{\mathcal{X},x} \otimes_R R' $ is finite, completion commutes with tensoring with $R'$, so we get that
$$ \widehat{\mathcal{O}}_{\mathcal{X}_{S'},x} = \widehat{\mathcal{O}}_{\mathcal{X},x} \otimes_R R', $$
and hence the $G$-action on $ \widehat{\mathcal{O}}_{\mathcal{X}_{S'},x} $ is induced from the action on $R'$ in the second factor. It follows that
$$ \widehat{\mathcal{O}}_{\mathcal{X}_{S'},x} \cong R'[[v_1,v_2]]/(\pi'^n - v_1^{m_1} v_2^{m_2}), $$
and that the $G$-action is given by $ [\xi](\pi') = \xi \pi' $ and $ [\xi](v_i) = v_i $, for any $ \xi \in \boldsymbol{\mu}_n $. 

Let $ \mathcal{X}' \rightarrow \mathcal{X}_{S'} $ be the normalization. Then $x'$ is the unique point mapping to $x$, and the induced map $ \widehat{\mathcal{O}}_{\mathcal{X}_{S'},x} \rightarrow \widehat{\mathcal{O}}_{\mathcal{X}',x'} $ is the normalization of $ \widehat{\mathcal{O}}_{\mathcal{X}_{S'},x} $. Furthermore, it follows also that the $G$-action on $ \widehat{\mathcal{O}}_{\mathcal{X}',x'} $ induced by the action on $ \mathcal{X}' $ is the unique lifting of the $G$-action on $ \widehat{\mathcal{O}}_{\mathcal{X}_{S'},x} $ to the normalization $ \widehat{\mathcal{O}}_{\mathcal{X}',x'} $.

Let $ \rho : \mathcal{Y} \rightarrow \mathcal{X}' $ be the minimal desingularization, and consider the fiber diagram
$$ \xymatrix{
\widehat{\mathcal{Y}} \ar[d]_{\hat{\rho}} \ar[r]^f & \mathcal{Y} \ar[d]^{\rho} \\
\Spec(\widehat{\mathcal{O}}_{\mathcal{X}',x'}) \ar[r] & \mathcal{X}' .}
$$
Then $ \hat{\rho} $ is the minimal desingularization of $ \Spec(\widehat{\mathcal{O}}_{\mathcal{X}',x'}) $ (cf. ~\cite{Lip}, Lemma 16.1, and use the fact that $ \mathcal{Y} $ is minimal), and hence the $G$-action on $ \Spec(\widehat{\mathcal{O}}_{\mathcal{X}',x'}) $ lifts uniquely to $ \widehat{\mathcal{Y}} $. 

We have that $f$ induces an isomorphism of the exceptional loci $ \hat{\rho}^{-1}(x') $ and $ \rho^{-1}(x') $. Let $E$ be an exceptional component. Then the $G$-action restricts to $E$, and it is easily seen that $f$, when restricted to $E$, is equivariant. 

Furthermore, for any closed point $ y \in \rho^{-1}(x') $, we have that $f$ induces an isomorphism $ \widehat{\mathcal{O}}_{\mathcal{Y},y} \cong \widehat{\mathcal{O}}_{\widehat{\mathcal{Y}},y} $ (one can argue in a similar way as in the proof of \cite{Liubook}, Lemma 8.3.49). If $y$ is a fixed point, it is easily seen that this isomorphism is equivariant.
We therefore conclude that in order to describe the action on $\mathcal{Y}$ locally at the exceptional locus over $x'$, it suffices to consider the minimal desingularization of $ \Spec(\widehat{\mathcal{O}}_{\mathcal{X}',x'}) $.

\subsection{Equivariant algebraization}
In order to find an algebraization of $ \widehat{\mathcal{O}}_{\mathcal{X}',x'} $, we consider first the polynomial ring $ V = R'[v_1,v_2]/(\pi'^n - v_1^{m_1} v_2^{m_2}) $. We let $G$ act on $ V $ by $ [ \xi ](\pi') = \xi \pi' $ and $ [ \xi ](v_i) = v_i $ for $ i = 1,2 $, for any $ \xi \in G $. Note that the maximal ideal $ P = (\pi', v_1,v_2) $ is fixed, and hence there is an induced action on the completion $ \widehat{V}_P = R'[[v_1,v_2]]/(\pi'^n - v_1^{m_1} v_2^{m_2}) $, given as above. This gives a $G$-equivariant algebraization of $ \widehat{\mathcal{O}}_{\mathcal{X}_{S'},x} $. 

Let us consider now the $R'$-homomorphism
$$ V = R'[v_1,v_2]/(\pi'^n - v_1^{m_1} v_2^{m_2}) \rightarrow T = R'[t_1,t_2]/(\pi' - t_1^{m_1} t_2^{m_2}), $$
given by $ v_i \mapsto t_i^n $. We let $ \boldsymbol{\mu}_n $ act on $T$, relatively to $R'$, by $ [ \eta ](t_1) = \eta t_1 $, $ [ \eta ](t_2) = \eta^r t_2 $, where $ r $ is the unique integer $ 0 < r < n $ such that $ m_1 + r m_2 \equiv_n 0 $. Note that this is an ad hoc action introduced to compute the normalization, which must not be confused with the natural $G$-action. Arguing as in Section 2 of \cite{Tame}, we find that the induced map $ V \rightarrow U := T^{\boldsymbol{\mu}_n} $ is the normalization of $V$. Furthermore, it is easily seen that there is a unique maximal ideal $Q \subset U $ mapping to $P$, corresponding to the origin $ (\pi, t_1,t_2) $ in $T$.

We will need the following lemma:

\begin{lemma}\label{lemma 4.1}
Consider the fiber diagrams
$$ \xymatrix{
V \ar[d] \ar[r] &  V_P \ar[d] \ar[r] & \widehat{V}_P \ar[d]\\ 
U \ar[r] & U_P \ar[r] & \widehat{U}_P, } $$
where $ U_P := U \otimes_V V_P $, $ \widehat{U}_P := U_P \otimes_{V_P} \widehat{V}_P $, and where $ V \rightarrow U $ is the normalization. 

Then we have that $ U_P = U_Q $, where $U_Q $ is the localization in the maximal ideal $ Q \subset U $. Furthermore, the maps $ V_P \rightarrow U_P $ and $ \widehat{V}_P \rightarrow  \widehat{U}_P $ are the normalizations.
\end{lemma}
\begin{pf}
Since $ V \rightarrow U $ is finite, it follows that $ V_P \rightarrow U_P $ is finite, and coincides with the normalization of $V_P$, since normalization commutes with localization. Furthermore, $ U_P $ is semi-local, with maximal ideals corresponding exactly to the maximal ideals of $U$ restricting to $P$. It follows that $ U_P = U_Q $, and that the map $ U \rightarrow U_Q $ above is the localization map.

Since $ V_P \rightarrow U_P $ is finite, we can identify the map $ U_P \rightarrow \widehat{U}_P $ with the completion of $ U_P $ in the radical ideal. We have that $V_P$ is reduced and excellent, since it is a localization of $V$. But then normalization commutes with completion, that is, $ \widehat{U}_Q = \widehat{U}_P = (\widehat{V}_P)' $, and the map $ \widehat{V}_P \rightarrow \widehat{U}_P = \widehat{U}_Q $ is the normalization (\cite{Liubook}, Proposition 8.2.41).
\end{pf}

\vspace{0.5cm}

It follows from this lemma that $U$ is an equivariant algebraization of $ \widehat{\mathcal{O}}_{\mathcal{X}',x'} $. Let $ \rho_U : \mathcal{Z} \rightarrow \Spec(U) $ be the minimal desingularization. Then it follows that we have a fiber diagram
$$ 
\xymatrix{
\widehat{\mathcal{Y}} \ar[d]_{\hat{\rho}} \ar[r] & \mathcal{Z} \ar[d]^{\rho_U} \\ 
\Spec(\widehat{\mathcal{O}}_{\mathcal{X}',x'}) \ar[r] & \Spec(U),}
$$
where all maps commute with the various $G$-actions. In particular, we can now conclude that to describe the $G$-action on $ \widehat{\mathcal{Y}} $ locally at fixed points or components in the exceptional locus, it suffices to compute the corresponding data for $ \mathcal{Z} $.



\subsection{Lifting the action to the normalization}


We will continue to work with the rings $V$, $U$ and $T$. We saw above that there was a $G$-action on $V$. Since $ U $ is the normalization of $V$, this action will lift uniquely, and we want to give a precise description of that action. Let us first remark that by arguing exactly as in Section 2 in \cite{Tame}, it follows that $ U = T^{\boldsymbol{\mu}_n} $ is generated as a $V$-module by the monomials $t_1^i t_2^j $ such that $ 0 \leq i,j < n $, where $ i + r j \equiv 0 $ modulo $n$.  

In order to describe the $G$-action on $U = T^{\boldsymbol{\mu}_n}$, it suffices to give the action on each of the $V$-module generators. This is done in the following lemma:

\begin{lemma}\label{lemma 3.1}
Let $ t_1^i t_2^j \in T $ be such that $ i + rj \equiv 0 $ modulo $n$. Then we have, for any $ \xi \in G $, that $ [ \xi ] (t_1^i t_2^j) = \xi^{i \alpha_1} t_1^i t_2^j $, where $ \alpha_1 $ is an inverse to $m_1$ modulo $n$.
\end{lemma}
\begin{pf}
In $T$ we have the relations $ t_1^{m_1} t_2^{m_2} = \pi' $, $ t_1^{n} = v_1 $ and $ t_2^{n} = v_2 $, hence all these elements lie in the image of $V$, and we therefore know how $ G $ acts on them. We have that $[ \xi ] (t_1^{m_1} t_2^{m_2}) = \xi t_1^{m_1} t_2^{m_2} $ and $ [ \xi ] (t_i^{n}) = t_i^{n} $, for any $ \xi \in G $.

For any $i,j$ as above, let $k(i,j)$ be the integer satisfying $ i + r j = k(i,j) n $. Furthermore, define $ K(i,j) := m_2 k(i,j) - k(m_1,m_2) j $. One computes easily that $ i m_2 = K(i,j) n + j m_1 $. From this, we get the relation
$$ (t_1^{m_1} t_2^{m_2})^i = (t_1^i t_2^j)^{m_1} (t_2^{n})^{K(i,j)} $$
in $\Frac(V)$, and hence
$$ (t_1^i t_2^j)^{m_1} = (t_1^{m_1} t_2^{m_2})^i (t_2^{n})^{ - K(i,j)}. $$
As $ \Gcd(m_1,n) = 1 $, we can find integers $ \alpha_1, \beta_1 $ such that $ \alpha_1 m_1 + \beta_1 n = 1 $. So we get
$$ t_1^i t_2^j = (t_1^i t_2^j)^{\alpha_1 m_1} (t_1^i t_2^j)^{\beta_1 n} = ((t_1^i t_2^j)^{m_1})^{\alpha_1} (t_1^{n})^{i \beta_1} (t_2^{n})^{j \beta_1} = $$
$$ (t_1^{m_1} t_2^{m_2})^{i \alpha_1} (t_2^{n})^{ - \alpha_1 K(i,j)} (t_1^{n})^{i \beta_1} (t_2^{n})^{j \beta_1}. $$
From this relation, it follows that $ [ \xi ] (t_1^i t_2^j) = \xi^{i \alpha_1} t_1^i t_2^j $. In particular, we note that every $ \boldsymbol{\mu}_n $-invariant monomial $ t_1^i t_2^j $ in $T$ is an ``eigenvector'' for the $G$-action, and that we can compute explicitly the eigenvalue. 
\end{pf}

\subsection{The action on the minimal desingularization}\label{section 3.4}
Now that we know how $G$ acts on $U$, we will use this description to lift the $G$-action to the minimal desingularization of $ \Spec(U) $. For this, we shall follow closely the inductive construction of the desingularization given in \cite{CED}.




Before we begin with that, let us first note that $ t_1^{i} t_2^{j} $ is $\boldsymbol{\mu_n}$-invariant if and only if $ i + r j = n k(i,j) $, for some integer $ k(i,j) $. From Lemma \ref{lemma 3.1}, it follows that $ [ \xi ] (t_1^i t_2^j) = \xi^{i \alpha_1} t_1^{i} t_2^{j} $, for any $ \xi \in G $, where $ m_1 \alpha_1 \equiv_n 1 $.


\subsection{Changing coordinates}\label{subs. 3.5}
Define elements $ z = t_1^{n} $ and $ w = t_2/t_1^r $. These are $\boldsymbol{\mu}_n$-invariant elements in the fraction field of $T$. Understanding how $G$ acts on these elements will be a key step in describing the action on the desingularization, so we work that out now.

Note that $ z = t_1^{n} = v_1 $, so $ [ \xi ](z) = z $. Moreover, since $ zw = t_1^{n-r} t_2 $ is an invariant monomial, we get that $ [ \xi ](zw) = [ \xi ](t_1^{n-r} t_2) = \xi^{(n-r) \alpha_1} t_1^{n-r} t_2 = \xi^{ -r \alpha_1} zw $, and it follows that $ [ \xi ](w) = [ \xi ](zw/z) = \xi^{- r \alpha_1} w $.


\subsection{Blowing up}
In the following, we shall describe the minimal desingularization of $ Z = \Spec(U) $, which will be denoted  $ \varrho : \widetilde{Z} \rightarrow Z $. This map will consist of a series of blow-ups 
$$ \ldots \rightarrow Z_{i+1} \rightarrow Z_i \rightarrow \ldots \rightarrow Z_1 \rightarrow Z_0 = Z, $$
where each step $ \varrho_{i+1} : Z_{i+1} \rightarrow Z_i $ is the blow-up in a certain ideal supported in the special fiber. Furthermore, each step will produce a regular, affine open chart of $ \widetilde{Z} $, that is stable under the $G$-action. For each of these charts we will be able to describe the $G$-action explicitly.

To do this, we shall use the new coordinates $z$ and $w$ of $U$. Let $i,j$ be any pair of integers such that $ 0 \leq j \leq (n/r)i $. Then we have that $ z^i w^j \in R'[t_1,t_2] $, and furthermore 
$$ R'[t_1,t_2]^{\boldsymbol{\mu}_n} = \oplus_{0 \leq j \leq (n/r)i} R'z^iw^j. $$
A trivial computation shows that $ t_1^{m_1} t_2^{m_2} = z^{\mu_1}w^{m_2} $, where $ m_1 + r m_2 = n \mu_1 $, so we can write
$$ U = \frac{\oplus_{0 \leq j \leq (n/r)i} R'z^iw^j}{(z^{\mu_1}w^{m_2} - \pi')}. $$

The first map $ \varrho_1 : Z_1 \rightarrow Z $ in the sequence above is obtained by blowing up in the ideal generated by $z$ and $zw$. There are two affine charts covering $Z_1$, namely $D_+(z)$ and $ D_+(zw) $, where we adjoin the quotients $ zw/z = w $ and $ z/zw = 1/w $ respectively. We shall treat these two charts separately.

\subsection{The chart $D_+(z)$.}
We have
$$ D_+(z) = \Spec(U[w]) = \Spec(R'[z,w]/(z^{\mu_1}w^{m_2} - \pi')). $$
This affine piece is already regular, and the special fiber has two components. The component $ w = 0 $ is the strict transform of the branch with multiplicity $m_2$ in $Z$. The component $ z = 0 $ is the exceptional curve $E_1$, which has multiplicity $ \mu_1 $ in the special fiber.

$ D_+(z) $ is birational to $Z$ via $ \varrho_1 $. The coordinates $z$ and $w$ are elements of the function field of $Z$, and we computed above that $ [ \xi ] (z) = z $ and $ [ \xi ] (w) = \xi^{-r \alpha_1} w $, for any $ \xi \in G $. So the same is true for the action on the function field of $ D_+(z) $. We also have that $ z^{\mu_1} w^{m_2} =  t_1^{m_1} t_2^{m_2} $, hence $ [ \xi ] (z^{\mu_1} w^{m_2}) = \xi z^{\mu_1} w^{m_2} $, and so it follows that $ D_+(z) $ is stable under the $G$-action.

On the chart $ D_+(z) $, the exceptional curve $E_1$ has affine ring $k[w]$, and the $G$-action restricted to $E_1$ is given by $ [ \xi ] (w) = \xi^{-r \alpha_1} w $. 

\subsection{The chart $D_+(zw)$.}
We have that $ D_+(zw) = \Spec(U[1/w]) $, where
$$ U[1/w] = \frac{\oplus_{j \leq (n/r)i, 0 \leq i} R'z^iw^j}{(z^{\mu_1}w^{m_2} - \pi')}. $$

We will now change coordinates. We start by performing the first step of the Jung-Hirzebruch continued fraction expansion. So we write $ n = b_1 r - r_1 $, where $ b_1 = \left \lceil n/r \right \rceil \geq 2 $. There are now two possibilities; either
\begin{enumerate}
\item $ r = 1 $ and $ r_1 = 0 $, or
\item $ r_1 > 0 $ and $ \Gcd(r,r_1) = 1 $.
\end{enumerate}

Consider first case (i), where $ r = 1 $ and $ r_1 = 0 $. Then we have $ b_1 = n $, so $ b_1 \mu_1 - m_2 = m_1 $. Let $ i_1 = b_1 i - j $ and $ j_1 = i $, and define new coordinates $ z_1 = 1/w $ and $ w_1 = z w^{b_1} $. By an easy computation, we see that $ z^i w^j = z_1^{i_1} w_1^{j_1} $. But then it follows that
$$ U[1/w] = R'[z_1,w_1]/(z_1^{b_1 \mu_1 - m_2} w_1^{\mu_1} - \pi') = R'[z_1,w_1]/(z_1^{m_1} w_1^{\mu_1} - \pi'), $$ so in particular $ \Spec(U[1/w]) $ is regular. The special fiber consists of $ z_1 = 0 $, which is the strict transform of the branch of $Z_k$ with multiplicity $m_1$, and the component $ w_1 = 0 $, the exceptional curve $ E_1 $.

Since $ \Spec(U[1/w]) $ is birational to $Z$, and since we know the action on the elements $z$ and $w$, the action on $z_1$ and $w_1$ follows immediately. It is $ [ \xi ](z_1) = \xi^{ r \alpha_1 } z_1 $ and $ [ \xi ](w_1) = \xi^{ - b_1 r \alpha_1 } w_1 = w_1 $, since $ b_1 = n $. As $ z^{\mu_1} w^{m_2} = z_1^{m_1} w_1^{\mu_1} $, we get that $ [ \xi ] (z_1^{m_1} w_1^{\mu_1}) = \xi z_1^{m_1} w_1^{\mu_1} $, so $ \Spec(U[1/w]) $ is stable under the $G$-action. Note also that $E_1$ has affine ring $k[z_1]$ on this chart, and therefore the $G$-action restricted to $E_1$ is given by $ [ \xi ](z_1) = \xi^{ r \alpha_1 } z_1 $.

Let us now consider case (ii). Then we have that $ r_1 > 0 $ and $ \Gcd(r,r_1) = 1 $. We shall still do the coordinate change $ z_1 = 1/w $, $ w_1 = z w^{b_1} $. With the equation $ n = b_1 r - r_1 $ in mind, it follows that the conditions $ j \leq (n/r)i $ and $ 0 \leq i $ can be rewritten as $ 0 \leq i \leq (r/r_1)(b_1 i-j) $, hence  $ 0 \leq j_1 \leq (r/r_1) i_1 $. So we may write 
$$ U_1 := U[1/w] = \frac{\oplus_{0 \leq j_1 \leq (r/r_1)i_1} R'z_1^{i_1} w_1^{j_1}}{(z_1^{b_1 \mu_1 - m_2} w_1^{\mu_1} - \pi')}.$$
Notice that this is a ring essentially of the same type that we started with, where in addition to the coordinate change,  the parameters are changed by
$$ (n, r, m_1, m_2, \mu_1) \mapsto (r, r_1, m_1, \mu_1, b_1 \mu_1 - m_2). $$ 

We will now perform a new blow up, in the ideal $(z_1, z_1 w_1)$. Similarly to the case in Section \ref{subs. 3.5}, we get two affine charts. The chart $ D_+(z_1) $ is the spectrum of the ring
$$ R'[z_1,w_1]/(z_1^{\mu_2} w_1^{\mu_1} - \pi'), $$
where we have put $ \mu_2 = b_1 \mu_1 - m_2 $. $G$ acts on this ring by $ [ \xi ](z_1) = \xi^{r \alpha_1} z_1 $ and $ [ \xi ](w_1) = \xi^{ - b_1 r \alpha_1} w_1 $. But as $ b_1 r = n + r_1 $, we get in fact that $ [ \xi ](w_1) = \xi^{ - r_1 \alpha_1} w_1 $, for any $ \xi \in G $. The new exceptional curve $E_2$ has affine ring $k[w_1]$ on this chart, and hence the restricted $G$-action is given by $ [ \xi ](w_1) = \xi^{ - r_1 \alpha_1} w_1 $.

The nature of the chart $ D_+(z_1w_1) $ depends entirely on our new set of parameters. That is, one computes the next step in the Jung-Hirzebruch continued fraction expansion $ r = b_2 r_1 - r_2 $. If $r_1 = 1$ (and hence $ r_2 = 0 $), also $ D_+(z_1w_1) $ is regular. In this case, we change coordinates in a similar fashion as above. Since $ z_1, w_1 $ are ``eigenelements'' under the $G$-action, the same will be true for the new coordinates $ z_2, w_2 $. In case $ r_2 > 0$, we change coordinates, and get a new ring of the same type as above, after changing our set of parameters.

\subsection{Induction step} 
As a preparation for the induction step, recall that when calculating the Jung-Hirzebruch continued fraction expansion of $ n/r $, we have a series of equations
\begin{equation} 
r_{l-1} = b_{l+1} r_l - r_{l+1}, 
\end{equation}
for all $ l \geq 0 $, where $ b_{l+1} = \left \lceil r_{l-1}/r_l \right \rceil $, and where we define $ r_{-1} = n $ and $ r_0 = r $. At each stage, we always have either $ r_l = 1 $ and $ r_{l+1} = 0 $, or $ r_{l+1} > 0 $ and $ \Gcd(r_l, r_{l+1}) = 1 $. Furthermore, $ r_{l-1} > r_l $, so this process will eventually stop after a finite number of steps. Let $L$ denote the \emph{length} of the resolution, by which we mean that $ r_{L-1} = 1 $, and $ r_L = 0 $.

We also have a system of equations 
\begin{equation} 
\mu_{l+1} = b_l \mu_l - \mu_{l-1} 
\end{equation}
for $ l \in \{ 1, \ldots, L \} $, where we define $ \mu_{0} = m_2 $, $ \mu_{1} = \mu = (m_1 + r m_2)/n $ and $ \mu_{L+1} = m_1 $.

Assume that we have reached the $(l-1)$-st step 
$$ Z_{l-1} \rightarrow \ldots \rightarrow Z_1 \rightarrow Z $$
in the blow-up procedure, and that $ Z_{l-1} $ is regular outside the open affine chart $\Spec(U_{l-1})$, where
$$ U_{l-1} = \frac{\oplus_{0 \leq j_{l-1} \leq (r_{l-2}/r_{l-1}) i_{l-1}} R' z_{l-1}^{i_{l-1}} w_{l-1}^{j_{l-1}}}{(z_{l-1}^{\mu_l} w_{l-1}^{\mu_{l-1}} - \pi')}. $$
Assume furthermore that the $G$-action is given by $ [ \xi ](z_{l-1}) = \xi^{\alpha_1 r_{l-2}} z_{l-1} $, $ [ \xi ](w_{l-1}) = \xi^{- \alpha_1 r_{l-1}} w_{l-1} $ and $ [ \xi ](z_{l-1}^{\mu_l} w_{l-1}^{\mu_{l-1}}) = \xi z_{l-1}^{\mu_l} w_{l-1}^{\mu_{l-1}} $.

We will now blow up in the ideal $ (z_{l-1}, z_{l-1} w_{l-1}) $. Then we get two affine charts $ D_+(z_{l-1}) $ and $ D_+(z_{l-1}w_{l-1}) $. The affine ring for the chart $ D_+(z_{l-1}) $ is 
$$ U_{l-1}[w_{l-1}] = R'[z_{l-1}, w_{l-1}]/( z_{l-1}^{\mu_l} w_{l-1}^{\mu_{l-1}} - \pi' ). $$
Notice that this chart is regular, and that by assumption, we already know the $G$-action here.

The affine ring for the chart $ D_+(z_{l-1}w_{l-1}) $ is
$$ U_{l-1}[1/w_{l-1}] = \frac{\oplus_{j_{l-1} \leq (r_{l-2}/r_{l-1}) i_{l-1}, 0 \leq i_{l-1}} R' z_{l-1}^{i_{l-1}} w_{l-1}^{j_{l-1}}}{(z_{l-1}^{\mu_l} w_{l-1}^{\mu_{l-1}} - \pi')}. $$
Let us introduce the new indexing $ i_l = b_l i_{l-1} - j_{l-1} $, $ j_l = i_{l-1} $, and define new coordinates $ z_l = 1/w_{l-1} $, $ w_l = z_{l-1} w_{l-1}^{b_l} $. There are two possibilities, that we will treat separately.

We shall first treat the case where $ r_{l-1} = 1 $ and $ r_l = 0 $. Then it is easy to compute that $ z_{l-1}^{\mu_l} w_{l-1}^{\mu_{l-1}} = z_{l}^{\mu_{l+1}} w_{l}^{\mu_{l}} $. We now have that $ b_l = r_{l-2} $, so the inequality $ j_{l-1} \leq (r_{l-2}/r_{l-1}) i_{l-1} $ can be written as $ i_l = b_l i_{l-1} - j_{l-1} \geq 0 $. It follows that $ i_l, j_l $ run over non-negative integers, and hence we get that
$$ U_{l-1}[1/w_{l-1}] = R'[z_{l}, w_{l}]/( z_{l}^{\mu_{l+1}} w_{l}^{\mu_{l}} - \pi' ), $$
which is regular. It follows that $ [ \xi ](z_{l}) = \xi^{\alpha_1 r_{l-1}} z_l $ and $ [ \xi ](w_{l}) = \xi^{\alpha_1 r_{l-2} - \alpha_1 b_l r_{l-1}} w_{l} = \xi^{- \alpha_1 r_l } w_{l} = w_{l} $, since $ r_l = 0 $. Furthermore, we see that $ \Spec(U_{l-1}[1/w_{l-1}]) $ is stable under the $ G $-action.

In the second case we have $ r_l > 0 $ and $ \Gcd(r_{l-1}, r_l) = 1 $. Since $ r_{l-2} = b_l r_{l-1} - r_l $, the inequalities $ j_{l-1} \leq (r_{l-2}/r_{l-1}) i_{l-1} $ and $ 0 \leq i_{l-1} $ can be written on the form
$$ 0 \leq j_l \leq (r_{l-1}/r_l) i_l. $$
It follows that we can write
$$ U_l := U_{l-1}[1/w_{l-1}] = \frac{\oplus_{0 \leq j_l \leq (r_{l-1}/r_l) i_l} R' z_l^{i_l} w_l^{j_l}}{(z_{l}^{\mu_{l+1}} w_{l}^{\mu_{l}} - \pi' )}. $$
Furthermore, $ [ \xi ](z_{l}) = \xi^{\alpha_1 r_{l-1}} z_l $ and $ [ \xi ](w_{l}) = \xi^{- \alpha_1 r_l } w_{l} $ for any $ \xi \in G $.

By induction, it now follows that we may cover $ \widetilde{Z} $ with a finite number of affine open charts $ \Spec(R'[z_l,w_l]/(z_{l}^{\mu_{l+1}} w_{l}^{\mu_{l}} - \pi')) $. These charts are stable under the $G$-action, and we have that $ [ \xi ] (z_l) = \xi^{\alpha_1 r_{l-1}} z_l $ and $ [ \xi ] (w_l) = \xi^{- \alpha_1 r_{l}} w_l $.

\subsection{}\label{3.11}

We sum up these results in Proposition \ref{prop. 3.3} below. This proposition will be important in later sections, when we consider the $G$-action on the cohomology groups $ H^i(\mathcal{Y}_k, \mathcal{O}_{\mathcal{Y}_k}) $.

\begin{prop}\label{prop. 3.3}
The minimal desingularization $ \widetilde{Z} $ of $ Z = \Spec(U) $ can be covered by the affine charts
$$ U_{l-1}[w_{l-1}] = R'[z_{l-1},w_{l-1}]/(z_{l-1}^{\mu_{l}} w_{l-1}^{\mu_{l-1}} - \pi'), $$
for $ l \in \{ 1, \ldots, L \} $, together with the affine chart 
$$ U_{L-1}[1/w_{L-1}] = R'[z_L,w_L]/(z_L^{\mu_{L+1}} w_L^{\mu_{L}} - \pi'). $$

These charts are $ G $-stable, and the $ G $-action is given by
$ [ \xi ] (\pi') = \xi \pi' $, $ [ \xi ] (z_{l-1}) = \xi^{ \alpha_1  r_{l-2} } z_{l-1} $ and $ [ \xi ] (w_{l-1}) = \xi^{ - \alpha_1 r_{l-1} } w_{l-1} $, for all $ l \in \{ 1, \ldots, L+1 \} $, for any $ \xi \in G $.

Let $ E_l $ be the $ l $-th exceptional component. On the chart $U_{l-1}[w_{l-1}]$, we have that the affine ring for $ E_l $ is $ k[w_{l-1}] $, and $ G $ acts by $ [ \xi ] (w_{l-1}) = \xi^{ - \alpha_1 r_{l-1} } w_{l-1} $, for any $ \xi \in G $. On the affine chart $U_{l}[w_{l}]$, the affine ring for $ E_l $ is $ k[z_l] $, and $ G $ acts by $ [ \xi ] (z_{l}) = \xi^{ \alpha_1  r_{l-1} } z_{l} $.
\end{prop}

\vspace{0.3cm}

Let us finally remark that the cotangent space to $ \widetilde{Z} $ at the fixed point that is the intersection point of $ E_l $ and $ E_{l+1} $ is generated by (the classes) of the local equations $ z_l $ and $ w_l $ for the curves. Therefore, Proposition \ref{prop. 3.3} gives a complete description of the action on the cotangent space. Furthermore, we can also read off the eigenvalues for the elements of this basis. Hence we immediately get an explicit description of the action on the cotangent space to the minimal desingularization of $\mathcal{X}'$ at the corresponding fixed point.

\section{Action on cohomology}
In the previous two sections, we have described the $G$-action on $ \mathcal{Y}/S' $. We will now begin the study of the induced $G$-action on the cohomology groups $ H^i(\mathcal{Y}_k, \mathcal{O}_{\mathcal{Y}_k}) $. 

We begin this section with recalling some generalities about coherent sheaf cohomology, and establish notation and terminology that will be used throughout the rest of the paper. 

Let $X$ be a scheme, and let $ \mathcal{F} $ be a coherent sheaf on $X$. Assume that we have an automorphism $ g : X \rightarrow X $, and an isomorphism $ u : g^* \mathcal{F} \rightarrow \mathcal{F} $. This data will, as we explain below, give an automorphism of the cohomology groups $ H^i(X,\mathcal{F}) $. If $X$ is projective over a field $k$, these groups are finite dimensional vector spaces over $k$, and hence the traces of these automorphisms are defined. However, we will typically work over fields of positive characteristic, and then traces can behave badly, in particular with respect to representation theoretical questions. In our case, we shall see below that this can be solved by lifting to characteristic $0$, using the theory of \emph{Brauer characters}. 

We end this section with considering the case where $X$ is a smooth curve, and where $ \mathcal{F} $ is an invertible sheaf. We give, in Proposition \ref{prop. 5.5} and Proposition \ref{prop. 5.6}, formulas for the alternating sum of the Brauer traces of automorphisms as mentioned above of the $ H^i(X,\mathcal{F}) $, in terms of certain data related to the fixed point locus. These formulas will be used numerous times in later sections.

\subsection{The general set-up}\label{inverse}
Let $ f : X \rightarrow Y $ be a morphism of schemes, and let $ \mathcal{G} $ be a sheaf on $Y$. Let us first recall the construction of the inverse image sheaf $ f^{-1} \mathcal{G} $ on $X$. It is obtained as the sheafification of the presheaf
$$ U \mapsto \underset{\rightarrow}{\Lim} ~ \mathcal{G}(V), $$
where the direct limit is taken over all open sets $ V \subseteq Y $ such that $ U \subseteq f^{-1}(V) $.

In particular, we have a canonical homomorphism of global sections
\begin{equation}\label{equation 5.1} 
\Gamma (Y, \mathcal{G}) \rightarrow \Gamma(X, f^{-1} \mathcal{G}). 
\end{equation}

\begin{rmk}We shall only use this construction when $f$ is an isomorphism, in which case the functor $ f^{-1} $ coincides with the direct image functor $ \varphi_* $, where $ \varphi $ denotes the inverse morphism to $ f $. Hence there is no need to sheafify.
\end{rmk}



\begin{lemma}\label{lemma 5.2}
The morphism $ f : X \rightarrow Y $ induces a natural and canonical homomorphism
$$ H^p(f) : H^p(Y,\mathcal{G}) \rightarrow H^p(X, f^* \mathcal{G}), $$  
for all $p \geq 0$. 
\end{lemma}
\begin{pf}
Let 
$$ 0 \rightarrow \mathcal{G} \rightarrow \mathcal{I}^{\bullet} $$  
be an injective resolution in $ \mathcal{A}b(Y) $. Since the functor $ f^{-1} $ is exact, we get that 
$$ 0 \rightarrow f^{-1} \mathcal{G} \rightarrow f^{-1} \mathcal{I}^{\bullet} $$
is a resolution. Let us now choose an injective resolution
$$ 0 \rightarrow f^{-1} \mathcal{G} \rightarrow \mathcal{J}^{\bullet} $$
in $ \mathcal{A}b(X) $. Then there is a map
$$ \xymatrix{0 \ar[r] & f^{-1} \mathcal{G} \ar[r] \ar[d] & f^{-1} \mathcal{I}^{\bullet} \ar[d] \\
             0 \ar[r] & f^{-1} \mathcal{G} \ar[r]        & \mathcal{J}^{\bullet}, } $$
unique up to homotopy. Hence, on global sections, we get
$$ \Gamma(Y, \mathcal{I}^{\bullet}) \rightarrow \Gamma(X, f^{-1} \mathcal{I}^{\bullet}) \rightarrow \Gamma(X, \mathcal{J}^{\bullet}), $$
which induces a map
\begin{equation}\label{eq. 5.2}
H^p(Y,\mathcal{G}) \rightarrow H^p(X, f^{-1} \mathcal{G}), 
\end{equation}
for every $ p \geq 0 $. This map is independent of the choice of injective resolutions. 

The canonical homomorphism
$$ f^{-1} \mathcal{G} \rightarrow f^* \mathcal{G} = f^{-1} \mathcal{G} \otimes_{f^{-1} \mathcal{O}_Y} \mathcal{O}_X $$
of sheaves on $X$ induces by functoriality a map
\begin{equation}\label{eq. 5.3} 
H^p(X, f^{-1} \mathcal{G}) \rightarrow H^p(X, f^* \mathcal{G}),  
\end{equation}
for every $ p \geq 0 $. The composition of the maps \ref{eq. 5.2} and \ref{eq. 5.3} is then the desired $H^p(f)$.
\end{pf}

\vspace{0.5cm}

The inverse image $ f^{-1} : \mathcal{A}b(Y) \rightarrow \mathcal{A}b(X) $ is an exact functor. Hence, the assignment $ \mathcal{G} \mapsto H^{\bullet}(X, f^{-1} \mathcal{G}) $ makes up a $ \delta $-functor from $ \mathcal{A}b(Y) $ to $ \mathcal{A}b $. Furthermore, the canonical map
$$ \phi^0 : \Gamma (Y, \mathcal{G}) \rightarrow \Gamma(X, f^{-1} \mathcal{G}), $$ 
induces a unique sequence $ \phi^p : H^p(Y,\mathcal{G}) \rightarrow H^p(X, f^{-1} \mathcal{G}) $ of morphisms, for every $ p \geq 0 $. Given any sequence
$$ 0 \rightarrow \mathcal{G}_1 \rightarrow \mathcal{G}_2 \rightarrow \mathcal{G}_3 \rightarrow 0 $$
in $ \mathcal{A}b(Y) $, the maps $ \phi^p $ commute with the differentials in the long exact sequences induced by $ H^p(Y,-) $ and $ H^p(X, f^{-1} (-)) $ respectively.

Consider now the case where $ f : X \rightarrow Y $ is an \emph{isomorphism}. Then the pullback $ f^* $ is an exact functor. Consequently, any exact sequence
$$ 0 \rightarrow \mathcal{G}_1 \rightarrow \mathcal{G}_2 \rightarrow \mathcal{G}_3 \rightarrow 0 $$
of $ \mathcal{O}_Y $-modules gives an exact sequence
$$ 0 \rightarrow f^* \mathcal{G}_1 \rightarrow f^* \mathcal{G}_2 \rightarrow f^* \mathcal{G}_3 \rightarrow 0 $$ 
of $ \mathcal{O}_X $-modules, and a commutative diagram

\begin{equation} \xymatrix{
0 \ar[r] & f^{-1} \mathcal{G}_1 \ar[r] \ar[d] & f^{-1} \mathcal{G}_2 \ar[r] \ar[d] & f^{-1} \mathcal{G}_3 \ar[r] \ar[d] & 0 \\
0 \ar[r] & f^* \mathcal{G}_1 \ar[r] & f^* \mathcal{G}_2  \ar[r] & f^* \mathcal{G}_2 \ar[r] & 0, }
\end{equation}
in $ \mathcal{A}b(X) $, where the vertical arrows are the canonical maps discussed in the proof of Lemma \ref{lemma 5.2}. Passing to cohomology, there are then induced maps 
$$ \psi^p : H^p(X, f^{-1} \mathcal{G}_i) \rightarrow H^p(X, f^* \mathcal{G}_i), $$ 
commuting with the differentials 
$$ H^p(X, f^{-1} \mathcal{G}_3) \rightarrow H^{p+1}(X, f^{-1} \mathcal{G}_1) $$
and 
$$ H^p(X, f^* \mathcal{G}_3) \rightarrow H^{p+1}(X, f^* \mathcal{G}_1), $$
for any $ p \geq 0 $ in the two long exact sequences.

Thus, if $ f : X \rightarrow Y $ is an isomorphism, and if 
$$ 0 \rightarrow \mathcal{G}_1 \rightarrow \mathcal{G}_2 \rightarrow \mathcal{G}_3 \rightarrow 0 $$
is a short exact sequence of $ \mathcal{O}_Y $-modules, the maps 
$$ H^p(f) : H^p(Y,\mathcal{G}_i) \rightarrow H^p(X, f^* \mathcal{G}_i) $$ 
constructed in Lemma \ref{lemma 5.2} commute with the differentials of the two long exact sequences. 

\subsection{Automorphisms of schemes}
Let $X$ be a scheme, let $ g : X \rightarrow X $ be a morphism, and let $ \mathcal{F} $ be a sheaf of $ \mathcal{O}_X $-modules. Assume that we are also given a homomorphism 
$$ u : g^* \mathcal{F} \rightarrow \mathcal{F} $$
of $ \mathcal{O}_X $-modules. By functoriality, $u$ induces a homomorphism on the cohomology groups 
$$ H^p(u) : H^p(X, g^* \mathcal{F}) \rightarrow H^p(X,\mathcal{F}), $$ 
for all $p \geq 0$.

In particular, if the map $ u : g^* \mathcal{F} \rightarrow \mathcal{F} $ is an isomorphism, then the induced map $ H^p(u) $ is also an isomorphism for all $p$. 

\begin{dfn}\label{dfn 5.4}
Let $ g : X \rightarrow X $ be an automorphism of a scheme $X$, $ \mathcal{F} $ a sheaf of $ \mathcal{O}_X $-modules, and $ u : g^* \mathcal{F} \rightarrow \mathcal{F} $ a homomorphism. The endomorphism 
$$ H^p(g,u) : H^p(X,\mathcal{F}) \rightarrow H^p(X,\mathcal{F}) $$
\emph{induced} by the couple $ (g,u) $ is defined as the composition of the maps $ H^p(g) $ and $ H^p(u) $.
\end{dfn}

\vspace{0.3cm}

In case $ \mathcal{F} = \mathcal{O}_X $, there is a canonical isomorphism $ g^* \mathcal{O}_X \cong \mathcal{O}_X $, associated to the morphism $g$. So we get naturally an endomorphism of the cohomology groups
 $$ H^p(g) : H^p(X,\mathcal{O}_X) \rightarrow H^p(X,\mathcal{O}_X), $$
for all $p \geq 0$.

\subsection{Short exact sequences}\label{5.4}
Let $ g : X \rightarrow X $ be an automorphism, and consider a commutative diagram
$$ \xymatrix{
 0 \ar[r] &  g^* \mathcal{F}_1 \ar[r] \ar[d]_{u_1} &  g^* \mathcal{F}_2 \ar[r] \ar[d]_{u_2} & g^* \mathcal{F}_3 \ar[r] \ar[d]_{u_3} & 0 \\
 0 \ar[r] & \mathcal{F}_1 \ar[r] & \mathcal{F}_2 \ar[r] & \mathcal{F}_3 \ar[r] & 0,} $$
where the horizontal lines are exact, and the maps $ u_i $ are homomorphisms of $ \mathcal{O}_X $-modules. The induced maps $ H^p(u_i) : H^p(X,g^* \mathcal{F}_i) \rightarrow H^p(X, \mathcal{F}_i) $ commute with the differentials of the induced long exact sequences in cohomology. 

We may now combine the maps $ H^p(u_i) $ with the maps $ H^p(g) $ constructed earlier, and obtain a commutative diagram

\begin{equation} 
\xymatrix{
\ldots \ar[r] & H^p(X, \mathcal{F}_2) \ar[r] \ar[d]^{H^p(g,u_2)} & H^p(X, \mathcal{F}_3) \ar[r]^{\delta} \ar[d]^{H^p(g,u_3)} & H^{p+1}(X, \mathcal{F}_1) \ar[r] \ar[d]^{H^{p+1}(g,u_1)} & \ldots \\
\ldots \ar[r] & H^p(X, \mathcal{F}_2) \ar[r] & H^p(X, \mathcal{F}_3) \ar[r]^{\delta} & H^{p+1}(X, \mathcal{F}_1) \ar[r] & \ldots .}
\end{equation}



\subsection{}
If our scheme $X$ in addition is projective over a field $k$, and $ \mathcal{F} $ is a coherent sheaf of $ \mathcal{O}_X $-modules, the cohomology groups $ H^p(X,\mathcal{F}) $ are finite dimensional $k$-vector spaces. We can then define the trace $ \Tr(H^p(g,u)) $ of the endomorphism $ H^p(g,u) $. We will use the notation
$$ \Tr(e(H^{\bullet}(g,u))) = \sum_{p=0}^{dim(X)} (-1)^p ~ \Tr(H^p(g,u)) $$
for the alternating sum of the traces.

\subsection{Witt vectors and Teichm\"uller liftings}
In the case where $ p = \Char(k) > 0 $, we let $W(k)$ denote the ring of \emph{Witt vectors} for $k$. Recall that $W(k)$ is a complete discrete valuation ring, $p$ is a uniformizing parameter in $W(k)$ and the residue field is $k$ (\cite{Serre}, Chap. II, par. 5). 

Let $ FW(k) $ denote the fraction field of $W(k)$. An important feature of the Witt vectors is the fact that $ \Char(FW(k)) = 0 $. Furthermore, there exists a unique multiplicative map $ w : k \rightarrow W(k) $ that sections the reduction map $ W(k) \rightarrow k $. The map $w$ is often referred to as the \emph{Teichm\"uller lifting} of $k$ to $ W(k) $. The existence of this map often makes it possible to lift computations from characteristic $p$ to characteristic $0$. 

Since $k$ is assumed to be algebraically closed, it follows that $k$ has a full set of $n$-th roots of unity, for any $n$ not divisible by $p$. As $W(k)$ is complete, these lift uniquely to $W(k)$. Furthermore, reduction modulo $p$ induces an isomorphism of $ \boldsymbol{\mu}_n(W(k)) $ onto $ \boldsymbol{\mu}_n(k) $. Hence, if $ \lambda \in \boldsymbol{\mu}_n(k) $ is any $n$-th root of unity, it follows that $ w(\lambda) \in W(k) $ is an $n$-th root of unity, of the same order as $ \lambda $.

\subsection{Brauer characters}
A few facts regarding \emph{modular characters} are needed, and are stated here in the case where $ G = \boldsymbol{\mu}_n $. We refer to \cite{SerreLin}, Chap. 18 for details. It should be mentioned that most of the constructions and properties below are actually valid for any finite group.  

If $ E $ is a $ k[G] $-module, we let $ g_E $ denote the endomorphism of $ E $ induced by $ g \in G $. Since the order of $g$ divides $n$, and $n$ is relatively prime to $p$, it follows that $g_E$ is diagonalizable, and that all the eigenvalues $ \lambda_1, \ldots, \lambda_{e=\Dim E} $ are $n$-th roots of unity. The \emph{Brauer character} is defined by assigning 
$$ \phi_E(g) = \sum_{i=1}^{e} w(\lambda_i), $$   
where $ w(\lambda_i) $ is the Teichm\"uller lift of $ \lambda_i $ to $ W(k) $. It can be seen that the function
$$ \phi_E : G \rightarrow W(k) $$
thus obtained is a class function on $G$. We shall call the element $ \phi_E(g) \in W(k) $ the \emph{Brauer trace} of $g_E$.  

Two properties of the Brauer character are of particular importance to us. First, we have that 
$$ \overline{\phi_E(g)} = \Tr(g_E). $$
That is, we obtain the ordinary trace from the Brauer trace by reduction modulo $p$. 

Another important property is that the Brauer character is additive on short exact sequences. That is, if 
$$ 0 \rightarrow E' \rightarrow E \rightarrow E'' \rightarrow 0 $$ 
is an exact sequence of $ k[G] $-modules, then $ \phi_E = \phi_{E'} + \phi_{E''} $. A useful consequence is that if 
$$ 0 \rightarrow E_0 \rightarrow \ldots \rightarrow E_i \rightarrow \ldots \rightarrow E_l \rightarrow 0 $$
is an exact sequence of $k[G]$-modules, we get that $ \sum_{i=0}^l (-1)^i \phi_{E_i}(g) = 0 $. That is, the alternating sum of the Brauer traces equals zero.



\begin{ntn}
If $ V $ is a finite dimensional vector space over $k$, and $ \psi : V \rightarrow V $ is an automorphism, we will use the notation $ \Tr_{\beta}(\psi) $ for the Brauer trace of $\psi$.
\end{ntn}

\subsection{Automorphisms of curves and invertible sheaves}
We will now consider the situation where $ X $ is a smooth, connected and projective curve $C$ over an algebraically closed field $k$, and where $\mathcal{F} = \mathcal{L}$ is an invertible sheaf on $C$. 

Let $ g : C \rightarrow C $ be an automorphism, and let $ u : g^*\mathcal{L} \rightarrow \mathcal{L} $ be an isomorphism, so that there are induced automorphisms $ H^p(g,u) $ for $ p = 0, 1 $. We would like to compute the alternating sum $ \sum_{p=0}^1 (-1)^p ~\Tr_{\beta} (H^p(g,u)) $ of the Brauer traces.  

\subsection{Computing the trace when $g$ is trivial}
Let us first consider the case when the automorphism $ g : C \rightarrow C $ is trivial, i.e., $ g = \Id_C $. Then $ H^p(g) $ is the identity, so we need only consider the isomorphism $ u : \mathcal{L} = g^*(\mathcal{L}) \rightarrow \mathcal{L} $. Hence $ H^p(g,u) = H^p(u) : H^p(C,\mathcal{L}) \rightarrow H^p(C,\mathcal{L}) $, where $ u \in \Aut_{\mathcal{O}_C}(\mathcal{L}) $. The results we shall need for this situation are listed in Proposition \ref{prop. 5.5} below.

\begin{prop}\label{prop. 5.5}
Let $ C $ be a smooth, projective and irreducible curve over an algebraically closed field $k$, and let $ \mathcal{L} $ be an invertible sheaf on $C$. Then we have that 
\begin{enumerate}
\item $ \Aut_{\mathcal{O}_C}(\mathcal{L}) = k^* $.
\item Let $ \lambda \in \Aut_{\mathcal{O}_C}(\mathcal{L}) $. Then the induced automorphism 
$$ H^p(\lambda) : H^p(C,\mathcal{L}) \rightarrow H^p(C,\mathcal{L}) $$ 
is given by multiplication by $ \lambda $ for all $ p \geq 0 $.
\item The following equality holds in $W(k)$:
$$ \Tr_{\beta}(e(H^{\bullet}(\lambda))) = w(\lambda) \cdot (\Deg_C(\mathcal{L}) + 1 - p_a(C)). $$
\end{enumerate}
\end{prop}
\begin{pf}
The statement in (i) is standard, and follows essentially from the fact that $ \mathcal{O}_C^*(C) = k^* $, and that any automorphism of $ \mathcal{O}_C $ is determined by $ \mathcal{O}_C^*(C) $.

Part (ii) can be seen using \v{C}ech cohomology.

Finally, (iii) follows easily, since according to (ii), all eigenvalues of $ H^p(\lambda) $ equal $ \lambda $, for each $p$. So we get that
$$  \Tr_{\beta}(H^p(\lambda)) = w(\lambda) \cdot h^p(C,\mathcal{L}), $$
and hence, by the Riemann-Roch theorem, 
$$ \Tr_{\beta}(e(H^{\bullet}(\lambda))) = w(\lambda) \cdot (h^0(C,\mathcal{L}) - h^1(C,\mathcal{L}) ) = w(\lambda) \cdot (\Deg_C(\mathcal{L}) + 1 - p_a(C)). $$ 
\end{pf}

\subsection{Computing the trace when $g$ is nontrivial}
If the automorphism $ g : C \rightarrow C $ is non-trivial, things become substantially harder.  In this case, we will need the \emph{Lefschetz-Riemann-Roch} formula. This formula is presented in Section \ref{LRR} below, and we also discuss the ingredients in the formula. We will follow the paper \cite{Don} by P. Donovan closely.   

\subsection{The Lefschetz-Riemann-Roch formula}\label{LRR}
Let $ Y $ be a smooth, projective variety over an algebraically closed field $k$. An automorphism $ y \in \Aut_k(Y) $ is said to be \emph{periodic} if $ y^n = \Id_Y $ for some $ n \in \mathbb{N} $. The smallest such $n$ will be denoted the period of $y$. We shall always make the assumption that the characteristic of $k$ does not divide $n$. It is a basic fact (\cite{Don}, Lemma 4.1) that any connected component $ Z $ of the fixed point set of $y$ is non-singular. We shall say that $Z$ is a \emph{fixed} component. 

Consider the set of isomorphism classes of homomorphisms $ \eta : y^* \mathcal{F} \rightarrow \mathcal{F} $, where $ \mathcal{F} $ is a coherent sheaf on $ Y $. Let $ M(y) $ denote the quotient of the free abelian group on this set, modulo certain relations that will not be written out explicitly here. In the special case where $ y = \Id_Y $, we write $ M(Y) = M(\Id_Y) $. 

\begin{rmk} 
$ M(y) $ can be seen to be generated by homomorphisms $ \eta : y^* \mathcal{F} \rightarrow \mathcal{F} $, where $ \mathcal{F} $ is a locally free sheaf. Also, $ M(y) $ has the structure of a $ \mathbb{Z}[k] $-algebra in a natural way. So one may define $ M'(y) := M(y) \otimes_{\mathbb{Z}[k]} W(k) $. We refer to \cite{Don} for details.
\end{rmk}

Recall from Definition \ref{dfn 5.4}, that a couple $ y, \eta $ induces endomorphisms $ H^p(y, \eta) $ on the cohomology groups $ H^p(Y,\mathcal{F}) $, for every $ p \geq 0 $. Let $ c_{!} \eta := \sum_p (-1)^p ~[H^p(y, \eta)] $ denote the alternating formal sum of these endomorphisms. The \emph{Lefschetz-Riemann-Roch} formula (\cite{Don}, Theorem 5.4, Corollary 5.5) computes the alternating sum of the Brauer traces of the endomorphisms $ H^p(y, \eta) $ as a sum of certain contributions over the fixed components of $y$.

\begin{thm}[The Lefschetz-Riemann-Roch formula, \cite{Don}, Corollary 5.5]\label{thm 5.7}
Let $ (Y,y) $ be as above. If $Z$ is a fixed component of $y$, write $ i_Z : Z \rightarrow Y $ for the injection. Furthermore, write $ c_Z : Z \rightarrow \Spec(k) $ and $ c : Y \rightarrow \Spec(k) $ for the structure morphisms. Then the following equality holds in the ring $W(k)$ for $ \eta \in M(y) $ :
$$ \Tr_{\beta}(c_{!} \eta) = \sum_Z c_{Z*}(Td(Z) \cdot ct (\lambda_Z)^{-1} \cdot ct(i_Z^{!} \eta)). $$
\end{thm}

We shall now discuss the terms appearing on the right hand side of this formula. First we need to explain the concept of the \emph{Chern trace}, as defined in \cite{Don}, Chapter 2. Let $ X/k $ be a smooth, projective and irreducible variety. Let $ \mathcal{E} $ be a locally free sheaf on $X$, and consider a homomorphism $ \psi : \mathcal{E} \rightarrow \mathcal{E} $. The characteristic equation of $ \psi $ has coefficients in $ k = \mathcal{O}_X(X) $, and hence all roots are in $k$. Let $ t $ be one of the roots, and define
$$ \mathcal{E}_t := \Ker(\psi - t)^N, $$
for $ N > \Rk(\mathcal{E}) $. This is a locally free subsheaf of $ \mathcal{E} $, and independent of the chosen $N$. Furthermore, $ \mathcal{E} = \oplus_t \mathcal{E}_t $, and $ \psi = \oplus_t \psi_t $. 

It can now be shown that $ M(X) $ is generated by the images of the constant homomorphisms $ t : \mathcal{E} \rightarrow \mathcal{E} $, where $ t \in k $ and $ \mathcal{E} $ is locally free (\cite{Don}, Lemma 2.7). Let $K(X)$ be the Grothendieck group of locally free sheaves. Then the map 
$$ \alpha : M(X) \rightarrow K(X) \otimes \mathbb{Z}[k] $$
defined by $ \alpha([\psi]) = \sum_t [\mathcal{E}_t] \otimes [t] $ is a natural isomorphism (\cite{Don}, Proposition 2.8). 

Recall the construction of the \emph{Chern character} 
$$ ch : K(X) \rightarrow A^*(X)_{\mathbb{Q}} $$
(\cite{Fulton}, p. 282). This is a homomorphism of rings, and for a line bundle $L$ on $X$ (which is the only case where we will need an explicit description), it is given by
$$ ch [L] = \sum_{ i \geq 0 } (1/i!) c_1(L)^i. $$
The Chern trace
$$ ct : M(X) \rightarrow A^*(X) \otimes FW(k) $$
is then defined as the composition of $ \alpha $ and the induced map $ ch \otimes w $. If $X$ is a reducible variety, then all constructions above can be done on each separate component of $X$, so the Chern trace may still be defined.

If $ Z $ is a fixed component for the periodic automorphism $ y : Y \rightarrow Y $, $ \lambda_Z $ is defined as follows: Let $ i_Z : Z \rightarrow Y $ be the inclusion. Since both $Z$ and $ Y $ are smooth, this is a regular immersion, and hence the conormal sheaf $ \mathcal{C}_{Z/Y} := \mathcal{J}/\mathcal{J}^2 $ is a locally free sheaf on $ Z $, of rank equal to the codimension $ q $ of $ Z $ in $ Y $. Furthermore, $y$ induces a canonical endomorphism $ \Phi : \mathcal{C}_{Z/Y} \rightarrow \mathcal{C}_{Z/Y} $. Hence there are induced endomorphisms $ \Lambda^t \Phi $ on $ \Lambda^t \mathcal{C}_{Z/Y} $ for every $ t \in \{ 0, \ldots, q \} $. Then one defines 
$$ \lambda_Z = \sum_{t=0}^q (-1)^t [\Lambda^t \Phi] \in M(Z). $$
In particular, $ ct (\lambda_Z) $ is now defined, and Lemma 4.3 in \cite{Don} asserts that this is a unit in $ A^*(Z) \otimes FW(k) $.

Let $ \eta : y^* \mathcal{F} \rightarrow \mathcal{F} $ be an element of $M(y) $. Then $ i_Z^{!} \eta $ is defined as (the class of) the restriction 
$$ i_Z^* \eta : i_Z^* y^* \mathcal{F} = i_Z^* \mathcal{F} \rightarrow i_Z^* \mathcal{F}, $$
and is an element of $ M(Z) $. Hence the Chern trace $ ct(i_Z^{!} \eta) \in A^*(Z) \otimes FW(k) $ is defined.

Finally, $ Td(Z) $ is the \emph{Todd class} of $Z$ in $ A^*(Z) \otimes_{\mathbb{Z}} FW(k) $ (\cite{Fulton}, p. 354). 

\subsection{Formula in the case where $g$ is non-trivial}
Let $C/k$ be a smooth, connected and projective curve, and let $ g \in \Aut_k(C) $ be a non-trivial automorphism of finite period $n$, where $n$ is not divisible by the characteristic of $k$. Furthermore, assume that a homomorphism $ u : g^* \mathcal{L} \rightarrow \mathcal{L} $ is given, where $ \mathcal{L} $ is an invertible sheaf on $C$. Denote by $ H^p(g,u) : H^p(C,\mathcal{L}) \rightarrow H^p(C,\mathcal{L}) $ the induced endomorphism, for any $ p \geq 0 $. 

Theorem \ref{thm 5.7} above gives a formula for the alternating sum of the Brauer traces of the endomorphisms $ H^p(g,u) $. In order for this formula to be useful, we need to be able to compute the contributions from the fixed point locus. Proposition \ref{prop. 5.6} below gives an explicit description of these terms. 

\begin{prop}\label{prop. 5.6}
Let $C$, $ \mathcal{L} $, $g$ and $u$ be as above, and denote by $ C^g $ the (finite) set of fixed points of $g$. Then the following equality holds in $W(k)$: 
$$ \Tr_{\beta}(e(H^{\bullet}(g,u))) = \sum_{p=0}^1 (-1)^p ~\Tr_{\beta}(H^p(g,u)) = \sum_{z \in C^g} w(\lambda_u(z))/(1- w(\lambda_{dg}(z))), $$
where $ \lambda_u(z) $ is the unique eigenvalue of the fiber $ u(z) $ of the map $ u $ at $z$, and $ \lambda_{dg}(z) $ is the unique eigenvalue of the cotangent map on the fiber at $z$.
\end{prop}
\begin{pf}
Let us first note that if $ z \in C^g $, then $ z = \Spec(k) $, and the structure map $ c_{ z } : \{ z \} \rightarrow \Spec(k) $ is the identity. Furthermore, the Chern trace
$$ ct : M(\{z\}) = M(\Spec(k)) \rightarrow A^*(\Spec(k)) \otimes FW(k) $$
reduces in this case to the Brauer trace defined for endomorphisms of finite dimensional $k$-vector spaces (\cite{Don}, p. 264).

The Todd class $ Td(\{z\}) $ is just the class of a point in $ A^*(\Spec(k)) $ (\cite{Fulton}, Theorem 18.3 (5)), hence it is the identity element of $ A^*(\Spec(k)) $. 

The conormal bundle $ \mathcal{C} := \mathcal{C}_{\{z\}/C} $ can be computed as $ m_z/m_z^2 $, where $ m_z \in \mathcal{O}_{C,z} $ is the maximal ideal corresponding to the point $z$. The endomorphism of $ \Lambda^0 \mathcal{C} $ is trivial, and $ dg : \Lambda^1 \mathcal{C} \rightarrow \Lambda^1 \mathcal{C} $ is just an endomorphism of $ 1 $-dimensional vectorspaces, hence determined by its unique eigenvalue $ \lambda_{dg}(z) \in k $. Let $ w(1) = 1 $ and $ w(\lambda_{dg}(z)) $ be the Teichm\"uller lifts. Then we get that
$$ ct(\lambda_{\{z\}}) = 1 - w(\lambda_{dg}(z)) \in A^*(\Spec(k)) \otimes FW(k). $$  

It remains to compute $ ct(i_z^! u) $. But $ i_z^! u $ is (the class of) the pullback $ i_z^*u $, which is
$$ i_z^*u : i_z^* \mathcal{L} = i_z^* g^* \mathcal{L} \rightarrow i_z^* \mathcal{L}. $$
But this is just the fiber of $u$ at $z$, and therefore 
$$ ct(i_z^!u) = \Tr_{\beta}(u(z)) = w(\lambda_u(z)). $$
\end{pf}

 \begin{rmk}
The reader might want to compare Proposition \ref{prop. 5.6} with the \emph{Woods-Hole}-formula (\cite{SGA5}, Exp. III, Cor. 6.12), that gives a formula for the ordinary trace, instead of the Brauer trace.  
\end{rmk}

\begin{rmk} 
Throughout the rest of the text we will, when no confusion can arise, continue to write $ \lambda $ instead of $w(\lambda) $ for the Teichm\"uller lift of a root of unity $\lambda$.
\end{rmk}

\section{Action on the minimal desingularization}\label{section 6}
Recall the set-up in Section \ref{extensions and actions}. We considered an SNC-model $ \mathcal{X}/S $, and a tamely ramified extension $ S'/S $ of degree $n$ that is prime to the least common multiple of the multiplicities of the irreducible components of $ \mathcal{X}_k $. The minimal desingularization of the pullback $ \mathcal{X}_{S'}/S' $ is an SNC-model $ \mathcal{Y}/S' $, and the Galois group $ G = \boldsymbol{\mu_n} $ of the extension $ S'/S $ acts on $ \mathcal{Y} $. 

Let $ g \in G $ be any element, and consider the induced automorphism $ g : \mathcal{Y} \rightarrow \mathcal{Y} $. In Corollary \ref{g^{-1}(Z) = Z}, we saw that $ g^{-1}(Z) = Z $ for any effective divisor $ 0 \leq Z \leq \mathcal{Y}_k $, and so the $G$-action restricts to $Z$. In particular, there is an induced action by $G$ on $ \mathcal{Y}_k $.

In the previous section, we saw that for each $ g \in G $, the induced automorphism $ g : \mathcal{Y}_k \rightarrow \mathcal{Y}_k $ gave an automorphism of the cohomology groups $ H^i(\mathcal{Y}_k, \mathcal{O}_{\mathcal{Y}_k}) $, and hence $G$ acts on $ H^i(\mathcal{Y}_k, \mathcal{O}_{\mathcal{Y}_k}) $ for all $ i \geq 0 $. Our ultimate goal is to compute the irreducible characters for this representation on $ H^1(\mathcal{Y}_k, \mathcal{O}_{\mathcal{Y}_k}) $. To do this, we would ideally compute the Brauer trace of the automorphism of $ H^1(\mathcal{Y}_k, \mathcal{O}_{\mathcal{Y}_k}) $ induced by $g$, for every group element $ g \in G $. This information would then be used to compute the Brauer character. However, we can not do this directly. Instead we will compute the Brauer trace of the automorphism induced by $g$ on the formal difference $ H^0(\mathcal{Y}_k, \mathcal{O}_{\mathcal{Y}_k}) - H^1(\mathcal{Y}_k, \mathcal{O}_{\mathcal{Y}_k}) $, for any $ g \in G $. In our applications, we know the character for $ H^0(\mathcal{Y}_k, \mathcal{O}_{\mathcal{Y}_k}) $, so this would suffice in order to compute the character for $ H^1(\mathcal{Y}_k, \mathcal{O}_{\mathcal{Y}_k}) $.

The fact that $ \mathcal{Y}_k $ is not in general smooth, prevents us from using Proposition \ref{prop. 5.5} and Proposition \ref{prop. 5.6} directly. On the other hand, the irreducible components of $ \mathcal{Y}_k $ are smooth and proper curves. So we shall in fact show that it is possible to reduce to computing Brauer traces on each individual component of $ \mathcal{Y}_k $, where Proposition \ref{prop. 5.5} and Proposition \ref{prop. 5.6} do apply. The key step in obtaining this is to introduce a certain filtration of the special fiber $\mathcal{Y}_k$.

\subsection{Filtration of the special fiber} 
Let $ \{ C_{\alpha} \}_{\alpha \in \mathcal{A}} $ denote the set of irreducible components of $ \mathcal{Y}_k $, and let $ m_{\alpha} $ denote the multiplicity of $ C_{\alpha} $ in $ \mathcal{Y}_k $. Then $ \mathcal{Y}_k $ can be written in Weil divisor form as
$$ \mathcal{Y}_k = \sum_{\alpha} m_{\alpha} C_{\alpha}. $$

\begin{dfn}\label{complete}
A \emph{complete} filtration of $\mathcal{Y}_k$ is a sequence 
$$ 0 < Z_m < \ldots < Z_j < \ldots < Z_1 = \mathcal{Y}_k $$
of effective divisors $Z_j$ supported on $\mathcal{Y}_k$, such that for each $1 \leq j \leq m - 1 $ there exists an $ \alpha_j \in \mathcal{A} $ with $ Z_j - Z_{j+1} = C_{\alpha_j} $. So $ m = \sum_{\alpha} m_{\alpha} $.
\end{dfn}

Loosely speaking, such a filtration of $\mathcal{Y}_k$ is obtained by removing the irreducible components of the special fiber one at the time (counted with multiplicity, of course). 

\subsection{The steps in the filtration}
At each step of a complete filtration, we can construct an exact sequence of sheaves. 

\begin{lemma}\label{lemma 6.1}
Let $ 0 \leq Z' < Z \leq \mathcal{Y}_k $ be divisors such that $ Z - Z' = C $, for some irreducible component $C$ of $ \mathcal{Y}_k $. Denote by $ \mathcal{I}_Z $ and $  \mathcal{I}_{Z'} $ the corresponding ideal sheaves in $ \mathcal{O}_{\mathcal{Y}} $. Let $ i_{Z} $, $ i_{Z'} $ and $ i_C $ be the canonical inclusions of $ Z $, $Z'$ and $ C $ in $\mathcal{Y}$. Furthermore, let $ \mathcal{L} = i_C^*(\mathcal{I}_{Z'}) $.

We then have an exact sequence 
$$ 0  \rightarrow (i_C)_* \mathcal{L} \rightarrow (i_{Z})_* \mathcal{O}_{Z} \rightarrow (i_{Z'})_* \mathcal{O}_{Z'} \rightarrow 0 $$
of $ \mathcal{O}_{\mathcal{Y}} $-modules. 
\end{lemma}

\begin{pf}
The inclusions of the ideal sheaves $ \mathcal{I}_Z \subset \mathcal{I}_{Z'} \subset \mathcal{O}_{\mathcal{Y}} $ give an exact sequence
$$ 0 \rightarrow \mathcal{K} \rightarrow \mathcal{O}_{\mathcal{Y}}/\mathcal{I}_Z \rightarrow \mathcal{O}_{\mathcal{Y}}/\mathcal{I}_{Z'} \rightarrow 0, $$
where $ \mathcal{K} $ denotes the kernel. We can identify $ \mathcal{O}_{\mathcal{Y}}/\mathcal{I}_Z = (i_{Z})_* \mathcal{O}_{Z} $ and $ \mathcal{O}_{\mathcal{Y}}/\mathcal{I}_{Z'} = (i_{Z'})_* \mathcal{O}_{Z'} $, so it remains to consider $ \mathcal{K} $.

Let $ U = \Spec(A) \subset \mathcal{Y} $ be an open affine set. Then $A$ is a regular domain, and the ideal sheaves $ \mathcal{I}_C $, $ \mathcal{I}_Z $ and $ \mathcal{I}_{Z'} $ restricted to $U$ correspond to invertible modules $ I_C $, $ I_{Z} $ and $ I_{Z'} $ in $\Frac(A)$ (and in fact in $A$). Since $ Z = Z' + C $, we have that $ I_{Z} = I_C I_{Z'} $. 

Restricted to $U$, the exact sequence above is associated to the sequence
$$ 0 \rightarrow I_{Z'}/I_{Z} \rightarrow A/I_{Z} \rightarrow A/I_{Z'} \rightarrow 0. $$
We can now make the identification $ I_{Z'}/I_{Z} = I_{Z'}/I_C I_{Z'} = I_{Z'} \otimes_A A/I_C $. But $ I_{Z'} \otimes_A A/I_C $ is exactly the pullback of the ideal sheaf of $Z'$ to $C$. It is easy to see that these isomorphisms glue to give the desired isomorphism of the sheaves.
\end{pf}

\subsection{$G$-sheaves}\label{6.4}
Let $ \mathcal{F} $ be an $ \mathcal{O}_{\mathcal{Y}} $-module. We will say that $ \mathcal{F} $ is a $G$-sheaf if there are isomorphisms $ g^*\mathcal{F} \rightarrow \mathcal{F} $, for every $g \in G$, satisfying certain natural compatibility properties (cf. \cite{Kock}, Chapter 1).

Let $ \mathcal{F} $ and $ \mathcal{G} $ be $G$-sheaves. A homomorphism of $G$-sheaves $ \mathcal{F} \rightarrow \mathcal{G} $ is an $ \mathcal{O}_{\mathcal{Y}} $-homomorphism that commutes with the $G$-sheaf structures of $ \mathcal{F} $ and $ \mathcal{G} $ respectively.

\begin{lemma}\label{G-sheaf} 
Let $ 0 \leq Z \leq \mathcal{Y}_k $ be an effective divisor, with ideal sheaf $ \mathcal{I}_Z $. Then we have that $ \mathcal{I}_Z $ is a $G$-sheaf.
\end{lemma}
\begin{pf}
Let $ g \in G $ be any group element, and consider the corresponding automorphism $ g : \mathcal{Y} \rightarrow \mathcal{Y} $. Applying the exact functor $ g^{-1} $ to the inclusion $ \mathcal{I}_Z \subset \mathcal{O}_{\mathcal{Y}} $ gives an inclusion 
\begin{equation}\label{equation 6.1} 
g^{-1} \mathcal{I}_Z \subset g^{-1} \mathcal{O}_{\mathcal{Y}}. 
\end{equation}
Composing this inclusion with the canonical map $ g^{\sharp} : g^{-1} \mathcal{O}_{\mathcal{Y}} \rightarrow \mathcal{O}_{\mathcal{Y}} $, we obtain a map 
$$ g^{-1} \mathcal{I}_Z \rightarrow \mathcal{O}_{\mathcal{Y}}. $$
Now, let $ \mathcal{J} $ be the sheaf of ideals generated by the image of $ g^{-1} \mathcal{I}_Z $ in $ \mathcal{O}_{\mathcal{Y}} $. We have that $ \mathcal{J} $ is the ideal sheaf of $ g^{-1}(Z) $. But in our case $ g^{-1}(Z) = Z $, and therefore $ \mathcal{J} = \mathcal{I}_Z $.

The inclusion in \ref{equation 6.1} above induces an injective map
$$ g^*\mathcal{I}_Z \rightarrow g^* \mathcal{O}_{\mathcal{Y}} \cong \mathcal{O}_{\mathcal{Y}} $$
of $ \mathcal{O}_{\mathcal{Y}} $-modules, whose image is $ \mathcal{I}_Z = g^{-1} \mathcal{I}_Z \cdot \mathcal{O}_{\mathcal{Y}} $ (\cite{Hart}, II.7.12.2). Hence we obtain a homomorphism
\begin{equation}\label{equation 6.2}
u_Z : g^*\mathcal{I}_Z \rightarrow \mathcal{I}_Z 
\end{equation}
of $ \mathcal{O}_{\mathcal{Y}} $-modules, which is necessarily an isomorphism.

It remains to check that the maps $ g^*\mathcal{I}_Z \rightarrow \mathcal{I}_Z $ for various elements $ g \in G $ satisfy the compatibility conditions, but this is immediate since the maps are derived from the action on the surface.
\end{pf}




\begin{rmk}
The isomorphism $ u_Z : g^* \mathcal{I}_Z \rightarrow \mathcal{I}_Z $ constructed in the proof of Lemma \ref{G-sheaf} may be specified on stalks. Indeed, for any $ y \in \mathcal{Y} $, we have the map 
$$ g^{\sharp}_y : (g^{-1} \mathcal{O}_{\mathcal{Y}})_y \cong \mathcal{O}_{\mathcal{Y},g(y)} \rightarrow \mathcal{O}_{\mathcal{Y},y}, $$
which is an isomorphism of rings, and the image of $ \mathcal{I}_{Z,g(y)} \subset \mathcal{O}_{\mathcal{Y},g(y)} $ under this map is precisely $ \mathcal{I}_{Z,y} $. The inclusion $ \mathcal{I}_{Z,g(y)} \subset \mathcal{O}_{\mathcal{Y},g(y)} $ gives, after tensoring with $ \mathcal{O}_{\mathcal{Y},y} $ via $ g^{\sharp}_y $, an injective map
$$ \mathcal{I}_{Z,g(y)} \otimes_{\mathcal{O}_{\mathcal{Y},g(y)}} \mathcal{O}_{\mathcal{Y},y} \rightarrow \mathcal{O}_{\mathcal{Y},g(y)} \otimes_{\mathcal{O}_{\mathcal{Y},g(y)}} \mathcal{O}_{\mathcal{Y},y} \cong \mathcal{O}_{\mathcal{Y},y}. $$ 
Let $ f \in \mathcal{I}_{Z,g(y)} $ be a generator. Then $ f \otimes 1 $ generates $ \mathcal{I}_{Z,g(y)} \otimes_{\mathcal{O}_{\mathcal{Y},g(y)}} \mathcal{O}_{\mathcal{Y},y} $ as an $ \mathcal{O}_{\mathcal{Y},y} $-module, and the image of this element under the map above is $ g^{\sharp}_y(f) $, which generates $ \mathcal{I}_{Z,y} $. It follows that $ u_{Z,y}(f \otimes 1) = g^{\sharp}_y(f) $.
\end{rmk}

\subsection{Exact sequence}\label{6.6}
The construction above gives, for every $ g \in G $, a commutative diagram 
$$ \xymatrix{
g^* \mathcal{I}_Z \ar[d]_{u_{Z}} \ar@{^{(}->}[r] & g^* \mathcal{O}_{\mathcal{Y}} \ar[d]^{\cong}\\
\mathcal{I}_{Z} \ar@{^{(}->}[r] &  \mathcal{O}_{\mathcal{Y}}. } $$
This induces a commutative diagram
$$ \xymatrix{
0 \ar[r] & g^* \mathcal{I}_Z \ar[r] \ar[d]_{u_{Z}} & g^* \mathcal{O}_{\mathcal{Y}} \ar[r] \ar[d]_{\cong} & g^*(i_Z)_* \mathcal{O}_{Z} \ar[r] \ar[d]_{v_Z} & 0\\
0 \ar[r] & \mathcal{I}_Z \ar[r] &  \mathcal{O}_{\mathcal{Y}} \ar[r] & (i_Z)_* \mathcal{O}_{Z} \ar[r]  & 0,} $$
where the horizontal sequences are exact. Here we have identified $ (i_Z)_* \mathcal{O}_{Z} \cong \mathcal{O}_{\mathcal{Y}}/\mathcal{I}_Z $, and $ v_Z $ denotes the induced isomorphism on the cokernel.

\subsection{Relative version}\label{6.7}
\begin{prop}\label{prop. 6.5}
Let us keep the hypotheses and notation from Lemma \ref{lemma 6.1}. The sequence 
$$ 0  \rightarrow (i_C)_* \mathcal{L} \rightarrow (i_{Z})_* \mathcal{O}_{Z} \rightarrow (i_{Z'})_* \mathcal{O}_{Z'} \rightarrow 0 $$
is an exact sequence of $ G $-sheaves.
\end{prop}
\begin{pf}
Let $ \mathcal{I}_{Z} \subset \mathcal{I}_{Z'} \subset \mathcal{O}_{\mathcal{Y}} $ be the inclusions of the ideal sheaves. From Lemma \ref{G-sheaf}, it follows that these maps are maps of $G$-sheaves. The result now follows from the fact that the category of $G$-modules on $\mathcal{Y}$ is an abelian category (\cite{Kock}, Lemma 1.3). 
\end{pf}

\begin{rmk}
We would like to point out the following fact, that will be useful later: Let $ u_{Z'/Z} : g^*(\mathcal{I}_{Z'}/\mathcal{I}_{Z}) \rightarrow \mathcal{I}_{Z'}/\mathcal{I}_{Z} $ be the map that is the kernel of the diagram 
$$ 
\xymatrix{
g^*(\mathcal{O}_{\mathcal{Y}}/\mathcal{I}_{Z}) \ar[r] \ar[d]_{v_Z} & g^*(\mathcal{O}_{\mathcal{Y}}/\mathcal{I}_{Z'})  \ar[d]^{v_{Z'}} \\
\mathcal{O}_{\mathcal{Y}}/\mathcal{I}_{Z} \ar[r] & \mathcal{O}_{\mathcal{Y}}/\mathcal{I}_{Z'}.}
$$
Then $ u_{Z'/Z} $ can also be described as the \emph{cokernel} of the diagram
$$
\xymatrix{
g^*\mathcal{I}_{Z} \ar[r] \ar[d]_{u_Z} & g^*\mathcal{I}_{Z'} \ar[d]^{u_{Z'}}\\
\mathcal{I}_{Z} \ar[r] & \mathcal{I}_{Z'}.}
$$
\end{rmk}





\subsection{}\label{6.8}
Proposition \ref{prop. 6.5} implies that for any $g \in G$, we have a commutative diagram
\begin{equation}\label{equation 6.3}
 \xymatrix{
0 \ar[r] &  g^* (i_C)_* \mathcal{L} \ar[d]_{u_{Z'/Z}}  \ar[r] &  g^* (i_Z)_* \mathcal{O}_{Z} \ar[r] \ar[d]_{v_Z} & g^* (i_{Z'})_* \mathcal{O}_{Z'}  \ar[r] \ar[d]_{v_{Z'}} & 0\\
0 \ar[r] & (i_C)_* \mathcal{L} \ar[r] & (i_Z)_* \mathcal{O}_{Z} \ar[r] & (i_{Z'})_* \mathcal{O}_{Z'} \ar[r] & 0,} 
\end{equation}
where the horizontal lines are exact sequences. To avoid unnecessarily complicated notation, we will write $ u $, $ v $ and $ v' $ for the vertical maps above, if no confusion can arise.

Consider now the long exact sequence in cohomology induced by the lower exact sequence in Diagram \ref{equation 6.3}
\begin{equation} 
0 \rightarrow H^0(\mathcal{Y},(i_C)_* \mathcal{L}) \rightarrow H^0(\mathcal{Y},(i_Z)_* \mathcal{O}_{Z}) \rightarrow H^0(\mathcal{Y},(i_{Z'})_* \mathcal{O}_{Z'}) \rightarrow \ldots 
\end{equation}
$$ \ldots \rightarrow H^p(\mathcal{Y},(i_C)_* \mathcal{L})\rightarrow H^p(\mathcal{Y}, (i_Z)_* \mathcal{O}_{Z}) \rightarrow H^p(\mathcal{Y}, (i_{Z'})_* \mathcal{O}_{Z'}) \rightarrow  \ldots. $$ 

The maps $ u $, $ v $ and $ v' $, together with the automorphism $ g : \mathcal{Y} \rightarrow \mathcal{Y} $ induce, for every $ p \geq 0 $, automorphisms $ H^p(g,u) $, $ H^p(g,v) $ and $ H^p(g,v') $ that commute with the differentials in the long exact sequence. That is, we obtain, for every $ g \in G $, an \emph{automorphism} of this long exact sequence. 

Note that all the sheaves appearing in the exact sequence above are supported on the special fiber of $ \mathcal{Y} $. We will now explain how we can ``restrict'' the endomorphisms $ H^p(g,u) $, $ H^p(g,v) $ and $ H^p(g,v') $ to the support of the various sheaves.

\subsection{Restriction}\label{6.9}
Let $X$ and $Y$ be schemes, and let $ i : X \rightarrow Y $ be a closed immersion. For any (quasi-coherent) $ \mathcal{O}_X $-module $ \mathcal{F} $, the canonical adjunction map $ \alpha_{\mathcal{F}} : i^* i_* \mathcal{F} \rightarrow \mathcal{F} $ is an isomorphism, since $ i $ is a closed immersion. Let $ \mathcal{G} $ be a quasi-coherent $ \mathcal{O}_Y $-module. The canonical adjunction map $ \beta_{\mathcal{G}} : \mathcal{G} \rightarrow i_* i^* \mathcal{G} $ is not necessarily an isomorphism. However, in the case where $ \mathcal{G} $ has support on $X$, $ \beta_{\mathcal{G}} $ is indeed an isomorphism.

Recall also that $ i^* $ and $ i_* $ are adjoint functors, the maps $ \alpha $ and $ \beta $ induce a natural bijection
$$ \Hom_{\mathcal{O}_X}(i^* \mathcal{G}, \mathcal{F}) \cong \Hom_{\mathcal{O}_Y}(\mathcal{G}, i_* \mathcal{F}), $$
for any $ \mathcal{O}_X $-module $ \mathcal{F} $ and any $ \mathcal{O}_Y $-module $ \mathcal{G} $ on $Y$. 



Consider now the case where an automorphism $ g : Y \rightarrow Y $ is given, that restricts to an automorphism $ f = g|_X : X \rightarrow X $. Hence there is a commutative diagram
\begin{equation}\label{diagram 6.1} 
\xymatrix{ X \ar[r]^i \ar[d]_f & Y \ar[d]^g \\
              X \ar[r]^i & Y. } 
\end{equation}

We will now apply these considerations to cohomology. Let us first note that if $ \mathcal{F} $ is a quasi-coherent sheaf on $X$, then the push-forward $ i_*(\mathcal{F}) $ is a quasi-coherent sheaf on $ Y $, since $i$ is a closed immersion. 

\begin{prop}\label{prop. 6.6}
Let $ \mathcal{F} $ be a quasi-coherent sheaf on $X$, and let $ u : g^* i_* \mathcal{F} \rightarrow i_* \mathcal{F} $ be a homomorphism of $ \mathcal{O}_Y $-modules. Then there is induced, for every $ p \geq 0 $, a commutative diagram
$$ \xymatrix{
H^p(Y, i_* \mathcal{F}) \ar[rr]^{H^p(g,u)} \ar[d]_{H^p(i)} & & H^p(Y, i_* \mathcal{F}) \ar[d]^{H^p(i)} \\
H^p(X, i^* i_*\mathcal{F}) \ar[rr]^{H^p(f,i^*u)} & & H^p(X, i^* i_*\mathcal{F}),} 
$$
where the vertical arrows are isomorphisms.
\end{prop}
\begin{pf}
Since $ g \circ i = i \circ f $, we may identify the induced homomorphisms
\begin{equation}\label{equation 6.7} 
H^p(Y, i_*\mathcal{F}) \rightarrow H^p(X, (g \circ i)^* i_* \mathcal{F}) 
\end{equation}
and
\begin{equation}\label{equation 6.8} 
H^p(Y, i_*\mathcal{F}) \rightarrow H^p(X, (i \circ f)^* i_* \mathcal{F}). 
\end{equation}
However, homomorphism \ref{equation 6.7} factors as 
$$ H^p(Y, i_*\mathcal{F}) \rightarrow H^p(Y, g^* i_*\mathcal{F}) \rightarrow H^p(X, i^* g^* i_* \mathcal{F}) = H^p(X, f^* i^* i_* \mathcal{F}), $$
and homomorphism \ref{equation 6.8} factors as 
$$H^p(Y, i_*\mathcal{F}) \rightarrow H^p(X, i^* i_* \mathcal{F}) \rightarrow H^p(X, f^* i^* i_* \mathcal{F}) $$
(\cite{EGAIII}, Chap. 0, 12.1.3.5). Consequently, there is a commutative diagram
\begin{equation} 
\xymatrix{
H^p(Y, i_*\mathcal{F}) \ar[r]^{H^p(g)} \ar[d]_{H^p(i)} & H^p(Y, g^* i_*\mathcal{F}) \ar[d]^{H^p(i)} \\
H^p(X, i^* i_* \mathcal{F}) \ar[r]^{H^p(f)} & H^p(X, f^* i^* i_* \mathcal{F}). }
\end{equation}

If we pull back the map $ u : g^* i_* \mathcal{F} \rightarrow i_* \mathcal{F} $ with $ i $, we get the map 
$$ i^*(u) : i^* g^* i_* \mathcal{F} = f^* i^* i_* \mathcal{F} \rightarrow i^* i_* \mathcal{F}. $$ 
By functoriality, we obtain a commutative diagram

\begin{equation} 
\xymatrix{
H^p(Y, g^* i_*\mathcal{F}) \ar[r]^{H^p(u)} \ar[d]_{H^p(i)} & H^p(Y, i_*\mathcal{F}) \ar[d]^{H^p(i)} \\
H^p(X, f^* i^* i_* \mathcal{F}) \ar[r]^{H^p(i^*u)} & H^p(X, i^* i_* \mathcal{F}). }
\end{equation}

To complete the proof, we need to check that the maps $ H^p(i) $ are isomorphisms. In order to see this, we first note that the adjunction map $ \alpha_{\mathcal{G}} : \mathcal{G} \rightarrow i_* i^* \mathcal{G} $, where we put $ \mathcal{G} := i_* \mathcal{F} $, is an isomorphism. For any open $ V \subset Y $, there is induced a map 
$$ \Gamma(V, \mathcal{G}) \rightarrow \Gamma(V, i_* i^* \mathcal{G}) = \Gamma(i^{-1}(V), i^* \mathcal{G}), $$
that is also an isomorphism. Choose now an open affine cover $ \mathcal{V} = (V_j)_{j \in J} $ of $Y$. As $i$ is a closed immersion, the inverse image $ \mathcal{U} := i^{-1}(\mathcal{V}) $ is an open affine cover of $X$. The map on sections above induces a map of Cech-complexes
$$ C^{\bullet}(\mathcal{V},\mathcal{G}) \rightarrow C^{\bullet}(\mathcal{U},i^*\mathcal{G}), $$
that is also an isomorphism. Passing on to cohomology, we therefore get isomorphisms
$$ H^p(\mathcal{V},\mathcal{G}) \rightarrow H^p(\mathcal{U},i^* \mathcal{G}), $$
for every $ p \geq 0 $.

The natural maps 
$$ H^p(\mathcal{V},\mathcal{G}) \rightarrow H^p(Y,\mathcal{G}) $$ 
and
$$ H^p(\mathcal{U},i^* \mathcal{G}) \rightarrow H^p(X,i^* \mathcal{G}) $$
are isomorphisms, by \cite{Hart}, Theorem III 4.5. By \cite{EGAIII}, Chap. 0, 12.1.4.2, the diagram
$$ \xymatrix{
H^p(\mathcal{V},\mathcal{G}) \ar[r] \ar[d] & H^p(\mathcal{U},i^* \mathcal{G}) \ar[d] \\
 H^p(Y,\mathcal{G}) \ar[r] & H^p(X,i^* \mathcal{G})}
$$
commutes. Therefore, we can conclude that 
$$ H^p(Y,\mathcal{G}) \rightarrow H^p(X,i^* \mathcal{G}) $$
is an isomorphism.
\end{pf}

\subsection{}\label{6.11}
We are now going to apply the technical results in Section \ref{6.9} to the setup in Section \ref{6.8}. To begin with, we note that the closed immersions $ i_C $, $ i_Z $ and $ i_{Z'} $ induce isomorphisms
$$ H^p(\mathcal{Y}, (i_C)_* \mathcal{L}) \cong H^p(C,\mathcal{L}), $$ 
$$ H^p(\mathcal{Y}, (i_Z)_* \mathcal{O}_Z) \cong H^p(Z, \mathcal{O}_Z) $$
and
$$ H^p(\mathcal{Y}, (i_{Z'})_* \mathcal{O}_{Z'}) \cong H^p(Z', \mathcal{O}_{Z'}), $$
for all $ p \geq 0 $. Here we have identified $ \mathcal{L} $ with $ (i_C)^* (i_C)_* \mathcal{L} $, and likewise for $ \mathcal{O}_Z $ and $ \mathcal{O}_{Z'} $.

Since $ C $, $ Z $ and $ Z' $ are projective curves over $k$, and since $ \mathcal{L} $, $ \mathcal{O}_Z $ and $ \mathcal{O}_{Z'} $ are coherent sheaves on the respective curves, the cohomology groups above are finite dimensional $k$-vector spaces, and nonzero only for $ p = 0 $ and $ p = 1 $.
So the long exact sequence in cohomology associated to the short exact sequence in Section \ref{6.8} is simply
\begin{equation}\label{sequence 6.7} 
0 \rightarrow H^0(C,\mathcal{L}) \rightarrow H^0(Z, \mathcal{O}_Z) \rightarrow H^0(Z', \mathcal{O}_{Z'}) \rightarrow 
\end{equation}
$$ H^1(C,\mathcal{L}) \rightarrow H^1(Z, \mathcal{O}_Z) \rightarrow H^1(Z', \mathcal{O}_{Z'}) \rightarrow 0. $$

Denote by $g_C$ the restriction of $g$ to $C$. By Proposition \ref{prop. 6.6}, restriction to $C$ gives a commutative diagram
$$ \xymatrix{ H^p(\mathcal{Y},(i_C)_* \mathcal{L}) \ar[d]_{H^p(g,u)} \ar[rr]^{\cong} & & H^p(C,\mathcal{L}) \ar[d]^{H^p(g_C, i_C^*u)} \\
H^p(\mathcal{Y},(i_C)_* \mathcal{L}) \ar[rr]^{\cong} & & H^p(C,\mathcal{L}). } $$
Also, we get similar diagrams for $ H^p(g_Z, i_Z^* v) $ and $ H^p(g_{Z'}, i_{Z'}^* v') $. 

Having made these identifications, we see that the automorphisms $ H^p(g_C, i_C^* u)$, $ H^p(g_Z, i_Z^* v) $ and $ H^p(g_{Z'}, i_{Z'}^* v') $, for $ p = 0, 1 $, fit together to give an automorphism of Sequence \ref{sequence 6.7} above. In particular, \ref{sequence 6.7} is an exact sequence of $ k[G] $-modules.



\begin{lemma}\label{lemma 6.7}
Let us keep the notation and constructions from the discussion above. The following equality holds in $W(k)$:
$$ \sum_{p=0}^1 (-1)^p ~\Tr_{\beta} (H^p(g_Z, i_Z^* v)) = $$
$$ \sum_{p=0}^1 (-1)^p ~\Tr_{\beta} (H^p(g_C,i_C^* u)) + \sum_{p=0}^1 (-1)^p ~\Tr_{\beta} (H^p(g_{Z'}, i_{Z'}^* v')). $$
\end{lemma}
\begin{pf}
The long exact sequence \ref{sequence 6.7} is an exact sequence of $k[G]$-modules. Therefore, for any $ g \in G $, the alternating sum of the Brauer traces equals zero. 
\end{pf}

\subsection{}
We will now generalize the formula in Lemma \ref{lemma 6.7}. Let us write $ \mathcal{Y}_k = \sum_{\alpha} m_{\alpha} C_{\alpha} $, where $ \alpha \in \mathcal{A} $, and put $ m = \sum_{\alpha} m_{\alpha} $. Fix a complete filtration
$$ 0 < Z_m < \ldots < Z_j < \ldots < Z_2 < Z_1 = \mathcal{Y}_k, $$ 
where $ Z_j - Z_{j+1} = C_j $ for some $ C_j \in \{ C_{\alpha} \}_{\alpha \in \mathcal{A} } $, for each $ j \in \{ 1, \ldots, m-1 \} $. At each step of this filtration, Lemma \ref{lemma 6.1} asserts that there is a short exact sequence
$$ 0 \rightarrow (i_{C_j})_* \mathcal{L}_j \rightarrow (i_{Z_j})_* \mathcal{O}_{Z_j} \rightarrow (i_{Z_{j+1}})_* \mathcal{O}_{Z_{j+1}} \rightarrow 0. $$
Note in particular that $ Z_m = C_m $, for some $ C_m \in \{ C_{\alpha} \}_{\alpha \in \mathcal{A} } $. So we write $ \mathcal{O}_{Z_m} = \mathcal{L}_m $, for the sake of coherence.

For every $ j \in \{ 1, \ldots , m \} $, we have isomorphisms  
$$ u_j : g^* (i_{C_j})_* \mathcal{L}_j \rightarrow (i_{C_j})_* \mathcal{L}_j, $$
and 
$$ v_j : g^* (i_{Z_j})_* \mathcal{O}_{Z_j} \rightarrow (i_{Z_j})_* \mathcal{O}_{Z_j}. $$

\begin{prop}\label{prop. 6.6}
The alternating sum of the Brauer traces of the automorphisms induced by $ g \in G $ on $ H^p(\mathcal{Y}_k, \mathcal{O}_{\mathcal{Y}_k}) $ can be computed by the formula
$$ \sum_{p=0}^1 (-1)^p ~\Tr_{\beta} (H^p(g_{\mathcal{Y}_k})) = \sum_{j=1}^m (\sum_{p=0}^1 (-1)^p ~\Tr_{\beta} (H^p(g_{C_j}, i_{C_j}^* u_j))). $$
\end{prop}
\begin{pf}
Let us first note that $ (i_{Z_j})^*v_j $ can be identified with the canonical isomorphism $ g_{Z_j}^* \mathcal{O}_{Z_j} \cong \mathcal{O}_{Z_j} $, for every $j$. So the $v_j$ will be dropped from the notation.

For $ j = 1 $, we have that $ Z_1 = \mathcal{Y}_k $, and so Lemma \ref{lemma 6.7} gives that
$$ \sum_{p=0}^1 (-1)^p ~\Tr_{\beta} (H^p(g_{\mathcal{Y}_k})) = $$
$$ \sum_{p=0}^1 (-1)^p ~\Tr_{\beta} (H^p(g_{C_1}, i_{C_1}^* u_1)) + \sum_{p=0}^1 (-1)^p \Tr_{\beta} ~(H^p(g_{Z_2})). $$

However, for any $ j \in \{ 1, \ldots , m-1 \} $, Lemma \ref{lemma 6.7} again gives that
$$ \sum_{p=0}^1 (-1)^p ~\Tr_{\beta} (H^p(g_{Z_j})) = $$
$$ \sum_{p=0}^1 (-1)^p ~\Tr_{\beta} (H^p(g_{C_j}, i_{C_j}^* u_j)) + \sum_{p=0}^1 (-1)^p ~\Tr_{\beta} (H^p(g_{Z_{j+1}})). $$
So the result follows by induction.
\end{pf}

\begin{rmk}
The importance of Proposition \ref{prop. 6.6} is that it reduces the problem of computing 
$$ \sum_{p=0}^1 (-1)^p ~\Tr_{\beta} (H^p(g_{\mathcal{Y}_k})) $$  
to instead computing the contributions 
$$ \sum_{p=0}^1 (-1)^p ~\Tr_{\beta} (H^p(g_{C_j}, i_{C_j}^* u_j)), $$
for certain invertible sheaves $ \mathcal{L}_j $, supported on the smooth irreducible components $ C_j $ of $ \mathcal{Y}_k $. For such computations, we may apply Proposition \ref{prop. 5.5} and Proposition \ref{prop. 5.6}. In order to apply these formulas, we need to understand the action on each irreducible component $C_j$ of $ \mathcal{Y}_k $. Furthermore, we need to understand the action on the fiber of the invertible sheaf $ \mathcal{L}_j $ at any fixed point on $C_j$.
\end{rmk}


\subsection{Fibers at fixed points and cotangent spaces}\label{7.1}
Let $ g \in G $, and consider the corresponding automorphism
$$ g : \mathcal{Y} \rightarrow \mathcal{Y}. $$
Let $ y \in \mathcal{Y}_k $ be an intersection point of two irreducible components $ C $ and $ C' $ of $ \mathcal{Y}_k $. Then $ y $ is a fixed point for $ g $, so there is an induced automorphism
$$ g^{\sharp}_y : \mathcal{O}_{\mathcal{Y},y} \rightarrow \mathcal{O}_{\mathcal{Y},y}. $$
As $ g^{\sharp}_y $ is a local homomorphism, it induces an automorphism of the cotangent space at $y$:
$$ dg(y) : m_y/m_y^2 \rightarrow m_y/m_y^2. $$

Let $ f $ and $ f' $ be local equations for $ C $ and $ C' $ at $y$. That is, we have 
$ \mathcal{I}_{C,y} = (f) \subset \mathcal{O}_{\mathcal{Y},y} $ and similarly $ \mathcal{I}_{C',y} = (f') \subset \mathcal{O}_{\mathcal{Y},y} $. Since $ \mathcal{Y} $ has strict normal crossings, it follows that  $ m_y = \mathcal{I}_{C,y} + \mathcal{I}_{C',y} = (f, f') $. In particular, the images of $f$ and $f'$ form a basis for $ m_y/m_y^2 $.

Since $C$ and $C'$ are stable under the $ G $-action, we have that $ \mathcal{I}_{C,y} $ and $ \mathcal{I}_{C',y} $ map to themselves via $ g^{\sharp}_y $. Hence, the images of $f$ and $f'$ are \emph{eigenvectors} for $ dg(y) $. Let $ \lambda $ and $ \lambda' $ be the corresponding eigenvalues. 

\subsection{}\label{section 7.2}

Consider the inclusion $ (f) = \mathcal{I}_{C,y} \subset m_y $. By tensoring with $ \mathcal{O}_{\mathcal{Y},y}/m_y $, we get a map
$$ \mathcal{I}_{C}(y) = \mathcal{I}_{C,y} \otimes_{\mathcal{O}_{\mathcal{Y},y}} \mathcal{O}_{\mathcal{Y},y}/m_y \rightarrow m_y \otimes_{\mathcal{O}_{\mathcal{Y},y}} \mathcal{O}_{\mathcal{Y},y}/m_y = m_y/m_y^2, $$
mapping $ f \otimes 1 $ to the image of $f$ in $ m_y/m_y^2 $. This is a non-degenerate map of $ k $-vector spaces, so we may identify 
$$ \mathcal{I}_{C}(y) = <f> \subset m_y/m_y^2. $$
Similarly, we may identify
$$ \mathcal{I}_{C'}(y) = <f'> \subset m_y/m_y^2. $$

\begin{lemma}\label{Lemma 7.1}
Let $ u_C : g^* \mathcal{I}_C \rightarrow \mathcal{I}_C $ be the map constructed in Lemma \ref{G-sheaf}. Then the eigenvalue of the induced map $ u_C(y) : (g^* \mathcal{I}_C)(y) \rightarrow \mathcal{I}_C(y) $ at the fiber in the fixed point $y$ is precisely the eigenvalue $\lambda$ of $f$ under the cotangent map $dg(y)$ at $y$.

Similarly, the eigenvalue of the map $ u_{C'}(y) : (g^* \mathcal{I}_{C'})(y) \rightarrow \mathcal{I}_{C'}(y) $ at the fiber in $y$ is precisely the eigenvalue $ \lambda'$ of $f'$ under $dg(y)$.
\end{lemma}
\begin{pf}
Let $ g^{\sharp}_y : \mathcal{O}_{\mathcal{Y},y} \rightarrow \mathcal{O}_{\mathcal{Y},y} $ be the map of local rings induced from $ g : \mathcal{Y} \rightarrow \mathcal{Y} $. We can describe the stalk of $ u_C $ at $y$ in the following way: Tensoring the inclusion $ \mathcal{I}_{C,y} \subset \mathcal{O}_{\mathcal{Y},y} $ with $ g^{\sharp}_y $ gives 
$$ \mathcal{I}_{C,y} \otimes^g \mathcal{O}_{\mathcal{Y},y} \rightarrow \mathcal{O}_{\mathcal{Y},y} \otimes^g \mathcal{O}_{\mathcal{Y},y}, $$ 
where $ M \otimes^g \mathcal{O}_{\mathcal{Y},y} $ denotes, for any $ \mathcal{O}_{\mathcal{Y},y} $-module $M$, tensorization with $ \mathcal{O}_{\mathcal{Y},y} $ via the homomorphism $ g^{\sharp}_y $.

By composing with the canonical isomorphism $ \mathcal{O}_{\mathcal{Y},y} \otimes^g \mathcal{O}_{\mathcal{Y},y} \cong \mathcal{O}_{\mathcal{Y},y} $, we get an injective map
$$ \mathcal{I}_{C,y} \otimes^g \mathcal{O}_{\mathcal{Y},y} \rightarrow \mathcal{O}_{\mathcal{Y},y} $$
with image $ \mathcal{I}_{C,y} $. The induced map 
$$ (g^* \mathcal{I}_C)_y = \mathcal{I}_{C,y} \otimes^g \mathcal{O}_{\mathcal{Y},y} \rightarrow \mathcal{I}_{C,y} $$
is the stalk $ (u_C)_y $. Let $ a \in \mathcal{I}_{C,y} $ be any element. Then we have that 
$$ (u_C)_y(a \otimes 1) = g^{\sharp}_y(a) \in \mathcal{I}_{C,y}. $$  

We will now consider the fiber $ (u_C)(y) $ of $ (u_C)_y $ at $y$. Let us first point out that $ g^{\sharp}_y $ induces the identity on the residue field of $y$. That is, we have a commutative diagram
$$ \xymatrix{
\mathcal{O}_{\mathcal{Y},y} \ar[r]^{g^{\sharp}_y} \ar[d] & \mathcal{O}_{\mathcal{Y},y} \ar[d] \\
\mathcal{O}_{\mathcal{Y},y}/m_y \ar[r]^{=} & \mathcal{O}_{\mathcal{Y},y}/m_y.}
$$
It therefore follows that 
$$ (g^* \mathcal{I}_{C})(y) = (g^* \mathcal{I}_C)_y \otimes_{\mathcal{O}_{\mathcal{Y},y}} \mathcal{O}_{\mathcal{Y},y}/m_y  = $$
$$ \mathcal{I}_{C,y} \otimes^g \mathcal{O}_{\mathcal{Y},y} \otimes_{\mathcal{O}_{\mathcal{Y},y}} \mathcal{O}_{\mathcal{Y},y}/m_y = \mathcal{I}_{C,y} \otimes_{\mathcal{O}_{\mathcal{Y},y}} \mathcal{O}_{\mathcal{Y},y}/m_y = \mathcal{I}_{C}(y). $$
With this identification, we see that the map $ (u_C)(y) $ on the fiber is 
$$ (u_C)(y) : \mathcal{I}_{C,y} \otimes_{\mathcal{O}_{\mathcal{Y},y}} \mathcal{O}_{\mathcal{Y},y}/m_y \rightarrow \mathcal{I}_{C,y} \otimes_{\mathcal{O}_{\mathcal{Y},y}} \mathcal{O}_{\mathcal{Y},y}/m_y, $$ where, for any $ a \in \mathcal{I}_{C,y} $, we have that $ (u_C)(y)(a \otimes 1) = g^{\sharp}_y(a) \otimes 1 $. 

If we replace $ \mathcal{I}_{C,y} $ with the maximal ideal, we can argue in exactly the same way, and get a map
$$ m_y \otimes_{\mathcal{O}_{\mathcal{Y},y}} \mathcal{O}_{\mathcal{Y},y}/m_y \rightarrow m_y \otimes_{\mathcal{O}_{\mathcal{Y},y}} \mathcal{O}_{\mathcal{Y},y}/m_y, $$
given by $ a \otimes 1 \mapsto g^{\sharp}_y(a) \otimes 1 $, for any $ a \in m_y $. Hence this is precisely the cotangent map $ dg(y) : m_y/m_y^2 \rightarrow m_y/m_y^2 $.



The inclusion $ \mathcal{I}_{C,y} \subset m_y $ now gives a commutative diagram 

$$ \xymatrix{
\mathcal{I}_C(y) \ar[d]_{u_{C}(y)} \ar[r] & m_y/m_y^2 \ar[d]^{dg(y)}\\
\mathcal{I}_{C}(y) \ar[r] &  m_y/m_y^2. } $$
Using the identifications above, we see that the eigenvalue of $ u_{C}(y) $ is precisely the eigenvalue $ \lambda $ of $f$ under the cotangent map $ dg(y) $.

In a similar way, we can show that the eigenvalue of $ u_{C'}(y) $ is precisely the eigenvalue $ \lambda' $ of $f'$ under the cotangent map $ dg(y) $. 
\end{pf}

\subsection{}
Let $ 0 \leq Z \leq \mathcal{Y}_k $ be an effective divisor on the form 
$$ Z = a C + a' C' + Z_0, $$
where the effective divisor $ Z_0 $ does not contain $ C $ or $C'$. Then the ideal sheaf of $Z$ can be written as
$$ \mathcal{I}_Z = \mathcal{I}_C^{\otimes a} \otimes \mathcal{I}_{C'}^{\otimes a'} \otimes \mathcal{I}_0, $$
where $ \mathcal{I}_0 $ is the ideal sheaf of $Z_0$.

\begin{lemma}\label{lemma 7.2}
Let $ u_Z : g^* \mathcal{I}_Z \rightarrow \mathcal{I}_Z $ be the map constructed in Lemma \ref{G-sheaf}. Then the eigenvalue of $ u_Z $ at the fiber in $y$ is given by
$$ \Tr(u_Z(y)) = \lambda^a \lambda'^{a'}. $$ 
\end{lemma}
\begin{pf}
One proves this in a similar way as Lemma \ref{Lemma 7.1}.
\end{pf}

\subsection{}
Finally, we consider divisors $ 0 \leq Z' < Z \leq \mathcal{Y}_k $, where $ Z - Z' = C $ for some irreducible component $ C $ of $ \mathcal{Y}_k $.

\begin{lemma}\label{lemma 7.3}
Let $ u_{Z'/Z} : g^* (i_C)_* \mathcal{L} \rightarrow (i_C)_* \mathcal{L} $ be the map constructed in Section \ref{6.8}. Then the eigenvalue of the pullback $ i_C^* (u_{Z'/Z}) $ in the fiber at the fixed point $y$ is 
$$ \Tr(i_C^* (u_{Z'/Z})(y)) = \Tr(u_{Z'}(y)), $$
where $ u_{Z'} : g^* \mathcal{I}_{Z'} \rightarrow \mathcal{I}_{Z'} $ is the corresponding map for $ \mathcal{I}_{Z'} $.
\end{lemma}
\begin{pf}
In the proof of Proposition \ref{prop. 6.5} we saw that $ u_{Z'/Z} $ fitted into a commutative diagram 
$$ \xymatrix{
g^*\mathcal{I}_{Z'} \ar[r]^-{g^* \pi} \ar[d]_{u_{Z'}} & g^*(\mathcal{I}_{Z'}/\mathcal{I}_{Z}) = g^* (i_C)_* \mathcal{L} \ar[d]^{u_{Z'/Z}} \\
\mathcal{I}_{Z'} \ar[r]^-{\pi} & \mathcal{I}_{Z'}/\mathcal{I}_{Z} = (i_C)_* \mathcal{L},} $$
where $ \pi $ is the quotient surjection. Let $ i_y : \{y\} \hookrightarrow \mathcal{Y} $ be the inclusion. Since $y$ is a fixed point, we have that $ g \circ i_y = i_y $, and therefore $ i_y^* g^* = i_y^* $. It follows that $ i_y^* \pi = i_y^* g^* \pi $. Moreover, $i_y$ factors via the inclusion $ i_C $, so pulling back with $i_y $, we get the diagram 
$$ \xymatrix{
\mathcal{L}(y) \ar[r]^-{i_y^* \pi} \ar[d]_{u_{Z'}(y)} & \mathcal{L}(y) \ar[d]^{u_{Z'/Z}(y)} \\
\mathcal{L}(y) \ar[r]^-{i_y^* \pi} & \mathcal{L}(y), } $$
and the result follows.



\end{pf}

\section{Resolution fibers of tame cyclic quotient singularities}
In this section we study certain properties of tame cyclic quotient singularities, which will be important in the following sections. As a motivation for the discussion below, let us consider the following set-up: Let $ \mathcal{X} $ be an SNC-model, and let $ S'/S $ be a tame extension of degree $n$ that is prime to the multiplicities of the irreducible components of $ \mathcal{X}_k $. The normalization $ \mathcal{X}' $ of $ \mathcal{X}_{S'} $ has tame cyclic quotient singularities, depending on the combinatorial structure of $ \mathcal{X}_k $ and on $n$. We are interested in how the exceptional locus of the desingularization $ \rho : \mathcal{Y} \rightarrow \mathcal{X}' $ varies as the degree $n$ grows. 

This is actually a local problem, so we will start by formalizing the situation, and identify singularities with certain parameters $m_1$, $m_2$ and $n$. We conclude in Proposition \ref{prop. 8.9}, that if $n$ is sufficiently large, the combinatorial structure of the exceptional locus in the minimal desingularization belongs, modulo chains of $(-2)$-curves, to only a finite number of types, depending only on $m_1$ and $m_2$. 

\subsection{Construction}\label{8.1}
Let $ m_1 $ and $ m_2 $ be positive integers, and let $ n $ be a positive integer that is not divisible by $p$, and that is relatively prime to $ \Lcm(m_1,m_2) $. The integers $m_1$, $m_2$ and $n$ define a tame cyclic quotient singularity in the following way:

Consider the ring $ A = R[[u_1,u_2]]/( \pi - u_1^{m_1} u_2^{m_2} ) $, and let $ R \rightarrow R' $ be a tamely ramified extension of degree $n$. Let $ B = A \otimes_R R' = R'[[u_1,u_2]]/( \pi'^n - u_1^{m_1} u_2^{m_2} ) $, and let $ \widetilde{B} $ be the normalization of $B$. Then $ Z = \Spec(\widetilde{B}) $ is a tame cyclic quotient singularity.

Note also that $B$ can be equipped with the obvious $ G = \boldsymbol{\mu}_n $-action $ [ \xi ](\pi') = \xi \pi' $, for every $ \xi \in \boldsymbol{\mu}_n $.

\begin{dfn}\label{Def. 8.1}
A tame cyclic quotient singularity arising as in \ref{8.1} above will be denoted as \emph{the} singularity $(m_1,m_2,n)$. That is, we identify the singularity with its parameters $ m_1 $, $ m_2 $ and $ n $. 
\end{dfn}

In this section, we will make a detailed study of the behaviour of the minimal resolution of a singularity with parameters $ (m_1,m_2,n) $, where we keep $m_1$ and $m_2$ fixed, but where $n$ varies. 


Let $ \rho : \widetilde{Z} \rightarrow Z $ be the minimal desingularization. We saw in Section \ref{3.11} that the action of the Galois group $G$ on the special fiber $ \widetilde{Z}_k $ was completely determined by the parameters $ (m_1,m_2,n) $ of the singularity.

\subsection{Data associated to the singularity}
Let us consider a singularity with parameters $ (m_1,m_2,n) $. Let $ m := \Gcd(m_1,m_2) $ and $ M := \Lcm(m_1,m_2) $, and let furthermore $r$ be the unique integer with $ 0<r<n $ such that $ m_1 + rm_2 = 0 $ modulo $n$. Write $ \frac{n}{r} = [b_1, \ldots , b_l, \ldots , b_L]_{JH} $ for the Jung-Hirzebruch expansion. As $L$ is the length of the expansion, there are $L$ exceptional components $ C_1, \ldots, C_L $, with self intersection numbers $ C_l^2 = - b_l $, for all $ l \in \{1, \ldots, L \} $. We let $ \mu_l $ denote the multiplicity of $C_l $, for all $l$. There are two series of numerical equations associated to the singularity. We have
\begin{equation}\label{equation 8.1} 
r_{l-1} = b_{l+1} r_l - r_{l+1}, 
\end{equation}
for $ 0 \leq l \leq L-1 $, where we put $ r_{-1} = n $ and $ r_0 = r $. Furthermore, we have
 \begin{equation}\label{equation 8.2} 
\mu_{l+1} = b_l \mu_l - \mu_{l-1}, 
\end{equation}
which is valid for $ 1 \leq l \leq L $. Here we define $ \mu_0 = m_2 $ and $ \mu_{L+1} = m_1 $. Note also that $ m_1 + r m_2 = n \mu_1 $ (see \cite{CED}, Corollary 2.4.3).

Let us also define $ C_0 $ to be the (formal) branch with multiplicity $ \mu_0 $, and $ C_{L+1} $ the (formal) branch with multiplicity $ \mu_{L+1} $. Let $ y_l $ be the node in the special fiber which is the intersection point of $ C_{l+1} $ and $ C_l $. Then we saw in Section \ref{3.11} that for any $ 0 \leq l \leq L $, we have
$$ \widehat{\mathcal{O}}_{\widetilde{Z},y_l} = R'[[z_l,w_l]]/(z_l^{\mu_{l+1}} w_l^{\mu_l} - \pi'). $$
By Proposition \ref{prop. 3.3}, the Galois group $G$ acts on this ring by $ [ \xi ]( \pi' ) = \xi \pi' $, $ [ \xi ](z_l) = \xi^{ \alpha_1 r_{l-1} } z_l $ and $ [ \xi ](w_l) = \xi^{ - \alpha_1 r_l } w_l $, where $ \alpha_1 $ is an inverse to $ m_1 $ modulo $n$. 

Since we have that $ I_{C_l} = (w_l) $ and $ I_{C_{l+1}} = (z_l) $, it follows that $ \widehat{\mathcal{O}}_{C_l, y_{l-1}} = k[[w_{l-1}]] $, and that $ \widehat{\mathcal{O}}_{C_l, y_{l}} = k[[z_l]] $. The cotangent space of $C_l$ at $y_{l-1}$ is generated by $ w_{l-1} $, and the eigenvalue is therefore $ \xi^{ - \alpha_1 r_{l-1} } $, and at $ y_{l} $ the cotangent space is generated by $ z_l $ with eigenvalue $ \xi^{ \alpha_1 r_{l-1} } $.

\subsection{Some general properties of the minimal resolution of quotient singularities}
The lemma below lists some properties of the exceptional locus of the minimal resolution of a tame cyclic quotient singularity. These properties will be used numerous times in the rest of this paper.

\begin{lemma}\label{lemma 8.1}
Let $ m_1 $ and $ m_2 $ be positive integers, and let $n$ be a positive integer not divisible by $p$ such that $ \Gcd(n,m_1) = \Gcd(n,m_2) = 1 $. Denote by $ (m_1,m_2,n) $ the associated singularity. Then the following properties hold:
\begin{enumerate}
\item If $d$ divides both $ \mu_k $ and $ \mu_{k+1} $ for some $ 0 \leq k \leq L $, then $d$ divides $ \mu_l $ for all $ 0 \leq l \leq L+1 $.
\item Let $ m = \Gcd(m_1,m_2) $. Then $m$ divides $ \mu_l $ for all $ 0 \leq l \leq L+1 $.
\item Let $ l $ be an integer such that $ 1 \leq l \leq L $. Then the pairs of inequalities $ \mu_{l-1} < \mu_l $ and $ \mu_{l+1} \leq \mu_l $, or $ \mu_{l-1} \leq \mu_l $ and $ \mu_{l+1} < \mu_l $, can not occur. 
\end{enumerate}
\end{lemma}
\begin{pf}
(i) Equation \ref{equation 8.2} gives that $ \mu_{k-1} = \mu_k b_k - \mu_{k+1} $. So if $d$ divides $ \mu_k $ and $ \mu_{k+1} $, it will also divide $ \mu_{k-1} $. Continuing inductively, we get that $d$ divides $ \mu_l $ for all $ 0 \leq l \leq k+1 $. The same argument shows that we also get that $d$ divides $ \mu_l $ for all $ k+1 \leq l \leq L+1 $ if we do induction for increasing $l$ instead.

(ii) Recall that we have the equation $ m_1 + r_0 m_2 = n \mu_1 $. So $m$ divides $ n \mu_1 $. But $n$ is relatively prime to both $m_1$ and $m_2$, and in particular to $m$, so $m$ must therefore divide $\mu_1$. But then we have that $m$ divides both $\mu_0$ and $\mu_1$, so by (i), we can conclude that $m$ divides $ \mu_l $ for all $ 0 \leq l \leq L+1 $.

(iii) Assume that $ \mu_{l-1} < \mu_l $ and $ \mu_{l+1} \leq \mu_l $ for some $l$. Using Equation \ref{equation 8.2}, we get that 
$$ 0 < b_l = (\mu_{l-1} + \mu_{l+1})/\mu_l < (\mu_l + \mu_l) /\mu_l = 2, $$
which implies that $ C_l^2 = - 1 $. Hence we get that $ C_l $ is a $(-1)$-curve, which contradicts the minimality of the resolution.  
\end{pf}

\begin{cor}\label{m_1=m_2}
If $ m_1 = m_2 $, then $ \mu_l = m_1 = m_2 $ for all $ l $.
\end{cor}

Let us also give a result about the action on the components in the exceptional locus.

\begin{lemma}\label{Lemma 8.4} Let $ \sigma = (m_1,m_2,n) $ be a singularity. Then the following property holds. Let $ \xi \in G $ be a primitive $n$-th root of unity, and denote by $g_{\xi}$ the induced automorphism on the minimal desingularization $ \widetilde{Z} $. Then the restriction $ g_{\xi}|_{C_l} $ is a non-trivial automorphism of each exceptional curve $C_l$, for all $ l =  1, \ldots , L $. Furthermore, the fixed points of $ g_{\xi}|_{C_l} $ are exactly the two points where $ C_l $ meets the rest of the special fiber.

\end{lemma}
\begin{pf}
Let $C_l$ be any of the exceptional components. Then we saw in Proposition \ref{prop. 3.3} that we may cover $C_l$ with the affine charts $ \Spec(k[w_{l-1}]) $ and $ \Spec(k[z_l]) $, where we glue along $ z_l = 1/w_{l-1} $. In these coordinates, the action is $ [\xi](z_l) = \xi^{\alpha_1 r_{l-1}} z_l $ and $ [\xi](w_{l-1}) = \xi^{ - \alpha_1 r_{l-1}} w_{l-1} $. Since $ \xi $ is a primitive root, it suffices to show that $ \alpha_1 r_{l} \not\equiv_n 0 $ for any $ l \in \{ 0, \ldots, L - 1 \}  $. As $\alpha_1$ is invertible modulo $n$, this is the same as showing that $ r_{l} \not\equiv_n 0 $. But the $r_l$ satisfy the inequalities 
$$ n = r_{-1} > r_0 > \ldots > r_{L-1} = 1, $$
which proves the first statement.



The last statement follows from the explicit description of the action on the affine charts of the exceptional components.
\end{pf}


\subsection{}
We will now study the minimal resolution of a singularity $ (m_1,m_2,n) $ in more detail. In particular, we are interested in what happens if we keep $m_1$ and $ m_2 $ fixed, but let $n$ grow to infinity. We shall see that, modulo a certain equivalence relation, there are finitely many possibilities, corresponding to the elements in $ (\mathbb{Z}/M)^* $.

\begin{lemma}\label{lemma 8.4}
Let $m_1, m_2$ be positive integers, and let $ n $ be a positive integer not divisible by $p$ such that $ \Gcd(n,M) = 1 $. Assume that the equation $ m_1 + r m_2 = n \mu $ holds, where $ r $ is an integer such that $ 0 < r < n $, and where $ \mu $ is some positive integer. Recall that $ \mu $ is the multiplicity $ \mu_1 $ of the first exceptional component in the resolution of the singularity $(m_1,m_2,n)$. Then the following properties hold:
\begin{enumerate}
\item Assume that $ \mu < m_2 $. Let $ n' = n + k M $ and let $ r' = r + k \mu \frac{M}{m_2} $, for any $ k \in \mathbb{N} $. Then we have that $ \Gcd(n',M) = 1 $, that $ 0 < r' < n' $, and that the equation $ m_1 + r' m_2 = n' \mu $ holds.
\item If $ m_2 > m_1 $, we automatically have that $ \mu < m_2 $. 
\end{enumerate}
\end{lemma}
\begin{pf}
(i) Assume that $t$ is a common divisor of $ n' $ and $M$. Since $ n = n' - kM $, it follows that $t$ also divides $n$, and hence $t=1$. In particular, $m_1$ and $m_2$ are invertible modulo $n'$. Furthermore, we see that
$$ n' \mu = (n + kM) \mu = n \mu + k \mu \frac{M}{m_2} m_2 = m_1 + (r + k \mu \frac{M}{m_2}) m_2 = m_1 + r' m_2.$$   
It is clear from this equation that $r'$ is invertible modulo $n'$. We now compute that
$$ 0 < r' = r + k \mu \frac{M}{m_2} < n + k m_2 \frac{M}{m_2} = n + k M = n', $$
so it follows that $ 0 < r' < n'  $.

(ii) Assume that $ \mu \geq m_2 $, so that $ n m_2 \leq n \mu $. Then we get
$$ m_1 + r m_2 < m_2 + r m_2 = (r+1) m_2 \leq n m_2 \leq n \mu, $$ 
a contradiction.
\end{pf}

\vspace{0.5cm}

We shall now consider the case where $ m_2 < m_1 $. Then we do not necessarily have that $ \mu \leq m_2 $, but we shall see that this is ``stably'' true. Note that the case $ m_1 = m_2 $ is treated in Corollary \ref{m_1=m_2}.

\begin{lemma}
Let us keep the hypotheses from Lemma \ref{lemma 8.4}. Assume in addition that $ m_2 \leq \mu $. Then there exists a positive integer $ K $ such that $ m_1 + r' m_2 = n' \mu' $, where $ n' = n + K M $, $ 0 < r' < n' $ and where $ \mu' \leq m_2 $. The equality occurs only for $ \mu' = m_2 = \Gcd(m_1,m_2) $.
\end{lemma}
\begin{pf}
In case $ \mu = m_2 $, we take $K = 0 $ and hence $ n' = n $, so we may assume that $ \mu > m_2 $.

Observe that if $ n' = n + kM $, and $ R' = r + k \mu \frac{M}{m_2} $, we get that 
$$ m_1 + R' m_2 = m_1 + (r + k \mu \frac{M}{m_2}) m_2 = m_1 + r m_2 + k \mu M = n \mu + k \mu M = n' \mu. $$
But $ n' - R' = (n-r) + k (M - \mu \frac{M}{m_2}) $, so if $k \gg 0$, this is negative. Let $k_0$ be the smallest integer such that $ (n-r) + k_0 (M - \mu \frac{M}{m_2}) $ is negative, and put $ n' = n + k_0 M $. Arguing as in Lemma \ref{lemma 8.4}, we see that $m_1$, $m_2$ and hence $R'$ are relatively prime to $n'$. Let then $l_0$ be the unique integer such that $ 0 < R' - l_0 n' < n' $, and define $ r' = R' - l_0 n' $. It is easily computed that
$$ m_1 + r' m_2 = n' (\mu - l_0 m_2). $$
We now claim that $ 0 < \mu' := \mu - l_0 m_2 \leq m_2 $. For if this was not the case, then $ m_2 < \mu - l_0 m_2 $, and we could find an integer $ l > l_0 $ such that $ 0 < \mu - l m_2 \leq m_2 $, giving the equation
$$ m_1 + r'' m_2 = n' (\mu - l m_2), $$ 
where $ r'' = R' - l n' $, and where we necessarily have $ r'' < 0 $, by assumption on $l_0$. But then $ m_1 + r'' m_2 < m_1 < n' $, and as $ n' (\mu - l m_2) \geq n' $, we obtain a contradiction.

Finally, if $ \mu' = m_2 $, then $ \mu' $ divides $m_1$, and hence $ \mu' = m_2 = \Gcd(m_1,m_2) $. 
\end{pf}

\begin{lemma}\label{lemma 8.6}
Let $ \sigma = (m_1,m_2,n) $ be a singularity, and let $ \mu_l $, $ b_l  $ and $r_l$ be the numerical data associated to the minimal resolution of $ \sigma $. Assume that there exists an integer $ \lambda > 0 $ such that 
$$ \mu_0 > \mu_1 > \ldots > \mu_{\lambda}. $$
 Let $ n' = n + k M $ for any positive integer $k$, and denote by $ \mu_l', b_l' $ and $ r_l' $ the numerical data associated to the resolution of the singularity $ \sigma' = (m_1,m_2,n') $.

Then we have that $ \mu_l' = \mu_l $ for all $ 0 \leq l \leq \lambda $, that $ b_l' = b_l $ for all $ 1 \leq l \leq \lambda - 1 $ and that $ r_l' = r_l + k \frac{M}{m_2} \mu_{l+1} $ for all $ -1 \leq l \leq \lambda - 1 $.
\end{lemma}
\begin{pf}
We shall prove this by induction. Notice that if $ \lambda = 1 $, this is just Lemma \ref{lemma 8.4}. So we may assume that $ \lambda > 1 $. Then the statement for $ \mu_1' $ and $ r_0' $ again follows from Lemma \ref{lemma 8.4}.

We have now that $ \mu_2 < \mu_1 $. With $ n' = n + kM $ and $ r_0' = r_0 + k \frac{M}{m_2} \mu_1 $, we need to show that $ \mu_2' = \mu_2 $, $ b_1' = b_1 $ and $ r_1' = r_1 + k \frac{M}{m_2} \mu_2 $. Let us define $ R_1' = r_1 + k \frac{M}{m_2} \mu_2 $. We now compute
$$ b_1 r_0' - n' = b_1 (r_0 + k \frac{M}{m_2} \mu_1) - (n + k M) = (b_1 r_0 - n) + k \frac{M}{m_2}(b_1 \mu_1 - \mu_0) $$ 
$$ = r_1 + k \frac{M}{m_2} \mu_2 = R_1'. $$
But by assumption $ 0 \leq r_1 < r_0 $ and $ 1 \leq \mu_2 < \mu_1 $, so we get that
$$ 0 < R_1' = r_1 + k \frac{M}{m_2} \mu_2 < r_0 + k \frac{M}{m_2} \mu_1 = r_0'. $$
In other words, $b_1$ is the unique integer such that $ 0 < b_1 r_0' - n' < r_0' $, and therefore $ b_1' = b_1 $ and $ r_1' = R_1' = r_1 + k \frac{M}{m_2} \mu_2 $. Furthermore, it follows that 
$$ \mu_2' = b_1' \mu_1' - \mu_0' = b_1 \mu_1 - \mu_0 = \mu_2. $$

Assume now that we have established that $ \mu_i' = \mu_i $ for all $ 0 \leq i \leq l $, that $ b_i' = b_i $ for all $ 1 \leq i \leq l-1 $ and that $ r_i' = r_i + k \frac{M}{m_2} \mu_{i+1} $ for all $ -1 \leq i \leq l-1 $, where $ 2 \leq l < \lambda $. We need to prove that $ \mu_{l+1}' = \mu_{l+1} $, that $ b_l' = b_l$ and that $ r_l' = r_l + k \frac{M}{m_2} \mu_{l+1} $. Define $ R_l' = r_l + k \frac{M}{m_2} \mu_{l+1} $. We then compute that
$$ b_l r_{l-1}' - r_{l-2}' = b_l(r_{l-1} + k \frac{M}{m_2} \mu_l) - (r_{l-2} + k \frac{M}{m_2} \mu_{l-1}) = $$
$$ (b_l r_{l-1} - r_{l-2}) + k \frac{M}{m_2} (b_l \mu_l - \mu_{l-1}) = r_l + k \frac{M}{m_2} \mu_{l+1} = R_l'. $$
Since $ 0 \leq r_l < r_{l-1} $ and $ 1 \leq \mu_{l+1} < \mu_l $, we get that 
$$ 0 < R_l' = r_l + k \frac{M}{m_2} \mu_{l+1} < r_{l-1} + k \frac{M}{m_2} \mu_l = r_{l-1}', $$
and so it follows that $ b_l' = b_l $ and that $ r_l' = r_l + k \frac{M}{m_2} \mu_{l+1} $. Finally, we also see that
$$ \mu_{l+1}' = b_l' \mu_l' - \mu_{l-1}' = b_l \mu_l - \mu_{l-1} = \mu_{l+1}, $$
which completes the proof.
\end{pf}

\begin{lemma}\label{lemma 8.7}
Consider again a singularity $ \sigma = (m_1,m_2,n) $. Assume that there exists a positive integer $ \lambda $ such that the multiplicities of the exceptional components in the minimal resolution of $ \sigma $ satisfy the inequalities 
$$ \mu_0 > \mu_1 > \ldots > \mu_{\lambda} < \mu_{\lambda + 1}. $$

Then there exists a positive integer $k_0$ such that the multiplicities of the components in the minimal resolution of the singularity $ \sigma' = (m_1,m_2,n') $, where $ n' = n + k_0 M $, satisfy $ \mu_0' = \mu_0, \ldots, \mu_{\lambda}' = \mu_{\lambda} $, and $ \mu_{\lambda + 1}' < \mu_{\lambda + 1} $.
\end{lemma}
\begin{pf}
We will use the notation and computations from the proof of Lemma \ref{lemma 8.6}.

We have that
\begin{equation}\label{equation 8.3}
R_{\lambda}' := r_{\lambda} + k \frac{M}{m_2} \mu_{\lambda + 1} = b_{\lambda} r_{\lambda - 1}' - r_{\lambda -2 }', 
\end{equation}
and $ r_{\lambda - 1}' = r_{\lambda - 1} + k \frac{M}{m_2} \mu_{\lambda} $. Observe that
$$ r_{\lambda - 1}' - R_{\lambda}' = (r_{\lambda - 1} + k \frac{M}{m_2} \mu_{\lambda}) - (r_{\lambda} + k \frac{M}{m_2} \mu_{\lambda + 1}) = (r_{\lambda - 1} - r_{\lambda}) + k \frac{M}{m_2} (\mu_{\lambda} - \mu_{\lambda + 1}). $$
By assumption, $ r_{\lambda - 1} - r_{\lambda} > 0 $ and $ \mu_{\lambda} - \mu_{\lambda + 1} < 0 $, so for $ k \gg 0 $, we get that $ r_{\lambda - 1}' < R_{\lambda}' $.

Let $ k_0 $ be the smallest integer such that $ r_{\lambda - 1}' < R_{\lambda}' $. Since $ r_{\lambda}' < r_{\lambda - 1}' $, it then follows from Equation \ref{equation 8.3} that $ b_{\lambda}' < b_{\lambda} $. Furthermore, we easily compute from Equation \ref{equation 8.3} that $ (b_{\lambda} - l) r_{{\lambda} - 1}' -  r_{\lambda - 2}' = R_{\lambda}' - l r_{\lambda - 1}' $. Now, let $ l_0 $ be the biggest integer such that $ 0 \leq R_{\lambda}' - l_0 r_{\lambda - 1}' $. We get that $ 0 \leq (b_{\lambda} - l_0) r_{\lambda - 1}' -  r_{\lambda - 2}' < r_{\lambda - 1}' $, so it follows that $ b_{\lambda}' = b_{\lambda} - l_0 $ and $ r_{\lambda}' = R_{\lambda}' - l_0 r_{\lambda - 1}' $. But then 
$$ \mu_{\lambda + 1}' = b_{\lambda}' \mu_{\lambda}' - \mu_{\lambda - 1}' = (b_{\lambda} - l_0) \mu_{\lambda} - \mu_{\lambda - 1} = \mu_{\lambda + 1} - l_0 \mu_{\lambda} < \mu_{\lambda + 1}. $$
The rest of the statement follows immediately from Lemma \ref{lemma 8.6}.
\end{pf}

\begin{cor}\label{cor. 8.8}
Let $ \sigma = (m_1,m_2,n) $ be a singularity, and assume that there exists a positive integer $\lambda$ such that the multiplicities of the components in the minimal resolution of $\sigma$ satisfy the inequalities 
$$ \mu_0 > \mu_1 > \ldots > \mu_{\lambda} < \mu_{\lambda + 1}. $$

Then there exists a positive integer $K$ such that the multiplicities of the components in the minimal resolution of the singularity $ \sigma' = (m_1,m_2,n') $, where $ n' = n + K M $, satisfy $ \mu_0' = \mu_0, \ldots, \mu_{\lambda}' = \mu_{\lambda} $, and $ \mu_{\lambda + 1}' \leq \mu_{\lambda}' $. The equality occurs only in the case where $ \mu_{\lambda + 1}' = \mu_{\lambda}' = \Gcd(m_1,m_2) $.
\end{cor}
\begin{pf}
According to Lemma \ref{lemma 8.7}, we can find an integer $k_0$ such that for the extension of degree $n + k_0 M$, the associated multiplicities satisfy the equalities $ \mu_0' = \mu_0, \ldots, \mu_{\lambda}' = \mu_{\lambda} $, and the inequality $ \mu_{\lambda + 1}' < \mu_{\lambda + 1} $. If we still have that $ \mu_{\lambda}' < \mu_{\lambda + 1}' $, we can apply the same lemma again, by replacing $n$ with $n + k_0 M$, and looking at the multiplicities of the components in the exceptional fibre of the minimal resolution of this singularity instead.  By repeated use of this procedure, the multiplicity of the $ (\lambda + 1) $-st component in the exceptional fiber will strictly decrease each time, so eventually we will have $ \mu_{\lambda + 1}' \leq \mu_{\lambda}' $. By then, we have made an extension of degree $n + KM$, which is the sought after extension.

If the inequality above is actually an equality, i.e., $ \mu_{\lambda + 1}' = \mu_{\lambda}' $, it follows from Lemma \ref{lemma 8.1} that $ \mu_{\lambda}' $ divides \emph{all} the multiplicities, and in particular it divides $m = \Gcd(m_1,m_2) $. On the other hand, $ m $ divides all multiplicities in the chain, so we get that $ \mu_{\lambda + 1}' = \mu_{\lambda}' = \Gcd(m_1,m_2) $.
\end{pf}

\vspace{0.5cm}

Combining the results developed in this section, we get the following result:

\begin{prop}\label{prop. 8.9}
Let $m_1, m_2$ be positive integers, let $ m = \Gcd(m_1,m_2) $, and let $ M = \Lcm(m_1,m_2) $. Let us furthermore fix a positive integer $ n_0 $ that is not divisible by $p$ and that is relatively prime to $M$. Then the following properties hold:
\begin{enumerate}
\item There exists an integer $ K \gg 0 $, such that the multiplicities of the components in the minimal resolution of the singularity $ \sigma = (m_1,m_2,n) $, where $ n = n_0 + K M $, satisfy
$$ \mu_0 > \mu_1 > \ldots > \mu_{l_0} =  \ldots = m = \ldots = \mu_{L + 1 - l_1} < \ldots < \mu_{L} < \mu_{L + 1}, $$
where $L$ denotes the \emph{length} of the singularity $ \sigma $.
\item The integers $ \mu_2 , \ldots , \mu_{l_0} $ are uniquely determined by $ \mu_0 $ and $ \mu_1 $, and similarly $ \mu_{L + 1 - l_1}, \ldots , \mu_{L-1} $ are uniquely determined by $ \mu_{L} $ and $ \mu_{L+1} $.
\item For any extension of degree $ n' = n + k M $, where $ k > 0 $, we have that the multiplicities $ \mu_l' $ of the components in the minimal resolution of the singularity $ \sigma' = (m_1,m_2,n') $ will only differ from the sequence of multiplicities associated to $ \sigma $ by inserting $m$'s ``in the middle''. 
\end{enumerate}
\end{prop}
\begin{pf}
(i) From Corollary \ref{cor. 8.8}, it follows that we can find a positive integer $K \gg 0$ such that the multiplicities associated to the singularity $ \sigma = (m_1,m_2,n) $, where $ n = n_0 + KM $, satisfy the inequalities
$$ \mu_0 > \mu_1 > \ldots > \mu_{l_0} = m $$
and
$$ \mu_{L+1} > \mu_{L} > \ldots > \mu_{L + 1 - l_1} = m, $$
for some integers $ l_0 \geq 0 $ and $ l_1 \geq 0 $.

For any $ \mu_l $, where $ l_0 < l < L + 1 - l_1 $, we have that $ \mu_l = m $. Indeed, by part (ii) of Lemma \ref{lemma 8.1}, we have that $ m \geq \mu_l $, and part (iii) of the same lemma implies that we actually have equalities. 

(ii) In case $ l_0 = 0 $ or $ l_0 = 1 $, the statement is empty. So we can assume that $ l_0 > 1 $. Let us first see that $ \mu_0 $ and $ \mu_1 $ determine $ \mu_2 $. By assumption, we have that $ \mu_0 > \mu_1 > \mu_2 \geq m $. These integers must satisfy the equation $ \mu_2 = b_1 \mu_1 - \mu_0 $, and since there is only one such integer $ \mu_2 $ satisfying these properties, it follows that $ \mu_2 $ is uniquely determined by $ \mu_0 $ and $ \mu_1 $.

Now let $ 1 < l < l_0 $, and assume that $ \mu_0 > \mu_1 > \ldots > \mu_{l} $ are given. Then $ \mu_{l+1} = b_l \mu_l - \mu_{l-1} $, and we have that $ \mu_{l+1} < \mu_l $. But then  $ \mu_{l+1} $ is uniquely determined by $ \mu_l $ and $ \mu_{l-1} $, which are in turn uniquely determined by $ \mu_0 $ and $ \mu_1 $ by the induction hypothesis, so the result follows. The same argument applies to the other end of the chain. 

(iii) This is a straightforward consequence of Lemma \ref{lemma 8.6}.
\end{pf}

\begin{rmk}
Proposition \ref{prop. 8.9} states that the minimal resolution of a singularity $(m_1,m_2,n)$ essentially only depends on $m_1$, $m_2$ and the residue class of $ n $ modulo $ M = \Lcm(m_1,m_2) $. That is, if $ n' \gg 0 $, and $ n'' \gg 0 $, and if $ n' \equiv_M n'' $, then the exceptional locus of the minimal desingularization of $(m_1,m_2,n')$ is equal to the exceptional locus of the minimal desingularization of $(m_1,m_2,n'')$, modulo chains of $(-2)$-curves ``in the middle''. 
\end{rmk}

\section{Special filtrations for trace computations}\label{special filtrations}
Let $ \mathcal{X}/S $ be an SNC-model, and let $ S' \rightarrow S $ be a tame extension of degree $n$, where $n$ is prime to the least common multiple of the multiplicities of the irreducible components of $ \mathcal{X}_k $. Let $ \mathcal{X}' $ be the normalization of $ \mathcal{X}_{S'} $, and let $ \mathcal{Y} $ be the minimal desingularization of $ \mathcal{X}' $. 

This section is devoted to computing, for any $ g \in G = \boldsymbol{\mu}_n $, the Brauer trace of the automorphism induced by $g$ on the formal difference $ H^0(\mathcal{Y}_k, \mathcal{O}_{\mathcal{Y}_k}) - H^1(\mathcal{Y}_k, \mathcal{O}_{\mathcal{Y}_k}) $. Hence a lot of our previous work will come together in this section. 

Our assumption on the degree of $ S'/S $ makes it possible to describe $ \mathcal{Y}_k $ in terms of $ \mathcal{X}_k $. In particular, since every component of $ \mathcal{Y}_k $ either is an exceptional curve, or dominates a component of $ \mathcal{X}_k $, it is natural to stratify the combinatorial structure of $ \mathcal{Y}_k $ according to the combinatorial structure of $ \mathcal{X}_k $.

This stratification proves to be very convenient for our trace computations. The section concludes with Theorem \ref{thm. 9.13}, wich gives a formula for the trace mentioned above as a sum of contributions associated in a natural way to the combinatorial structure of $ \mathcal{X}_k $. 

\subsection{The graph $ \Gamma(\mathcal{X}_k) $}\label{9.1}
We will associate a graph $ \Gamma(\mathcal{X}_k) $ to $ \mathcal{X}_k $ in the following way: The set of vertices, $ \mathcal{V} $, consists of the irreducible components of $ \mathcal{X}_k $. The set of edges, $ \mathcal{E} $, consists of the intersection points of $ \mathcal{X}_k $, and two distinct vertices $\upsilon$ and $\upsilon'$ are connected by $ \Card (\{ D_{\upsilon} \cap D_{\upsilon'} \}) $ edges, where $D_{\upsilon}$ is the irreducible component corresponding to $\upsilon$.

We define two natural functions on the set of vertices $ \mathcal{V} $. First, define the \emph{genus}
$$ \mathfrak{g} : \mathcal{V} \rightarrow \mathbb{N}_0, $$ 
by $ \mathfrak{g}(\upsilon) = p_a(D_{\upsilon}) $. We also define the \emph{multiplicity} 
$$ \mathfrak{m} : \mathcal{V} \rightarrow \mathbb{N}, $$ 
by $ \mathfrak{m}(\upsilon) = \Mult_{\mathcal{X}_k}(D_{\upsilon}) $.

The graph $ \Gamma(\mathcal{X}_k) $, together with the functions $\mathfrak{g}$ and $\mathfrak{m}$, encode all the combinatorial and numerical properties of $ \mathcal{X}_k $.

\subsection{A partition of the set of irreducible components of $ \mathcal{Y}_k $}\label{9.2}

Let $ \mathcal{S} $ denote the set of irreducible components of $ \mathcal{Y}_k $. If $ C \in \mathcal{S} $, then we have either: 
\begin{enumerate}
\item $ C $ dominates a component $ D_{\upsilon} $ of $ \mathcal{X}_k $, or  

\item $ C $ is a component of the exceptional locus of the minimal desingularization $ \rho : \mathcal{Y} \rightarrow \mathcal{X}' $.
\end{enumerate}

In the first case, we have that $ p_a(C) = \mathfrak{g}(\upsilon) $, and $ \Mult_{\mathcal{Y}_k}(C) = \mathfrak{m}(\upsilon) $. Furthermore, $G$ acts trivially on $C$. Since $C$ is the unique component of $ \mathcal{Y}_k $ corresponding to $ \upsilon $, we write $C = C_{\upsilon}$.

In the second case, we have that $ C $ is part of a chain of exceptional curves, corresponding uniquely to an edge $ \varepsilon \in \mathcal{E} $. Hence $ p_a(C) = 0 $. By choosing an ordering (or direction) of this chain, we can index the components in the chain by $ l $, for $ 1 \leq l \leq L(\varepsilon) $, where $ L(\varepsilon) $ is the length of the chain. So we can write $ C = C_{\varepsilon,l} $, for some $l \in \{1, \ldots, L(\varepsilon) \} $. By Lemma \ref{Lemma 8.4}, $G$ acts nontrivially on $ C $, with fixed points exactly at the two points where $C$ meets the rest of the special fiber.

The special fiber $ \mathcal{Y}_k $ can now be written, as an effective divisor on $ \mathcal{Y} $, in the form 
$$ \mathcal{Y}_k = \sum_{\varepsilon \in \mathcal{E}} \sum_{l = 1}^{L(\varepsilon)} \mu_{\varepsilon,l} C_{\varepsilon,l} + \sum_{\upsilon \in \mathcal{V}} m_{\upsilon} C_{\upsilon}, $$
where $ \mu_{\varepsilon,l} $ denotes the multiplicity of the component $ C_{\varepsilon,l} $, and $m_{\upsilon}$ is the multiplicity of $C_{\upsilon}$.

\subsection{Special filtrations}\label{9.3}
We will now consider \emph{special} filtrations of $ \mathcal{Y}_k $, inspired by the partition of the set of irreducible components of $ \mathcal{Y}_k $ introduced above. 

Let us to begin with choosing an ordering of the elements in $ \mathcal{V} $. We can then define the following sequence:
$$ 0 < \ldots < Z_{\mathcal{E}} =: Z_{\upsilon_{|\mathcal{V}| + 1}} < Z_{\upsilon_{|\mathcal{V}|}} < \ldots < Z_{\upsilon_i} < \ldots < Z_{\upsilon_1} = \mathcal{Y}_k, $$
where $ Z_{\mathcal{E}} := \mathcal{Y}_k - \sum_{\upsilon \in \mathcal{V}} m_{\upsilon} C_{\upsilon} $. The $ Z_{\upsilon_i} $ are defined inductively, for every $ i \in \{ 1, \ldots, |\mathcal{V}| \} $, by the refinements
$$ Z_{\upsilon_{i+1}} = Z_{\upsilon_i}^{m_{\upsilon_i} + 1} < \ldots < Z_{\upsilon_i}^j < \ldots < Z_{\upsilon_i}^1 = Z_{\upsilon_i}, $$ 
where $ Z_{\upsilon_i}^{j + 1} = Z_{\upsilon_i} - j C_{\upsilon_i} $ for every $ j \in \{ 0, \ldots, m_{\upsilon_i} \} $.

Next, we choose an ordering of the elements in $ \mathcal{E} $. We can then define the following sequence:
$$ 0 =: Z_{\varepsilon_{|\mathcal{E}|+1}} < Z_{\varepsilon_{|\mathcal{E}|}} < \ldots < Z_{\varepsilon_i} < \ldots < Z_{\varepsilon_1} := Z_{\mathcal{E}}. $$  
 The $ Z_{\varepsilon_i} $ are defined inductively, for any $ i \in \{ 1, \ldots, |\mathcal{E}| \} $, by the refinements
$$ Z_{\varepsilon_{i+1}} := Z_{\varepsilon_i, L(\varepsilon_i) + 1} < \ldots < Z_{\varepsilon_i,l} < \ldots < Z_{\varepsilon_i,1} := Z_{\varepsilon_i}, $$
which in turn are defined inductively, for every $ l \in \{ 1, \ldots, L(\varepsilon_i) \} $, by the further refinements
$$ Z_{\varepsilon_i, l+1} := Z_{\varepsilon_i, l}^{\mu_l + 1} < \ldots < Z_{\varepsilon_i, l}^j < \ldots < Z_{\varepsilon_i, l}^1 := Z_{\varepsilon_i, l}, $$
where $ Z_{\varepsilon_i, l}^{j+1} := Z_{\varepsilon_i, l} - j C_{\varepsilon_i, l} $, for every $ j \in \{ 0, \ldots, \mu_l \} $.

\subsection{} 
In the rest of this paper, we shall always choose complete filtrations of $ \mathcal{Y}_k $ that are of the form
\begin{equation}
0 < \ldots < Z_{\varepsilon_i} < \ldots < Z_{\varepsilon_1} = Z_{\mathcal{E}} < \ldots < Z_{\upsilon_i} < \ldots < Z_{\upsilon_1} = \mathcal{Y}_k, 
\end{equation}
where $ Z_{\upsilon_{i+1}} < Z_{\upsilon_i} $ and $ Z_{\varepsilon_{i+1}} < Z_{\varepsilon_i} $ are subfiltrations as described above. 

We shall soon see that the chosen orderings of the sets $ \mathcal{E} $ and $ \mathcal{V} $ are irrelevant. The nice feature of working with filtrations like this becomes evident when one wants to do trace computations \`a la Section \ref{section 6}. Then we may actually reduce to considering subfiltrations $ Z_{\upsilon_{i+1}} < Z_{\upsilon_i} $, which we interpret as contributions from the vertices of $\Gamma$, and subfiltrations $ Z_{\varepsilon_{i+1}} < Z_{\varepsilon_i} $, which we interpret as contributions from the edges.

\subsection{Contribution to the trace from vertices in $ \mathcal{V}$}
Let us fix a vertex $\upsilon \in \mathcal{V}$. We shall now define and calculate the \emph{contribution} to the trace from $\upsilon$. To do this, we choose a filtration of $ \mathcal{Y}_k $ as in Section \ref{9.3}. Then there will be a subfiltration of the form:
 $$ Z_{\mathcal{E}} \leq Z_{\upsilon}^{m_{\upsilon}+1} < \ldots < Z_{\upsilon}^k < \ldots < Z_{\upsilon}^1 = Z_{\upsilon} \leq \mathcal{Y}_k, $$
where $ Z_{\upsilon}^{k} - Z_{\upsilon}^{k + 1} = C_{\upsilon} $, for all $ 1 \leq k \leq m_{\upsilon} $. The invertible sheaf associated to the $k$-th step in this filtration is $ \mathcal{L}_{\upsilon}^k := j_{\upsilon}^*(\mathcal{I}_{Z_{\upsilon}^{k+1}}) $, where $ j_{\upsilon} : C_{\upsilon} \hookrightarrow \mathcal{Y} $ is the canonical inclusion.

We will use the following easy lemma, whose proof is omitted.
\begin{lemma}\label{lemma 9.1}
Assume that $ S'/S $ is a nontrivial extension. If $C_1$ and $C_2$ are two distinct components of $ \mathcal{Y}_k $, corresponding to elements in $ \mathcal{V} $, then they have empty intersection.
\end{lemma}

Let $ D_1, \ldots, D_f $ be the components of $ Z_{\upsilon} $ that intersect $ C_{\upsilon} $ non-trivially, and that are not equal to $ C_{\upsilon} $. Let $ a_i $ denote the multiplicity of $D_i$. It follows from Lemma \ref{lemma 9.1} that the $D_i$ are exceptional components. Moreover, it follows from the way we constructed the filtration that the $D_i$ are precisely the components of $ \mathcal{Y}_k $ different from $ C_{\upsilon} $ that have non-empty intersection with $ C_{\upsilon} $. We can then write 
$$ Z_{\upsilon}^{k + 1} = (m_{\upsilon} - k) C_{\upsilon} + a_1 D_1 + \ldots + a_f D_f + Z_0, $$
where all components of $ Z_0 $ have empty intersection with $ C_{\upsilon} $. So we get that
\begin{equation}\label{equation 9.2}
\mathcal{L}_{\upsilon}^k  =  j_{\upsilon}^*(\mathcal{I}_{Z_{\upsilon}^{k+1}}) = (\mathcal{I}_{C_{\upsilon}}|_{C_{\upsilon}})^{ \otimes m_{\upsilon} - k } \otimes (\mathcal{I}_{D_1}|_{C_{\upsilon}})^{\otimes a_1} \otimes \ldots \otimes (\mathcal{I}_{D_f}|_{C_{\upsilon}})^{\otimes a_f}. 
\end{equation}            

Let $ g $ be an element of $ G = \boldsymbol{\mu}_n $, corresponding to a root of unity $ \xi $. Note that the restriction of the automorphism $ g $ to $C_{\upsilon}$ is $ \Id_{C_{\upsilon}}$. Let 
$$ u_{\upsilon}^k : g^* (i_{C_{\upsilon}})_* \mathcal{L}_{\upsilon}^k \rightarrow (i_{C_{\upsilon}})_* \mathcal{L}_{\upsilon}^k $$
be the map constructed in Section \ref{6.8}, and let $ j_{\upsilon}^* (u_{\upsilon}^k) : \mathcal{L}_{\upsilon}^k \rightarrow \mathcal{L}_{\upsilon}^k $ be the pullback of this map to $C_{\upsilon}$. 

\begin{dfn}\label{contribution from vertex}
We define \emph{the contribution to the trace} from the vertex $ \upsilon \in \mathcal{V} $ as the sum
\begin{equation}
\Tr_{\upsilon}(\xi) = \sum_{k=1}^{m_{\upsilon}} \Tr_{\beta}(e(H^{\bullet}(\Id_{C_\upsilon},j_{\upsilon}^*(u_{\upsilon}^k)))) .
\end{equation}
\end{dfn}

\begin{prop}\label{prop. 9.7} The contribution to the trace from the vertex $ \upsilon $ is 
$$ \Tr_{\upsilon}(\xi) = \sum_{k=1}^{m_{\upsilon}} (\xi^{\alpha_{m_{\upsilon}}})^{m_{\upsilon} - k} (k C_{\upsilon}^2 + 1 - p_a(C_{\upsilon})) $$
$$ = \sum_{k=0}^{m_{\upsilon}-1} (\xi^{\alpha_{m_{\upsilon}}})^{k} ((m_{\upsilon} - k) C_{\upsilon}^2 + 1 - p_a(C_{\upsilon})), $$
where $ \alpha_{m_{\upsilon}} $ is an inverse to $ m_{\upsilon} $ modulo $n$. 
\end{prop}
\begin{pf}
By Proposition \ref{prop. 5.6}, we have that
$$ \Tr_{\beta}(e(H^{\bullet}(\Id_{\upsilon}, j_{\upsilon}^*(u_{\upsilon}^k)))) = \lambda_k(\Deg_{C_{\upsilon}}(\mathcal{L}_{\upsilon}^k) + 1 - p_a(C_{\upsilon})), $$
where $ \lambda_k = \Tr(j_{\upsilon}^*(u_{\upsilon}^k)(y)) $, for any point $ y \in C_{\upsilon} $.
The proof will consist of specifying precisely the terms appearing in this formula. 

Let us first compute $ \Deg_{C_{\upsilon}}(\mathcal{L}_{\upsilon}^k) $. Since $ \mathcal{I}_{C_{\upsilon}} = \mathcal{O}_{\mathcal{Y}}(- C_{\upsilon}) $, it follows that
$$ \Deg_{C_{\upsilon}}( \mathcal{I}_{C_{\upsilon}}|_{C_{\upsilon}}) = \Deg_{C_{\upsilon}}(\mathcal{O}_{\mathcal{Y}}(- C_{\upsilon})|_{C_{\upsilon}}) = - \Deg_{C_{\upsilon}}(\mathcal{O}_{\mathcal{Y}}(C_{\upsilon})|_{C_{\upsilon}}) = - C_{\upsilon}^2.$$
Furthermore, for any $ i \in \{ 1, \ldots, f \} $, we have that $ \mathcal{I}_{D_i} = \mathcal{O}_{\mathcal{Y}}(- D_i) $, and hence
$$ \Deg_{C_{\upsilon}}(\mathcal{I}_{D_i}|_{C_{\upsilon}}) = \Deg_{C_{\upsilon}}(\mathcal{O}_{\mathcal{Y}}(- D_i)|_{C_{\upsilon}}) = - 1, $$
since $ C_{\upsilon} $ and $ D_i $ intersect transversally at exactly one point. Since $ \Deg_{C_{\upsilon}}(-) $ is additive on tensor products, it follows from Equation \ref{equation 9.2} that
$$ \Deg_{C_{\upsilon}}(\mathcal{L}_{\upsilon}^k) = - (m_{\upsilon} - k) C_{\upsilon}^2 - (a_1 + \ldots + a_f). $$ 
On the other hand, we have that $ - C_{\upsilon}^2 = (a_1 + \ldots + a_f)/m_{\upsilon} $, and therefore $ \Deg_{C_{\upsilon}}(\mathcal{L}_{\upsilon}^k) = k C_{\upsilon}^2 $.

To see that $ \lambda_k = (\xi^{\alpha_{m_{\upsilon}}})^{m_{\upsilon} - k} $, let $D$ be one of the components of $ \mathcal{Y}_k $ meeting $ C_{\upsilon} $, and denote by $y$ the unique point where they intersect. Then $D$ is part of a chain of exceptional curves. Let $ C_{\upsilon'} $ be the component at ``the other end'' of this chain, and denote by $L$ the length of this chain. Then, using the notation and computations in Proposition \ref{prop. 3.3}, with $ C_{\upsilon} = C_{L+1} $ and $ D = C_L $, we can identify the fiber of  $ \mathcal{I}_{C_{\upsilon}} $ at $ y = y_L $ with $ <z_{L+1}> $. And the eigenvalue of $ z_{L+1} $ for the automorphism induced by $\xi$ was precisely equal to $ \xi^{\alpha_{m_{\upsilon}}} $. It follows that $ \lambda_k = (\xi^{\alpha_{m_{\upsilon}}})^{m_{\upsilon} - k} $. 
\end{pf}

\begin{rmk}
In particular, it is clear that this formula is independent of how we have chosen to order the elements in $ \mathcal{V} $.
\end{rmk}


\subsection{Contribution to the trace from edges in $ \mathcal{E} $}\label{9.6}
Let us now choose an edge $ \varepsilon \in \mathcal{E} $. In the filtration of $ \mathcal{Y}_k $, we can find a subfiltration $ 0 < Z_{\varepsilon} \leq Z_{\mathcal{E}} < \mathcal{Y}_k $, with the refinements
$$ Z_{\varepsilon, L(\varepsilon) + 1} < \ldots < Z_{\varepsilon,l} < \ldots < Z_{\varepsilon,1} = Z_{\varepsilon}, $$
for any $ l \in \{ 1, \ldots, L(\varepsilon) \} $, and further refinements
$$ Z_{\varepsilon,l+1} = Z_{\varepsilon,l}^{\mu_l + 1} < \ldots < Z_{\varepsilon,l}^{k} < \ldots < Z_{\varepsilon,l}^{1} = Z_{\varepsilon,l}, $$ 
where $ Z_{\varepsilon,l}^{k} - Z_{\varepsilon,l}^{k+1} = C_{\varepsilon,l} $, for any $ k \in \{ 1, \ldots, \mu_l \} $.


 As we are working with a fixed $ \varepsilon $, we will for the rest of this section suppress the index $ \varepsilon $, to simplify the notation. Take now an integer $ l \in \{1, \ldots, L-1 \} $, and let $ j_l : C_l \hookrightarrow \mathcal{Y} $ be the canonical inclusion. Consider then the subfiltration involving the component $C_l$:
$$ \ldots < Z^{\mu_l + 1}_l < \ldots < Z^k_l < \ldots < Z^1_l < \ldots . $$
At the $k$-th step in this filtration, we have $ Z_l^k - Z_l^{k+1} = C_l $ for all $ 1 \leq k \leq \mu_l $. The associated invertible sheaf at the $k$-th step is
$$ \mathcal{L}_l^k := j_l^*( \mathcal{I}_{Z_l^{k+1}}) = (\mathcal{I}_{C_l}|_{C_l})^{\otimes (\mu_l - k)} \otimes (\mathcal{I}_{C_{l+1}}|_{C_l})^{ \otimes \mu_{l+1} }. $$

For $ l = L $, we have 
$$ \ldots < Z^{\mu_L + 1}_L < \ldots < Z^k_L < \ldots < Z^1_L < \ldots, $$
and at the $ k $-th step, we have $ Z_L^k - Z_L^{k+1} = C_L $. As all components in $ Z_{L} $ other than $ C_L $ have empty intersection with $ C_L $, we get that
$$ \mathcal{L}_L^k := j_L^*( \mathcal{I}_{Z_L^{k+1}}) = (\mathcal{I}_{C_L}|_{C_L})^{\otimes (\mu_L - k)}. $$

Let $ g \in G $ be a group element corresponding to a root of unity $ \xi $. The restriction $ g|_{C_l} $ is either the identity on $C_l$, or has fixed points exactly at the two points $ y_l $ and $ y_{l-1} $ where $C_l$ meets the rest of the special fiber. We need to compute the fibers at $ y_l $ and $ y_{l-1} $ of $ \mathcal{L}_l^k $, and the corresponding eigenvalues for the automorphisms induced by $ g $ at these fibers.

Let $ u_l^k : g^* (j_l)_* \mathcal{L}_l^k \rightarrow (j_l)_* \mathcal{L}_l^k $ be the map defined in Section \ref{6.8}. We will now compute $ \Tr_{\beta}(e(H^{\bullet}(g|_{C_l},(j_l)^*(u_k^l)))) $.

\begin{lemma}\label{fiberlemma}
For any $ l \in \{1, \ldots, L \} $, we have that 
$$ \Tr((j_l)^*(u_l^k)(y_{l-1})) = (\xi^{\alpha_1 r_{l-2}})^{\mu_l - k}. $$
For any $ l \in \{1, \ldots, L-1 \} $, we have that 
$$ \Tr((j_l)^*(u_l^k)(y_{l})) = (\xi^{ - \alpha_1 r_l})^{\mu_l - k} (\xi^{ \alpha_1 r_{l-1}})^{\mu_{l+1}}, $$
and for $ l = L $, we have that 
$$ \Tr((j_L)^*(u_L^k)(y_{L})) = (\xi^{- \alpha_1 r_L})^{\mu_L - k}. $$
\end{lemma}
\begin{pf}
In order to compute the eigenvalue of the map $ (j_l)^*(u_k^l)(y) $ on the fiber at $y \in C_l$, it suffices, by Lemma \ref{lemma 7.3}, to compute the eigenvalue of $ v_l^k(y) $, where 
$$ v_l^k : g^* \mathcal{I}_{Z_l^{k+1}} \rightarrow \mathcal{I}_{Z_l^{k+1}} $$ 
is the map constructed as in Section \ref{6.4}. 

If $ l \in \{1, \ldots, L-1\} $, we can write
$$ \mathcal{I}_{Z_l^{k+1}} = \mathcal{I}_{C_l}^{\otimes(\mu_l - k)} \otimes \mathcal{I}_{C_{l+1}}^{\otimes \mu_{l+1}} \otimes \mathcal{I}_0, $$
where $ \mathcal{I}_0 $ has support away from $ C_l $. If $l = L $, we have 
$$ \mathcal{I}_{Z_L^{k+1}} = \mathcal{I}_{C_L}^{\otimes(\mu_L - k)} \otimes \mathcal{I}_0. $$
Recall the notation from Proposition \ref{prop. 3.3}. For any $ l \in \{1, \ldots, L \} $, the fiber at $ y_{l-1} $ is
$$ \mathcal{I}_{Z_l^{k+1}}(y_{l-1}) = \mathcal{I}_{C_l}^{\otimes(\mu_l - k)}(y_{l-1}) = <z_{l-1}>^{\otimes (\mu_l - k)}, $$
and the fiber at $y_l$ is, for any $ l \in \{1, \ldots, L-1 \} $, 
$$ \mathcal{I}_{Z_l^{k+1}}(y_l) = \mathcal{I}_{C_l}^{\otimes(\mu_l - k)}(y_l) \otimes \mathcal{I}_{C_{l+1}}^{\otimes \mu_{l+1}}(y_l) = <w_l>^{\otimes (\mu_l - k)} \otimes <z_l>^{ \otimes \mu_{l+1} }. $$  
In the case $ l = L $, we get
$$ \mathcal{I}_{Z_L^{k+1}}(y_L) = \mathcal{I}_{C_L}^{\otimes(\mu_L - k)}(y_L) = <w_L>^{\otimes (\mu_L - k)}. $$

By Lemma \ref{lemma 7.2}, $ \Tr(v_l^k(y_{l-1})) = \lambda^{\mu_l - k} $, where $ \lambda $ is the eigenvalue of the eigenvector $z_{l-1}$ under the automorphism induced by $ g $ on the cotangent space of $ \mathcal{Y} $ at $y_{l-1}$. Proposition \ref{prop. 3.3} then shows that $ \lambda = \xi^{\alpha_1 r_{l-2}} $, where $ \alpha_1 m_1 \equiv_n 1 $, so therefore 
\begin{equation}\label{equation 9.4} 
\Tr((j_l)^*(u_l^k)(y_{l-1})) = \Tr(v_l^k(y_{l-1})) = (\xi^{\alpha_1 r_{l-2}})^{\mu_l - k}.  
\end{equation}

For $ l \in \{1, \ldots, L-1 \} $, Lemma \ref{lemma 7.2} shows again that $ \Tr(v_l^k(y_{l})) = \lambda^{\mu_l - k} \lambda'^{\mu_{l+1}} $, where $ \lambda $ and $ \lambda' $ are the eigenvalues of the eigenvectors $ w_l $ and $ z_l $ under the automorphism induced by $ g $ on the cotangent space of $ \mathcal{Y} $ at $y_{l}$. Using Proposition \ref{prop. 3.3}, it is then easy to see that 
\begin{equation}\label{equation 9.5} 
\Tr((j_l)^*(u_l^k)(y_{l})) = \Tr(v_l^k(y_{l})) = (\xi^{ - \alpha_1 r_l})^{\mu_l - k} (\xi^{ \alpha_1 r_{l-1}})^{\mu_{l+1}}. 
\end{equation}
In the case $ l = L $, we get by the same methods
$$ \Tr((j_L)^*(u_L^k)(y_{L})) = \Tr(v_L^k(y_{L})) = (\xi^{- \alpha_1 r_L})^{\mu_L - k}. $$
\end{pf}

\begin{ntn}
We will write $ \chi = \xi^{\alpha_1} $.
\end{ntn}

\begin{prop}\label{lemma 9.8}
Let $ g \in G $ be a group element  corresponding to a root of unity $ \xi $. If $ g|_{C_l} = \Id_{C_l} $, then we have, for any $ l \in \{1, \ldots, L - 1 \} $, that
$$ \Tr_{\beta}(e(H^{\bullet}(g|_{C_l},(j_l)^*(u_l^k)))) = \chi^{ r_{l - 2} (\mu_l - k)} ((\mu_{l} - k) b_l - \mu_{l+1}+ 1). $$
If $l=L$, we have that
$$ \Tr_{\beta}(e(H^{\bullet}(g|_{C_L},(j_L)^*(u_L^k)))) = \chi^{ r_{L - 2} (\mu_L - k)} ((\mu_{L} - k) b_L + 1). $$
\end{prop}
\begin{pf}
Since $ g|_{C_l} = \Id_{C_l} $, we can apply Proposition \ref{prop. 5.5} (iii), which states that
$$ \Tr_{\beta}(e(H^{\bullet}(g|_{C_l},(j_l)^*(u_k^l)))) = \lambda(\Deg_{C_l}(\mathcal{L}_l^k) + 1 - p_a(C_l)). $$
From Lemma \ref{fiberlemma}, it follows that $ \lambda = \chi^{r_{l - 2} (\mu_l - k)}$. Furthermore, we have that $ p_a(C_l) = 0 $. So it remains to compute the degree of $ \mathcal{L}_l^k $. For $ l \neq L $, we have that $ \mathcal{L}_l^k = (\mathcal{I}_{C_l}|_{C_l})^{\otimes (\mu_l - k)} \otimes (\mathcal{I}_{C_{l+1}}|_{C_l})^{ \otimes (\mu_{l+1})} $, and we need therefore to compute $ \Deg_{C_l}(\mathcal{I}_{C_l}|_{C_l}) $ and $ \Deg_{C_l}(\mathcal{I}_{C_{l+1}}|_{C_l}) $. But it is easily seen that $ \Deg_{C_l}(\mathcal{I}_{C_l}|_{C_l}) = - C_l^2 = b_l $, and that $ \Deg_{C_l}(\mathcal{I}_{C_{l+1}}|_{C_l}) = - 1 $. Hence we get that $ \Deg_{C_l}(\mathcal{L}_l^k) = (\mu_l - k) b_l - \mu_{l+1} $. If $l=L$, we get instead that $ \Deg_{C_L}(\mathcal{L}_L^k) = (\mu_L - k) b_L $.
\end{pf}


\begin{prop}\label{lemma 9.10}
Assume that $ g|_{C_l} $ is not the identity on $C_l$. 
\begin{enumerate}
\item If $ l \in \{2, \ldots, L - 1 \} $, we have that 
$$ \Tr_{\beta}(e(H^{\bullet}(g|_{C_l},(j_l)^*(u_k^l)))) = \frac{\chi^{r_{l-2}(\mu_l - k)}}{1 - \chi^{- r_{l-1}}} + \frac{\chi^{ - r_l (\mu_l - k) + r_{l-1} \mu_{l+1} }}{1 - \chi^{r_{l-1}}}, $$
for any $ k = 1, \ldots , \mu_l $.

\item If $ l = 1 $, we get
$$ \Tr_{\beta}(e(H^{\bullet}(g|_{C_1},(j_1)^*(u_k^1)))) = \frac{1}{1 - \chi^{ - r_0 }} + \frac{ \chi^{ - r_1 (\mu_1 - k) + r_0 \mu_2}}{1 - \chi^{r_0}} $$ 
for any $ k = 1, \ldots , \mu_1 $. 

\item Finally, if $ l = L $, we get
$$ \Tr_{\beta}(e(H^{\bullet}(g|_{C_L},(j_L)^*(u_k^L)))) = \frac{\chi^{ r_{L-2}(\mu_L - k)}}{1 - \chi^{ - r_{L-1}}} + \frac{1}{1 - \chi^{r_{L-1}}} $$
for any $ k = 1, \ldots , \mu_L $.
\end{enumerate}
\end{prop}
\begin{pf}
Let us first assume that $ l \in \{2, \ldots, L - 1 \} $. We are going to apply Proposition \ref{prop. 5.6}. The only fixed points of the automorphism $ g|_{C_l} : C_l \rightarrow C_l $ are the two points where $C_l$ meets the other components, denoted as usual by $ y_{l-1}$ and $y_l$. From Lemma \ref{fiberlemma}, we have that
$$ \Tr((j_l)^*(u_k^l)(y_{l-1})) = \chi^{r_{l-2}(\mu_l - k)}, $$
and that
$$ \Tr((j_l)^*(u_k^l)(y_{l})) = \chi^{ - r_l (\mu_l - k) + r_{l-1} \mu_{l+1}}. $$

Let $ dg(y) : \Omega(y) \rightarrow \Omega(y) $ denote the induced automorphism of the cotangent space at $y$ on $ C_l $. It follows immediately from Proposition \ref{prop. 3.3} that $ \Tr(dg(y_{l-1})) = \chi^{-r_{l-1}} $, and that $ \Tr(dg(y_{l})) = \chi^{r_{l-1}} $. 

Proposition \ref{prop. 5.6} now gives that
$$ \Tr_{\beta}(e(H^{\bullet}(g|_{C_l},(j_l)^*(u_l^k)))) = \frac{\chi^{r_{l-2}(\mu_l - k)}}{1 - \chi^{- r_{l-1}}} + \frac{\chi^{ - r_l (\mu_l - k) + r_{l-1} \mu_{l+1} }}{1 - \chi^{r_{l-1}}},$$
which is precisely the sought after formula. 

The cases $ l = 1 $ and $ l = L $ are treated in the same fashion, when we recall that $ r_{-1} \equiv_n 0 $, and $ r_L := 0 $.
\end{pf}

\begin{ntn}
We will write $ \Tr_{\varepsilon}^{l , k}(\xi) := \Tr_{\beta}(e(H^{\bullet}(g|_{C_l},(j_l)^*(u_l^k)))) $. 
\end{ntn}

We can then make the following definition:

\begin{dfn}\label{contribution from edge} 
Let 
 \begin{equation} \Tr_{\varepsilon}(\xi) : = \sum_{l = 1}^{L(\varepsilon)} \sum_{k = 1}^{\mu_l} \Tr_{\varepsilon}^{l , k}(\xi). 
\end{equation}

We say that $ \Tr_{\varepsilon}(\xi) $ is the \emph{contribution} from $ \varepsilon \in \mathcal{E} $ to the trace 
$$ \Tr_{\beta}(e(H^{\bullet}(g|_{\mathcal{Y}_k})). $$
\end{dfn}

\begin{rmk}
Let us note that $ \Tr_{\varepsilon}(\xi) $ is defined entirely in terms of the intrinsic data of the singularity associated to $ \varepsilon $. Furthermore, this expression does not depend on the order in which we chose $ \varepsilon $. It is also clear that this expression does not depend on the chosen subfiltration of the divisor $ \sum_{i = 1}^{L(\varepsilon)} m_{\varepsilon_i} C_{\varepsilon_i} $. In the next section we shall define this trace \emph{intrinsically} for the singularity.
\end{rmk}

We will now show that we obtain a formula for the Brauer trace $ \Tr_{\beta}(e(H^{\bullet}(g|_{\mathcal{Y}_k}))) $, in terms of the vertex and edge contributions discussed above. 

\begin{thm}\label{thm. 9.13}
Let $ g \in G $ be a group element corresponding to a root of unity $\xi \in \boldsymbol{\mu}_n$. Then we have that
$$ \Tr_{\beta}(e(H^{\bullet}(g|_{\mathcal{Y}_k}))) = \sum_{\upsilon \in \mathcal{V}} \Tr_{\upsilon}(\xi) + \sum_{\varepsilon \in \mathcal{E}} \Tr_{\varepsilon}(\xi). $$ 
Furthermore, this expression depends only on $ \Gamma(\mathcal{X}_k) $ and the functions $ \mathfrak{g} $ and $ \mathfrak{m} $ defined in Section \ref{9.1}. 
 \end{thm}
\begin{pf}
We begin with choosing a special filtration
$$ 0 < \ldots < Z_{\varepsilon_i} < \ldots < Z_{\varepsilon_1} = Z_{\mathcal{E}} < \ldots < Z_{\upsilon_i} < \ldots < Z_{\upsilon_1} = \mathcal{Y}_k. $$
It then follows from Proposition \ref{prop. 6.6} that 
$$ \Tr_{\beta}(e(H^{\bullet}(g|_{\mathcal{Y}_k}))) = \sum_{\upsilon \in \mathcal{V}} \Tr_{\upsilon}(\xi) + \sum_{\varepsilon \in \mathcal{E}} \Tr_{\varepsilon}(\xi), $$
where $ \Tr_{\upsilon}(\xi) $ is the expression defined in Definition \ref{contribution from vertex} and $ \Tr_{\varepsilon}(\xi) $ is the expression defined in Definition \ref{contribution from edge}.

It follows from Proposition \ref{prop. 9.7} that $ \Tr_{\upsilon}(\xi) $ only depends on the combinatorial structure of $ \mathcal{X}_k $. Likewise, Proposition \ref{lemma 9.8} and Proposition \ref{lemma 9.10} give that $ \Tr_{\varepsilon}(\xi) $ only depends on the combinatorial structure of $ \mathcal{X}_k $. Therefore, the same is true for the sum of these expressions.
\end{pf}



\section{Trace formulas for singularities}\label{trace formula}

\subsection{}
In Section \ref{special filtrations}, given $ \epsilon \in \mathcal{E}(\mathcal{X}_k) $ and a tame extension $S'/S$ of degree $n$, we defined the expression $ \Tr_{\epsilon}(\xi) $ for any $ \xi \in \boldsymbol{\mu}_n $, and we were able to compute it using the global geometry of the surface $ \mathcal{Y} $, combined with the local description of the $ \boldsymbol{\mu}_n $-action on $ \mathcal{Y} $. However, the formula we obtain is expressed in terms of intrinsic data for the unique singularity of $ \mathcal{Y} $ associated to $ \epsilon $. So it makes sense to try to define this expression for any tame cyclic quotient singularity, without thinking of the global situation. Moreover, we will investigate the properties of the trace formulas, and in particular express them in a closed, polynomial form.

\subsection{}
Let $ \sigma = (m_1,m_2,n) $ be a singularity, as defined in Definition \ref{Def. 8.1}, and let  
$$ \rho_{\sigma} : \widetilde{\mathcal{Z}} \rightarrow \mathcal{Z} $$ 
be the minimal desingularization. We denote by $ Z $ the exceptional fiber $ \rho_{\sigma}^{-1}(z) $, where $ z \in \mathcal{Z} $ is the unique singular point of $\mathcal{Z}$. Then $ G = \boldsymbol{\mu}_n $ acts on $ \widetilde{\mathcal{Z}} $, and in particular on $ Z $. So every $ g \in G $ induces an automorphism $ e(H^{\bullet}(g|_Z)) $ of $ e(H^{\bullet}(Z, \mathcal{O}_Z)) $.

In order to compute the trace of $ e(H^{\bullet}(g|_Z)) $, it is easily seen that we can apply the same methods as in Section \ref{special filtrations}. Let $ L = L(\sigma) $ denote the length of the resolution chain, and let the irreducible components of $ \rho_{\sigma}^{-1}(z) $ be denoted $C_l$, and let ~$ \mu_l = \Mult(C_l) $. So again, we can choose a filtration
$$ 0 < Z_m < \ldots < Z_i < \ldots < Z_2 < Z_1 = Z, $$ 
where $ m = \sum_{l=1}^{L} \mu_l $, and $ Z_i - Z_{i+1} = C_{l_i} $, for every $i$. At the $i$-th step in this filtration, there is as usual associated an invertible sheaf $ \mathcal{L}_i $, supported on $ C_{l_i} $, and for every $ g \in G $, an isomorphism 
$$ u_i : g^* (j_{C_{l_i}})_* \mathcal{L}_i \rightarrow (j_{C_{l_i}})_* \mathcal{L}_i, $$  
where $ j_{C_{l_i}} : C_{l_i} \hookrightarrow \widetilde{\mathcal{Z}} $ is the canonical inclusion. This set of data gives, for every $i$, an automorphism $ e(H^{\bullet}(g|_{C_{l_j}}, (j_{C_{l_i}})^*u_i)) $ such that
$$ \Tr(e(H^{\bullet}(g|_Z))) = \sum_{j=1}^m  \Tr(e(H^{\bullet}(g|_{C_{l_j}}, (j_{C_{l_i}})^*u_i))). $$

The trace is independent of which filtration we choose. So in particular, we can use the \emph{special} filtration
$$ 0 =: Z_{L+1} < Z_L < \ldots < Z_l < \ldots < Z_1 := Z, $$
where the $Z_l$ are defined by the refinements
$$ Z_{l+1} := Z_l^{\mu_l + 1} < \ldots < Z_l^k < \ldots < Z_l^1 := Z_l, $$
and $ Z_l^{k} - Z_l^{k+1} = C_l $, for every $ k \in \{ 1, \ldots, \mu_l \} $.

We immediately get the following result:

\begin{prop}Let $ g \in G $ be a group element corresponding to a root of unity $ \xi $. Let us write $ \Tr_{\sigma}(\xi) := \Tr(e(H^{\bullet}(g|_Z))) $ and $ \Tr_{\sigma}^{l,k} (\xi) := \Tr(e(H^{\bullet}(g|_{C_l}, (j_{C_l})^* u_l^k))) $. Then we have that
$$ \Tr_{\sigma}(\xi) = \sum_{l=1}^L \sum_{k=1}^{\mu_l} \Tr_{\sigma}^{l,k} (\xi). $$
\end{prop} 


\begin{rmk}
The connection to Section \ref{special filtrations} is as follows: Let $ \varepsilon \in \mathcal{E}(\mathcal{X}_k) $ be an edge, and let $ \sigma = (m_1,m_2,n) $ be the singularity of $ \mathcal{Y} $ corresponding to $ \varepsilon $, after a base change $S'/S$ of degree $n$. Let $ \xi \in \boldsymbol{\mu}_n $. Then we have that $ \Tr_{\sigma}(\xi) = \Tr_{\varepsilon}(\xi) $.
\end{rmk}

We are now going to investigate closer the terms in the formula 
$$ \Tr_{\sigma}(\xi) = \sum_{l=1}^L \sum_{k=1}^{\mu_l} \Tr_{\sigma}^{l,k} (\xi). $$




\begin{lemma}\label{sigmasing}
Consider a singularity $ \sigma = (m_1,m_2,n) $, and let $ \xi \in \boldsymbol{\mu}_n $ be a primitive root of unity. Let $ \chi = \xi^{\alpha_1} $, where $ \alpha_1 m_1 \equiv_n 1 $. Then we have that
\begin{enumerate} 
\item $$ \sum_{k=1}^{\mu_l} \Tr_{\sigma}^{l,k}(\xi) = \sum_{k=1}^{\mu_l} \frac{ \chi^{r_{l-2}(\mu_l - k)}}{1 - \chi^{- r_{l-1}}} + \sum_{k=1}^{\mu_l} \frac{ \chi^{ - r_l (\mu_l - k) + r_{l-1} \mu_{l+1} }}{1 - \chi^{r_{l-1}}}, $$
for any $ l \in \{2, \ldots,  L - 1\} $, that
 \item $$ \sum_{k=1}^{\mu_1} \Tr_{\sigma}^{1,k}(\xi) = \sum_{k=1}^{\mu_1} \frac{1}{1 - \chi^{ - r_0 }} + \sum_{k=1}^{\mu_1} \frac{ \chi^{ - r_1 (\mu_1 - k) + r_0 \mu_2}}{1 - \chi^{r_0}}, $$
and that 
\item $$ \sum_{k=1}^{\mu_L} \Tr_{\sigma}^{L,k}(\xi) = \sum_{k=1}^{\mu_L} \frac{\chi^{ r_{L-2}(\mu_L - k)}}{1 - \chi^{ - r_{L-1}}} + \sum_{k=1}^{\mu_L} \frac{1}{1 - \chi^{r_{L-1}}}. $$
\end{enumerate}
\end{lemma}
\begin{pf}
The proof is similar to the proof of Lemma \ref{lemma 9.10}.
\end{pf}

\subsection{Formal manipulations}
Consider a singularity $ \sigma = (m_1,m_2,n) $. We have, for any $ \xi \in \boldsymbol{\mu}_n $, a trace $ \Tr_{\sigma}(\xi) = \sum_{l=1}^L \sum_{k=1}^{\mu_l} \Tr_{\sigma}^{l,k} (\xi) $. We would like to have an explicit formula for this sum. In particular, we would like to know if it is a polynomial in $\xi$, and in which way it depends on the parameters of the singularity. In order to deal with these questions, and to obtain a closed formula, some formal manipulations of the expressions in Lemma \ref{sigmasing} will be necessary. Unless otherwise mentioned, $ \xi $ will be a primitive root of unity.

The following lemma is an easy computation, whose proof is omitted.
\begin{lemma}\label{Lemma 10.3}
We have that 
$$ \chi^{ r_{l-1} \mu_{l+1} } \sum_{k=1}^{\mu_l} \frac{ \chi^{ - r_l (\mu_l - k) }}{1 - \chi^{r_{l-1}}} + \sum_{k=1}^{\mu_{l+1}} \frac{ \chi^{r_{l-1} (\mu_{l+1} - k)}}{1 - \chi^{- r_l}} = \frac{1 - \chi^{ r_{l-1} \mu_{l+1} - r_l \mu_l}}{(1 - \chi^{r_{l-1}})(1 - \chi^{- r_l})}, $$
for $ l \in \{ 1, \ldots, L-1 \}$.
\end{lemma}

Let us make the following definition:

\begin{dfn} We define
$$ \Tr_{y_l}(\xi) := \frac{1 - \chi^{ r_{l-1} \mu_{l+1} - r_l \mu_l}}{(1 - \chi^{r_{l-1}})(1 - \chi^{- r_l})}, $$
for every $ l \in \{ 1, \ldots, L-1 \}$.
We also define 
$$ \Tr_{y_0}(\xi) := \frac{\mu_1}{1 - \chi^{-r_0}} $$
and
$$ \Tr_{y_L}(\xi) := \frac{\mu_L}{1 - \chi^{r_{L-1}}}. $$
\end{dfn}

\vspace{0.5cm}

Note that it follows that $ \Tr_{\sigma}(\xi) = \sum_{l=0}^L \Tr_{y_l}(\xi) $. We shall find a more convenient way to write $ \Tr_{y_l} $. The first step is Lemma \ref{lemma 10.5} below, whose proof is omitted.

\begin{lemma}\label{lemma 10.5}
For every $ l \in \{ 1, \ldots, L-1 \}$, we have
$$ 1 - \chi^{ r_{l-1} \mu_{l+1} - r_l \mu_l} = $$
$$ (1 - \chi^{- r_l \mu_l})(1 - \chi^{r_{l-1} \mu_{l+1}}) + \chi^{r_{l-1} \mu_{l+1}} (1 - \chi^{- r_l \mu_l}) + \chi^{- r_l \mu_l}(1 - \chi^{r_{l-1} \mu_{l+1}}). $$  
\end{lemma}

\vspace{0.5cm}

\begin{cor}\label{cor. 10.6}
For every $ l \in \{ 1, \ldots, L-1 \}$, we can write $ \Tr_{y_l} $ in the following form:
$$ \Tr_{y_l}(\xi) = \sum_{k=0}^{\mu_l - 1} (\chi^{- r_l})^k \sum_{k=0}^{\mu_{l+1} - 1} (\chi^{r_{l-1}})^k   +   $$
$$ \frac{\chi^{ r_{l-1} \mu_{l+1}}}{1 - \chi^{ r_{l-1}}} \sum_{k=0}^{\mu_l - 1} (\chi^{- r_l})^k  + \frac{\chi^{- r_l \mu_l}}{1 - \chi^{- r_l}} \sum_{k=0}^{\mu_{l+1} - 1} (\chi^{r_{l-1}})^k. $$
\end{cor}

We introduce some notation for the terms appearing in $ \Tr_{y_l}(\xi) $. 
\begin{dfn}\label{dfn 10.7} We define
\begin{enumerate}
\item $$ \Tr_{\mu_l, \mu_{l+1}}(\xi) := \sum_{k=0}^{\mu_l - 1} (\chi^{- r_l})^k \sum_{k=0}^{\mu_{l+1} - 1} (\chi^{r_{l-1}})^k, $$ 
for all $ l \in \{ 0, \ldots, L \} $,
\item $$ \Tr_{y_l}^{\mu_l}(\xi) := \frac{\chi^{ r_{l-1} \mu_{l+1}}}{1 - \chi^{ r_{l-1}}} \sum_{k=0}^{\mu_l - 1} (\chi^{- r_l})^k, $$
for all $ l \in \{ 1, \ldots, L \} $, and 
\item $$ \Tr_{y_l}^{\mu_{l+1}}(\xi) := \frac{\chi^{- r_l \mu_l}}{1 - \chi^{- r_l}} \sum_{k=0}^{\mu_{l+1} - 1} (\chi^{r_{l-1}})^k, $$
for all $ l \in \{ 0, \ldots, L-1 \} $.
\end{enumerate}
\end{dfn}

\begin{lemma}\label{lemma blabla} For any $ l $ such that $ 0 \leq l \leq L - 1 $, we have that
$$ \Tr_{y_l}^{\mu_{l+1}}(\xi) + \Tr_{y_{l + 1}}^{\mu_{l+1}}(\xi) = - \sum_{k=0}^{\mu_{l+1} - 1} \chi^{r_{l-1} k - r_l (\mu_l - 1)} (1 + (\chi^{r_{l}}) + \ldots + (\chi^{r_l})^{b_{l+1}(\mu_{l+1} - k) - 2}). $$ 

Furthermore, we may also write this sum as 
$$ \sum_{k=0}^{\mu_{l+1} - 1} \chi^{r_{l-1} k} ((\chi^{ - r_{l}})^{ \mu_l - 1 } + (\chi^{ - r_{l}})^{ \mu_l - 2 } + \ldots + (\chi^{ - r_{l}})^{b_{l+1} k - ( \mu_{l + 2} - 1 )}). $$
\end{lemma}
\begin{pf}
Observe that we have 
$$ (\sum_{k=0}^{\mu_{l+1} - 1} (\chi^{r_{l-1}})^k) \frac{\chi^{- r_l \mu_l}}{1 - \chi^{- r_l}} + (\sum_{k=0}^{\mu_{l+1} - 1} (\chi^{ - r_{l+1}})^k) \frac{\chi^{ r_{l} \mu_{l+2}}}{1 - \chi^{ r_{l}}} $$
$$ = \frac{1}{1 - \chi^{ r_{l}}} \sum_{k=0}^{\mu_{l+1} - 1} (\chi^{ - r_{l+1} k + r_{l} \mu_{l+2} } - \chi^{r_{l-1} k - r_l (\mu_l - 1 ) } ). $$ 
We can rewrite $ \chi^{ - r_{l+1} k + r_{l} \mu_{l+2} } - \chi^{r_{l-1} k - r_l (\mu_l - 1 ) } $ as 
$$ - \chi^{r_{l-1} k - r_l (\mu_l - 1 ) } ( 1 - \chi^{  - r_{l+1} k + r_{l} \mu_{l+2} - r_{l-1} k + r_l (\mu_l - 1 ) }). $$
We then compute that
$$ - r_{l+1} k + r_{l} \mu_{l+2} - r_{l-1} k + r_l (\mu_l - 1 ) = r_l (b_{l+1} ( \mu_{l+1} - k ) - 1), $$
since $ \mu_l + \mu_{l+2} = b_{l+1} \mu_{l+1} $ and $ r_{l-1} = b_{l+1} r_{l} - r_{l+1} $.

This gives that 
$$ \chi^{ - r_{l+1} k + r_{l} \mu_{l+2} } - \chi^{r_{l-1} k - r_l (\mu_l - 1 ) } = - \chi^{r_{l-1} k - r_l (\mu_l - 1 ) } ( 1 - (\chi^{r_l})^{b_{l+1} ( \mu_{l+1} - k ) - 1}). $$
Recall that $ b_{l+1} \geq 2 $, so for all $ 0 \leq k \leq \mu_{l+1} - 1 $, we have that 
$$ b_{l+1} ( \mu_{l+1} - k ) - 1 \geq 1. $$ 
But then it follows that 
$$ \frac{1 - (\chi^{r_l})^{b_{l+1} ( \mu_{l+1} - k ) - 1}}{1 - \chi^{ r_{l}}} = 1 + (\chi^{r_l}) + \ldots + (\chi^{r_l})^{b_{l+1} ( \mu_{l+1} - k ) - 2}. $$

The last statement follows from observing that 
$$ (\mu_l - 1) - (b_{l+1} ( \mu_{l+1} - k ) - 2) = b_{l+1} k - ( \mu_{l + 2} - 1 ). $$
\end{pf}

\vspace{0.3cm}

It is convenient to introduce some notation for the terms of this form.

\begin{dfn} For any $l$ such that $ 0 \leq l \leq L-1 $, we define
$$ \Tr_{\mu_{l+1}}(\xi) := \Tr_{y_l}^{\mu_{l+1}}(\xi) + \Tr_{y_{l + 1}}^{\mu_{l+1}}(\xi). $$ 
\end{dfn}

It is now possible to express $ \Tr_{\sigma}(\xi) $ in the following form:

\begin{prop}\label{prop. 10.14}
Let $ \sigma = (m_1, m_2, n) $ be a singularity, and let $ \xi \in \boldsymbol{\mu}_n $ be a primitive root. Then we have that
$$ \Tr_{\sigma}(\xi) = \sum_{l=0}^{L} \Tr_{\mu_l, \mu_{l+1}}(\xi) + \sum_{l=1}^{L} \Tr_{\mu_{l}}(\xi). $$
\end{prop}
\begin{pf}
This is rather immediate, with the exception of the appearance of the terms $ \Tr_{\mu_0, \mu_1} $, $ \Tr_{\mu_1} $, $ \Tr_{\mu_L} $ and $ \Tr_{\mu_L,\mu_{L+1} } $ in the formula. But an easy computation shows that
$$ \sum_{k=0}^{\mu_1 - 1} \frac{1}{1 - \chi^{ - r_0 }} + \Tr_{y_1}^{\mu_1} = (\sum_{k=0}^{\mu_1 - 1} \frac{1}{1 - \chi^{ - r_0 }} - \Tr_{y_0}^{\mu_1}) + (\Tr_{y_0}^{\mu_1} + \Tr_{y_1}^{\mu_1}) = \Tr_{\mu_0, \mu_1} + \Tr_{\mu_1}. $$ 
In a similar way, we compute
$$ \Tr_{y_{L-1}}^{\mu_L} + \sum_{k=1}^{\mu_L} \frac{1}{1 - \chi^{r_{L-1}}} = \Tr_{\mu_L} + \Tr_{\mu_L,\mu_{L+1} }. $$
\end{pf}

Furthermore, the following is now automatic from our description:

\begin{cor}\label{nonprimformula} 
Let $ \xi \in G $ be any root of unity. Then we have that
$$ \Tr_{\sigma}(\xi) = \sum_{l=0}^{L} \Tr_{\mu_l, \mu_{l+1}}(\xi) + \sum_{l=1}^{L} \Tr_{\mu_{l}}(\xi). $$
\end{cor}
\begin{pf}
Since $p$ does not divide the order of $ \boldsymbol{\mu}_n $, it follows that the action of $ \boldsymbol{\mu}_n $ on $ H^i(Z, \mathcal{O}_Z) $, where $ i \in \{0,1\} $, is diagonalizable. In particular, the irreducible Brauer characters are all of the form $ \xi \mapsto \xi^j $, for some $ j \geq 0 $. Therefore, the Brauer trace is given by the same polynomial for any $ \xi \in \boldsymbol{\mu}_n $. 
\end{pf}

\section{Trace formula}\label{Trace formula}
Let $ \sigma = (m_1,m_2,n) $ be a singularity. In this section, we will prove an explicit formula for $\Tr_{\sigma}(\xi)$, under the assumption that $ n \gg 0 $. We will see that the ``shape'' of this formula is closely related to the properties of the exceptional locus of the minimal desingularization of $ \sigma $.   

In Proposition \ref{prop. 8.9}, we saw that the exceptional fiber of the minimal desingularization of $ \sigma = (m_1, m_2, n) $, for $ n \gg 0 $, only depended on the residue class $ [n]_{M} $ of $ n $ modulo $ M = \Lcm(m_1,m_2)$, modulo chains of curves with multiplicity $ m = \Gcd(m_1,m_2) $. More precisely, the multiplicities satisfy:
$$ m_2 = \mu_0 > \mu_1 > \ldots > \mu_{l_0} = \ldots = m = \ldots = \mu_{L+1-l_1} < \ldots < \mu_L < \mu_{L+1} = m_1. $$
The integer $l_0$ and the sequence of multiplicities $ \mu_1, \ldots, \mu_{l_0} $, as well as the integer $l_1$ and the sequence of multiplicities $ \mu_{L+1-l_1}, \ldots,  \mu_L $ depend only on $ [n]_{M} $. 

In fact, we shall see that $\Tr_{\sigma}(\xi)$ can be written as a polynomial, where the degree and the coefficients only depend on $ \mu_0 $, $ \mu_1 $, $ \mu_L $ and $ \mu_{L+1} $. Since $ \mu_0 $ and $ \mu_{L+1} $ are fixed, and $ \mu_1 $ and $ \mu_L $ only depend on the residue class of $n$ modulo $M$, it will follow that this also holds for $\Tr_{\sigma}(\xi)$. 

As a further motivation for what we will do in this section, we remark that as $n$ goes to infinity (but with a fixed residue class modulo $M$), any $ \xi \in \boldsymbol{\mu}_n $ behaves ``less'' like a root of unity, and more like an independent variable. This suggests that the cancellations occuring in our formulas are of a formal nature, and that we should substitute $ \xi $ with a variable.   


\subsection{Formal substitution}
For any $ n \gg 0 $ with a fixed residue class $ [n]_{M} $, the self intersection numbers $ b_1, \ldots, b_{l_0} $ are constant, since they are computed in terms of the multiplicities via the formula $ b_l \mu_l = \mu_{l-1} + \mu_{l+1} $. However, the integers $ r_l $ may vary as $n$ varies. But on the other hand, these integers are related in terms of the $b_l$ via the equations $ r_{l-1} = b_{l+1} r_l - r_{l+1} $. 

We will now define universal polynomials $ P_l $ inductively, by the following procedure: Put $ P_{-1} = 0 $, and $ P_0 = 1 $. Then we define $ P_l = b_l P_{l-1} - P_{l-2} $, for $ l \geq 1 $. Note that $ P_l = P_l(b_1, \ldots, b_l) $, when $ l \geq 1 $. For instance, we have $ P_1 = b_1 $, $ P_2 = b_2 b_1 - 1 $, $ P_3 = b_3(b_2b_1 - 1) - b_1 $, and so on.

The importance of these polynomials is that $ r_l \equiv_n  P_l r_0 $. So for any $n$-th root of unity $\eta$, we have that $ \eta^{r_l} = \eta^{P_l r_0} $. Recall that $ \chi := \xi^{\alpha_1} $, and that $ - \alpha_1 r_0 \equiv_n \alpha_2 $. It follows that we can write
\begin{equation}\label{equation 10.5}
\chi^{r_l} = \xi^{\alpha_1 r_l} = \xi^{ \alpha_1 r_0 P_l} = (\xi^{ \alpha_2})^{ - P_l}. 
\end{equation}

\subsection{}\label{lemmas}

Recall from Proposition \ref{prop. 10.14} that we had
$$ \Tr_{\sigma}(\xi) = \sum_{l=1}^{L+1} \Tr_{\mu_{l-1}, \mu_{l}}(\xi) + \sum_{l=1}^{L} \Tr_{\mu_l}(\xi), $$  
for any $ \xi \in \boldsymbol{\mu}_n $. It turns out that a small reformulation of the expressions $ \Tr_{\mu_{l-1}, \mu_{l}}(\xi) $ and $ \Tr_{\mu_l}(\xi) $ is suitable for our computations later on. Let us also here remark that we will drop the reference to $ \xi $ from the notation.

Let $ q \in \{ 1, \ldots, L \} $ be an integer such that $ \mu_q = m $. For all $ l \in \{1, \ldots, q \} $ such that $ \mu_{l+1} > 1 $, we then define:
\begin{equation}
\Tr_{l-1,l} := \sum_{k=1}^{\mu_{l} - 1} (\chi^{r_{l-2}})^k ( (\chi^{-r_{l-1}})^{b_l k -1} + \ldots + 1), 
\end{equation}
and
\begin{equation}
\Tr_l := - \sum_{k=0}^{\mu_{l} - 1} (\chi^{r_{l-2}})^k ( (\chi^{-r_{l-1}})^{b_l k -1} + \ldots + (\chi^{-r_{l-1}})^{b_l k - (\mu_{l+1} - 1)}). 
\end{equation}

We can replace $ \Tr_{\mu_{l-1}, \mu_{l}} + \Tr_{\mu_l} $ with $ \Tr_{l-1,l} + \Tr_l $ in the expression for $ \Tr_{\sigma}(\xi) $, as the lemma below shows. 

\begin{lemma}\label{Lemma 11.3}
For all $ l \in \{1, \ldots, q \} $ such that $ \mu_{l+1} > 1 $, we have that 
$$ \Tr_{\mu_{l-1}, \mu_{l}} + \Tr_{\mu_l} = \Tr_{l-1,l} + \Tr_l. $$
\end{lemma}
\begin{pf}
We will write $ \Tr_{\mu_l}^{(k)} $ for the $k$-th summand in the expression for $ \Tr_{\mu_l} $, and likewise for $ \Tr_{\mu_{l-1}, \mu_{l}}^{(k)} $, $ \Tr_l^{(k)} $ and $ \Tr_{l-1,l}^{(k)} $.

Let us first consider the case where $ b_l k - 1 > \mu_{l-1} - 1 $. Then we have that 
$$ \Tr_{l-1,l}^{(k)} = (\chi^{r_{l-2}})^k ( (\chi^{-r_{l-1}})^{b_l k -1} + \ldots + 1) $$
and
$$ \Tr_l^{(k)} = - (\chi^{r_{l-2}})^k ( (\chi^{-r_{l-1}})^{b_l k -1} + \ldots + (\chi^{-r_{l-1}})^{b_l k - (\mu_{l+1} - 1)}), $$
that is, we add some monomials to $ \Tr_{\mu_{l-1}, \mu_{l}}^{(k)} $ and subtract exactly the same monomials from $ \Tr_{\mu_l}^{(k)} $. It follows easily that
$$ \Tr_{\mu_{l-1}, \mu_{l}}^{(k)} + \Tr_{\mu_l}^{(k)} = \Tr_{l-1,l}^{(k)} + \Tr_l^{(k)}. $$

For $ k = 0 $, one easily computes that 
$$ \Tr_{\mu_l}^{(0)} = - \Tr_{\mu_{l-1}, \mu_{l}}^{(0)} + \Tr_l^{(0)}, $$
and for all $ k \geq 1 $ such that $ b_l k \leq \mu_{l-1} - 1 $, it is trivial to see that 
$$ \Tr_{\mu_{l-1}, \mu_{l}}^{(k)} + \Tr_{\mu_l}^{(k)} = \Tr_{l-1,l}^{(k)} + \Tr_l^{(k)}. $$


\end{pf}

\vspace{0.5cm}

We will now replace $ \chi^{-r_0} $ in the expressions $ \Tr_{l-1,l} $ and $ \Tr_l $ above with a formal variable $y$. There are two reasons for doing this. First, we can write these expressions in a compact form as polynomials in $y$. Second, when we consider various sums of these expressions, it is manageable to keep track of the formal cancellations that occur.

\begin{lemma}\label{Omskriving} Let $ l \in \{ 1, \ldots, q \} $ be such that $ \mu_{l+1} > 1 $. Put $ y = \chi^{-r_0} $. Then we can write 
\begin{enumerate}
\item $$ \Tr_{l-1,l} = \sum_{k=1}^{\mu_l - 1} \sum_{m=1}^{b_l k} y^{k P_l - m P_{l-1}}, $$
and
\item $$ \Tr_l = - \sum_{k=0}^{\mu_l - 1} \sum_{m=1}^{\mu_{l+1} - 1} y^{k P_l - m P_{l-1}}. $$
\end{enumerate}
\end{lemma}
\begin{pf}
We have in case (i) that
$$ \Tr_{l-1,l} = \sum_{k=1}^{\mu_l - 1} (\chi^{r_{l-2}})^k ( (\chi^{-r_{l-1}})^{b_l k - 1} + \ldots + 1) = $$
$$ \sum_{k=1}^{\mu_l - 1} y^{-P_{l-2} k} (y^{b_l P_{l-1} k - P_{l-1}} + \ldots + y^{b_l P_{l-1} k - b_l k P_{l-1}}) =$$
$$ \sum_{k=1}^{\mu_l - 1} \sum_{m=1}^{b_l k} y^{(b_l P_{l-1} - P_{l-2})k - m P_{l-1}} = \sum_{k=1}^{\mu_l - 1} \sum_{m=1}^{b_l k} y^{k P_l - m P_{l-1}}. $$

The proof of (ii) is similar, and is omitted.
\end{pf}

\subsection{}
This section consists of three lemmas that we will use when proving the trace formula in Theorem \ref{Formula}. The reader might want to skip this section for now, and refer back to these results when needed in Section \ref{star} and Section \ref{explicit formula}.

\begin{lemma}\label{Lemma 11.4}
Let $ \sigma = (m_1,m_2,n) $ be a singularity where $ n \gg 0 $, and let $ L $ be the length of the resolution of $\sigma$. Let $ q \in \{1, \ldots, L\} $ be such that $ \mu_q = m $ and choose an integer $ l \in \{1, \ldots, L\} $ such that $ l + 1 \leq q $. Consider the inequality
$$ (*) ~ b_l t - \left \lceil t \frac{\mu_{l+1}}{\mu_l} \right \rceil \geq s, $$
where $ 1 \leq s \leq \mu_{l-1} $, and where $ t \geq 0 $ is an integer. Then we have that $ t(s) := \left \lceil s \frac{\mu_{l}}{\mu_{l-1}} \right \rceil $ is the smallest positive integer that satisfies the inequality $(*)$. Note in particular that $ 1 \leq t(s) \leq \mu_l $.

\end{lemma}
\begin{pf}
Let us first show that $ t(s) $ satisfies the inequality. By definition we have that $ s \frac{\mu_{l}}{\mu_{l-1}} \leq t(s) $, so it follows that $ s \leq t(s) \frac{\mu_{l-1}}{\mu_{l}} $. From the equality $ b_l \mu_l = \mu_{l+1} + \mu_{l-1} $ it follows that $ \frac{\mu_{l-1}}{\mu_{l}} = b_l - \frac{\mu_{l+1}}{\mu_{l}} $. So we get that 
$$ s \leq t(s)(b_l - \frac{\mu_{l+1}}{\mu_{l}}) = t(s) b_l - t(s) \frac{\mu_{l+1}}{\mu_{l}}. $$
From this, it follows that $ t(s) \frac{\mu_{l+1}}{\mu_{l}} \leq t(s) b_l - s $, and since $ t(s) $, $ b_l $ and $s$ are integers, we actually have that $ \left \lceil t(s) \frac{\mu_{l+1}}{\mu_{l}} \right \rceil \leq t(s) b_l - s $, and hence $t(s)$ satisfies $ (*) $.

Assume now that $ t \geq 0 $ is an integer such that $ b_l t - \left \lceil t \frac{\mu_{l+1}}{\mu_l} \right \rceil \geq s $, and that $ t < t(s) $. This implies that $ t < s \frac{\mu_{l}}{\mu_{l-1}} $, and so it follows that 
$$ t \mu_{l-1} < s \mu_{l} \leq \mu_{l} (b_l t - \left \lceil t \frac{\mu_{l+1}}{\mu_l} \right \rceil) = t \mu_{l-1} + t \mu_{l+1} - \mu_{l} \left \lceil t \frac{\mu_{l+1}}{\mu_l} \right \rceil. $$

So in fact, we get that $ \mu_{l} \left \lceil t \frac{\mu_{l+1}}{\mu_l} \right \rceil < t \mu_{l+1} $, and hence $ \left \lceil t \frac{\mu_{l+1}}{\mu_l} \right \rceil < t \frac{\mu_{l+1}}{\mu_l} $, a contradiction.
\end{pf}

\begin{lemma}\label{Lemma Key2}
Let $ \sigma = (m_1,m_2,n) $ be a singularity where $ n \gg 0 $. Let $ q \in \{1, \ldots, L\} $ be such that $ \mu_q = m $, and choose $ l \in \{ 0, \ldots, L-1 \} $ such that $ l + 2 \leq q $. Let us assume that $ \mu_{l+1} \geq 2 $. For all integers $ s $ and $ k $ such that $ 1 \leq s < k \leq \mu_{l+1} $, the strict inequality 
$$ b_{l+1} k - \left \lceil k \frac{\mu_{l+2}}{\mu_{l+1}} \right \rceil > b_{l+1} s - \left \lceil s \frac{\mu_{l+2}}{\mu_{l+1}} \right \rceil $$
holds. In particular, by taking $ k = \mu_{l+1} $, we get that 
$$ \mu_l > b_{l+1} s - \left \lceil s \frac{\mu_{l+2}}{\mu_{l+1}} \right \rceil, $$
for all $ s \in \{ 1, \ldots, \mu_{l+1} - 1 \} $.
\end{lemma}
\begin{pf}
We first note that $ \mu_{l+1} \geq \mu_{l+2} $. Observe that
$$ b_{l+1} k - b_{l+1} s = b_{l+1} (k - s) \geq 2 (k - s), $$
and likewise
$$ k \frac{\mu_{l+2}}{\mu_{l+1}} - s \frac{\mu_{l+2}}{\mu_{l+1}} = (k-s) \frac{\mu_{l+2}}{\mu_{l+1}} \leq (k-s). $$
It follows that
$$ k \frac{\mu_{l+2}}{\mu_{l+1}} \leq (k-s) + s \frac{\mu_{l+2}}{\mu_{l+1}} \leq (k-s) + \left \lceil s \frac{\mu_{l+2}}{\mu_{l+1}} \right \rceil. $$
Hence
$$ \left \lceil k \frac{\mu_{l+2}}{\mu_{l+1}} \right \rceil - \left \lceil s \frac{\mu_{l+2}}{\mu_{l+1}} \right \rceil \leq (k-s). $$
By combining these inequalities, we get
$$ \left \lceil k \frac{\mu_{l+2}}{\mu_{l+1}} \right \rceil - \left \lceil s \frac{\mu_{l+2}}{\mu_{l+1}} \right \rceil \leq (k-s) < 2 (k-s) \leq b_{l+1} k - b_{l+1} s, $$
and therefore
$$ b_{l+1} k - \left \lceil k \frac{\mu_{l+2}}{\mu_{l+1}} \right \rceil > b_{l+1} s - \left \lceil s \frac{\mu_{l+2}}{\mu_{l+1}} \right \rceil. $$
\end{pf}

\begin{lemma}\label{Lemma 11 Key}
Let $ \sigma = (m_1,m_2,n) $ be a singularity where $ n \gg 0 $. Let $ q \in \{1, \ldots, L\} $ be such that $ \mu_q = m $, and choose $ l \in \{ 0, \ldots, L-1 \} $ such that $ l + 2 \leq q $. Consider the inequalities $ \mu_l \geq \mu_{l+1} \geq \mu_{l+2} $. Let us assume that either  
\begin{enumerate}
\item $ \mu_{l+1} > \mu_{l+2} $, or that
\item $ \mu_{l+1} = \mu_{l+2} \geq 3 $. 
\end{enumerate}
Then the inequality $ \mu_l - 1 \geq b_{l+1} $ holds.
\end{lemma}
\begin{pf}
The equality $ \mu_l = b_{l+1} \mu_{l+1} - \mu_{l+2} $ can be written as 
$$ \mu_l = b_{l+1} + b_{l+1} (\mu_{l+1} - 1) - \mu_{l+2}. $$  
Let us first consider case (i), where $  \mu_{l+1} > \mu_{l+2} $. Then $ \mu_{l+1} - 1 \geq \mu_{l+2} $, and $ b_{l+1} \geq 2 $, so 
$$ b_{l+1} (\mu_{l+1} - 1) - \mu_{l+2} \geq 2 \mu_{l+2} - \mu_{l+2} = \mu_{l+2} \geq 1, $$
and therefore $ \mu_l \geq b_{l+1} + 1 $, as desired.

In case (ii), where $ \mu_{l+1} = \mu_{l+2} \geq 3 $, we see that
$$ b_{l+1} (\mu_{l+1} - 1) - \mu_{l+2} \geq 2 (\mu_{l+1} - 1) - \mu_{l+1} = \mu_{l+1} - 2 \geq 1, $$
and so again, it follows that $ \mu_l \geq b_{l+1} + 1 $.
\end{pf}

\subsection{Formal cancellation}\label{star} 
We will now consider sums of the expressions $ \Tr_{l-1,l} $ and $ \Tr_{l} $. It turns out that there will be cancellations occurring in these sums, following a certain pattern. We will eventually, in Proposition \ref{Computation}, end up with a polynomial in $y$ derived from $ \Tr_{0,1} $, which we can compute precisely.    

To set this up, let $ q \in \{ 1, \ldots , L-1 \} $ be an integer such that $ \mu_q = m $, and assume that $ \lambda \in \{ 0, \ldots , q \} $ is an index such that either $ \mu_{\lambda + 1} > \mu_{\lambda + 2} $ holds or that $ \mu_{\lambda + 1} = \mu_{\lambda + 2} \geq 3 $ holds. From Lemma \ref{Lemma 11 Key} above, it follows that $ \mu_l - 1 \geq b_{l+1} $ for all $ l \in \{ 0, \ldots, \lambda \} $, and hence it makes sense to define
\begin{equation} 
\Tr_{l}^* := - \sum_{k = b_{l+1}}^{\mu_{l} - 1} \sum_{m=1}^{\left \lceil k \frac{\mu_{l+1}}{\mu_{l}} \right \rceil - 1} y^{k P_{l}  - m P_{l-1} }, 
\end{equation}
and
\begin{equation} 
\Tr_{l-1,l}^* := \sum_{k = 1}^{\mu_{l} - 1} \sum_{m = \left \lceil k \frac{\mu_{l+1}}{\mu_{l}} \right \rceil}^{b_l k} y^{k P_{l}  - m P_{l-1}. } 
\end{equation}

\begin{prop}\label{Prop Key} 
Let us keep the assumptions above. For any index $l$ such that $ 0 \leq l < \lambda $, the following identities hold:
\begin{enumerate}
\item $ \Tr_l^* = \Tr_l + \Tr_{l,l+1}^* $,
and
\item $ \Tr_{l-1,l}^* = \Tr_{l-1,l} + \Tr_l^* $. 
\end{enumerate}
\end{prop}
\begin{pf}
By definition, we have that
$$ \Tr_{l, l + 1}^* = \sum_{s = 1}^{\mu_{l + 1} - 1} \sum_{t = \left \lceil s \frac{\mu_{l + 2}}{\mu_{l + 1}} \right \rceil}^{b_{l + 1} s } y^{s P_{l + 1}  - t P_{l}}. $$
Furthermore, it follows from Lemma \ref{Omskriving} that
$$ \Tr_{l} = - \sum_{k = 0}^{\mu_{l} - 1} \sum_{m=1}^{\mu_{l + 1} - 1} y^{k P_{l}  - m P_{l - 1} }, $$
and it is easily seen by re-indexing the expression for $ \Tr_{l, l + 1}^* $ that we may write 
$$ \Tr_{l, l + 1}^* = \sum_{m=1}^{\mu_{l + 1} - 1} \sum_{k = 0}^{b_{l + 1} m - \left \lceil m \frac{\mu_{l + 2}}{\mu_{l + 1}} \right \rceil} y^{k P_{l}  - m P_{l - 1} }. $$

Let us put
$$ Q_k := - \sum_{m=1}^{\mu_{l + 1} - 1} y^{k P_{l}  - m P_{l - 1} }, $$
and 
$$ S_m := \sum_{k = 0}^{b_{l + 1} m - \left \lceil m \frac{\mu_{l + 2}}{\mu_{l + 1}} \right \rceil} y^{k P_{l}  - m P_{l - 1} }. $$

Consider now $S_m$ for a fixed $ m \in \{1, \ldots, \mu_{l+1} - 1 \} $. We let the monomials in $S_m$ cancel monomials in $Q_k$ in the following systematic way: We let the $k$-th term $ y^{k P_{l}  - m P_{l - 1} } $ in $S_m$, where $ 0 \leq k \leq b_{l + 1} m - \left \lceil m \frac{\mu_{l + 2}}{\mu_{l + 1}} \right \rceil $, cancel the $m$-th term $ y^{k P_{l}  - m P_{l - 1} } $ in $Q_k$. Note that all terms in $S_m$ are cancelled in this way, since 
$$ \mu_l - 1 \geq b_{l + 1} m - \left \lceil m \frac{\mu_{l + 2}}{\mu_{l + 1}} \right \rceil, $$
for all $ m \in \{ 1, \ldots, \mu_{l+1} - 1 \} $, by Lemma \ref{Lemma Key2}. On the other hand, if we now fix $ k \in \{ 0, \ldots, \mu_l - 1 \} $, we see that $S_m$ will annihilate the $m$-th term in $Q_k$ precisely when
$$ S_m(1) = b_{l + 1} m - \left \lceil m \frac{\mu_{l + 2}}{\mu_{l + 1}} \right \rceil + 1 \geq k + 1. $$

Let now $ m_k $ be the smallest positive integer such that
$$ b_{l + 1} m_k - \left \lceil m_k \frac{\mu_{l + 2}}{\mu_{l + 1}} \right \rceil \geq k. $$
From Lemma \ref{Lemma 11.4}, we know that $ m_k = \left \lceil k \frac{\mu_{l + 1}}{\mu_{l}} \right \rceil $. Then it follows that the monomials $ y^{k P_{l}  - m P_{l - 1} } $, with $ m_k \leq m \leq \mu_{l + 1} - 1 $, are cancelled in $Q_k $, provided that $ m_k \leq \mu_{l + 1} - 1 $. In the extremal case where $ m_k = \left \lceil k \frac{\mu_{l + 1}}{\mu_{l}} \right \rceil = \mu_{l + 1} $ nothing gets cancelled.

Note in particular that all monomials in $ Q_1, \ldots, Q_{b_{l+1} - 1} $ are cancelled, since in all these cases we have $ m_k = 1 $. We also immediately see that all monomials in $Q_0$ are cancelled.

So we put $ Q_k^* := 0 $, for all $ k \in \{0, \ldots, b_{l+1} - 1 \} $, and 
$$ Q_k^* := - \sum_{m=1}^{\left \lceil k \frac{\mu_{l + 1}}{\mu_{l}} \right \rceil - 1} y^{k P_{l}  - m P_{l - 1} }, $$
for $ k \in \{ b_{l+1}, \ldots, \mu_l - 1 \} $. Note that from Lemma \ref{Lemma 11 Key} we have that $ \mu_l - 1 \geq b_{l+1} $, and  note also that $ Q_k = Q_k^* $ for all $k$ such that $ \left \lceil k \frac{\mu_{l + 1}}{\mu_{l}} \right \rceil = \mu_{l + 1} $. It follows that
$$ \Tr_l + \Tr_{l,l+1}^* = - \sum_{k = b_{l + 1}}^{\mu_{l} - 1} \sum_{m=1}^{\left \lceil k \frac{\mu_{l + 1}}{\mu_{l}} \right \rceil - 1} y^{k P_{l}  - m P_{l - 1} } = \Tr_{l}^*. $$

It remains to prove the statement for $ \Tr_{l - 1, l}^* $. We have that
$$ \Tr_{l - 1, l} = \sum_{k = 1}^{\mu_{l} - 1} \sum_{m = 1}^{b_{l} k } y^{k P_{l}  - m P_{l - 1}}. $$
Let us put $ R_k := \sum_{m = 1}^{b_{l} k } y^{k P_{l}  - m P_{l - 1}} $. If $ k \in \{1, \ldots, b_{l+1} - 1 \} $, one computes that
$$ R_k^* := R_k + Q_k^* = R_k - 0 = \sum_{m = 1}^{b_{l} k } y^{k P_{l}  - m P_{l - 1}} $$
$$ = \sum_{m = \left \lceil k \frac{\mu_{l + 1}}{\mu_{l}} \right \rceil}^{b_{l} k } y^{k P_{l}  - m P_{l - 1}}, $$
since $ \left \lceil k \frac{\mu_{l + 1}}{\mu_{l}} \right \rceil = 1 $ in these cases. As $ b_l \mu_l = \mu_{l + 1} + \mu_{l - 1} > \mu_{l + 1} $, it follows that $ b_l k > k \frac{\mu_{l + 1}}{\mu_{l}} $ for any $ k \in \{ 1, \ldots, \mu_l - 1 \} $, and therefore, the inequality
$$ b_l k \geq \left \lceil k \frac{\mu_{l + 1}}{\mu_{l}} \right \rceil $$
holds. Therefore,
$$ R_k^* := R_k + Q_k^* = \sum_{m = 1}^{b_{l} k } y^{k P_{l}  - m P_{l - 1}} - \sum_{m=1}^{\left \lceil k \frac{\mu_{l + 1}}{\mu_{l}} \right \rceil - 1} y^{k P_{l}  - m P_{l - 1} } $$
$$ = \sum_{m = \left \lceil k \frac{\mu_{l + 1}}{\mu_{l}} \right \rceil}^{b_{l} k} y^{k P_{l}  - m P_{l - 1} }. $$

Consequently, we get that
$$ \Tr_{l - 1, l} + \Tr_l^* = \sum_{k = 1}^{\mu_{l} - 1} R_k^* = \sum_{k = 1}^{\mu_{l} - 1} \sum_{m = \left \lceil k \frac{\mu_{l + 1}}{\mu_{l}} \right \rceil}^{b_{l} k} y^{k P_{l}  - m P_{l - 1} } = \Tr_{l - 1, l}^*, $$
and the proof is complete.
\end{pf}

\vspace{0.5cm}

By induction, we get the following result:

\begin{cor}\label{Cor Key}
$$ \Tr_{0,1}^* = \Tr_{0,1} + \Tr_1 + \ldots + \Tr_{\lambda - 1, \lambda}^*. $$
\end{cor}

We can compute the polynomial $ \Tr_{0,1}^* $ explicitly: 

\begin{prop}\label{Computation}
We have that
$$ \Tr_{0,1}^* = \sum_{r=0}^{\mu_0 - 1} c_r y^r, $$
where $ c_0 = \mu_1 - 1 $, and $ c_r = \mu_1 - \left \lceil r \frac{\mu_{1}}{\mu_{0}} \right \rceil $ for all $ r \in \{ 1, \ldots, \mu_0 - 1\} $.
\end{prop}
\begin{pf}
We have that
$$ \Tr_{0,1}^* = \sum_{k = 1}^{\mu_{1} - 1} \sum_{m = \left \lceil k \frac{\mu_{2}}{\mu_{1}} \right \rceil}^{b_{1} k} y^{k P_{1}  - m P_{0} }. $$
Since $ P_1 = b_1 $ and $ P_0 = 1 $, we have 
$$ \sum_{m = \left \lceil k \frac{\mu_{2}}{\mu_{1}} \right \rceil}^{b_{1} k} y^{k P_{1}  - m P_{0} } = y^{b_{1} k - \left \lceil k \frac{\mu_{2}}{\mu_{1}} \right \rceil} + y^{b_{1} k - \left \lceil k \frac{\mu_{2}}{\mu_{1}} \right \rceil - 1} + \ldots + 1. $$
Let us denote this expression by $F_k$. From Lemma \ref{Lemma Key2}, we know that if $ l,k$ are integers such that $ 1 \leq l < k \leq \mu_1 - 1 $, then
$$ b_1 k - \left \lceil k \frac{\mu_{2}}{\mu_{1}} \right \rceil > b_1 l - \left \lceil l \frac{\mu_{2}}{\mu_{1}} \right \rceil. $$
Therefore $ F_k(1) > F_l(1) $. That is, the number of monomials in $F_k$ strictly increases with $k$. This fact makes it easy to calculate the coefficients of the polynomial $ \sum_{k = 1}^{\mu_{1} - 1} F_k $. Let $ r \in \{ 1, \ldots, \mu_0-1 \} $. As the monomial $y^r$ appears at most once in $ F_k $, we see that its coefficient $c_r$ will equal the number of $F_k$ in which $y^r$ appears. And $ y^r $ will appear in $F_k$ exactly when $ F_k(1) \geq r + 1 $. Let $K_r$ denote the smallest positive integer such that the inequality
$$ b_1 K_r - \left \lceil K_r \frac{\mu_{2}}{\mu_{1}} \right \rceil \geq r $$ 
holds. Then $y^r$ does not appear in any of the polynomials $ F_1, \ldots, F_{K_r - 1} $, but it does appear in all the polynomials $ F_{K_r}, \ldots, F_{\mu_1 - 1} $. Hence $ c_r = \mu_1 - K_r $. By Lemma \ref{Lemma 11.4}, it then follows that $ c_r = \mu_1 - \left \lceil r \frac{\mu_{1}}{\mu_{0}} \right \rceil $. Finally, it is easy to see that $c_0 = \mu_1 - 1 $, which concludes the proof.
\end{pf}

\subsection{Explicit trace formula}\label{explicit formula}
We are now ready to prove the formula for $ \Tr_{\sigma} $. The idea of the proof is as follows: We consider the resolution chain for $ \sigma = (m_1,m_2,n) $. Since we assume that $ n \gg 0 $, we have that the multiplicities $ \mu_l $ descend strictly from $  l = 0 $ to some $ l_0 $. After that, the multiplicities are constant equal to $ m = \Gcd(m_1,m_2) $, and will then strictly increase from $ l = L + 1 - l_1 $ up to $ l = L + 1 $. 

The idea is to ``cut'' the resolution chain somewhere in the constant locus, and treat the two halves independently. For each of the two parts of the chain, we have trace expressions $ \Tr_{l-1,l} $ and $ \Tr_l $ that we can sum up as in Section \ref{star}. The rest of the proof consists of computing the correction term.

Let us remark that in order for the results in Section \ref{star} to apply to both parts of the chain, we will need to perform a certain coordinate change, corresponding to switching the order of the formal branches of the singularity. This process is explained in Section \ref{coordinate change}, and the reader might want to consult this section while reading the proof of Theorem \ref{Formula}.   

\begin{thm}\label{Formula} 
Let $ \sigma = (m_1,m_2,n) $ be a singularity as in Definition \ref{Def. 8.1}, where $ n \gg 0 $. Let $ m = \Gcd(m_1,m_2) $, and let $ \alpha $ resp.~$ \alpha_1 $, resp.~$ \alpha_2 $ be inverses to $m$ resp.~$m_1$, resp.~$m_2$. For any root of unity $ \xi \in \boldsymbol{\mu}_n $, we let $ y = \xi^{\alpha_2} $, $ z = \xi^{\alpha_1} $ and $ w = \xi^{\alpha_1 (m_1/m)} = \xi^{\alpha_2 (m_2/m)} = \xi^{\alpha } $.

Then we have that
$$ \Tr_{\sigma}(\xi) = \sum_{r=0}^{\mu_0 - 1} (\mu_1 - \left \lceil r \frac{\mu_{1}}{\mu_{0}} \right \rceil) y^r + 
\sum_{r=0}^{\mu_{L+1} - 1} (\mu_L - \left \lceil r \frac{\mu_{L}}{\mu_{L+1}} \right \rceil) z^r - 
\sum_{r=0}^{m-1} w^r. $$
The coefficients in this expression depend only on the residue class of $n$ modulo $ \Lcm(m_1,m_2) $.
\end{thm}
\begin{pf}
It suffices to give the proof for a \emph{primitive} root of unity (see the discussion in the proof of Corollary \ref{nonprimformula}). So throughout the proof, $ \xi $ will denote a primitive $n$-th root of unity.

Let $ q \in \{1, \ldots, L \} $ be an index such that $ \mu_{q-2} = \mu_{q-1} = \mu_{q} = \mu_{q+1} = \mu_{q+2} = m $. Since $ \xi $ is primitive, we have that
$$ \Tr_{\sigma} = (\Tr_{y_0} + \ldots + \Tr_{y_q}) + (\Tr_{y_{q+1}} + \ldots + \Tr_{y_L}). $$
Recall from Corollary  \ref{cor. 10.6} that $ \Tr_{y_l} = \Tr_{\mu_l, \mu_{l+1}} + \Tr_{y_l}^{\mu_l} + \Tr_{y_l}^{\mu_{l+1}} $, for $ 1 \leq l \leq L-1 $, and from Lemma \ref{lemma blabla} that $ \Tr_{y_l}^{\mu_{l+1}} + \Tr_{y_{l+1}}^{\mu_{l+1}} = \Tr_{\mu_{l+1}} $. 

We will begin with assuming that $ m = \Gcd(m_1,m_2) \geq 3 $. Then we have that $ \Tr_{\mu_{l}, \mu_{l+1}} + \Tr_{\mu_{l+1}} = \Tr_{l,l+1} + \Tr_{l+1} $ for all $ l \leq q $, and we get the equality
$$ \Tr_{y_0} + \ldots + \Tr_{y_q} = \Tr_{0,1} + \Tr_1 + \ldots + \Tr_q + \Tr_{\mu_q, \mu_{q+1}} + \Tr_{y_q}^{\mu_{q+1}}. $$

We will now use the calculations and the notation from Section \ref{coordinate change}. Using Lemma \ref{lemma 10.25}, one computes easily that $ \Tr_{y_q}^{\mu_{q+1}} = \Tr_{y'_{L-q}}^{\mu'_{L-q}} $. Furthermore, from Corollary \ref{cor. 10.26}, it follows that
$$ \Tr_{q+1} + \ldots + \Tr_{y_L} = \Tr_{y'_0} + \ldots + \Tr_{y'_{L-1-q}}. $$
Since $ \Tr_{y'_{L-1-q}} = \Tr_{\mu'_{L-1-q}, \mu'_{L-q}} + \Tr_{y'_{L-1-q}}^{\mu'_{L-1-q}} + \Tr_{y'_{L-1-q}}^{\mu'_{L-q}} $, we can write
$$ \Tr_{y'_0} + \ldots + \Tr_{y'_{L-1-q}} = \Tr_{0,1}' + \ldots + \Tr_{L-1-q}' + \Tr_{\mu'_{L-1-q}, \mu'_{L-q}} + \Tr_{y'_{L-1-q}}^{\mu'_{L-q}}. $$
So all in all, we get that
$$ \Tr = \Tr_{0,1} + \Tr_1 + \ldots + \Tr_q + \Tr_{\mu_q, \mu_{q+1}} + \Tr_{y'_{L-q}}^{\mu'_{L-q}} + $$
$$ \Tr_{0,1}' + \ldots + \Tr_{L-1-q}' + \Tr_{\mu'_{L-1-q}, \mu'_{L-q}} + \Tr_{y'_{L-1-q}}^{\mu'_{L-q}} = $$
$$ \Tr_{0,1} + \ldots + \Tr_q + \Tr_{\mu_q, \mu_{q+1}} + \Tr_{0,1}' + \ldots + \Tr_{L-q-1}' + \Tr_{L-1-q,L-q}' + \Tr_{L-q}'. $$
Corollary \ref{Cor Key} gives that 
$$ \Tr_{0,1} + \Tr_1 + \ldots + \Tr_q^* = \Tr_{0,1}^*, $$
and that
$$ \Tr_{0,1}' + \ldots + \Tr_{L-q-1}' + \Tr_{L-1-q,L-q}' + (\Tr_{L-q}')^* = (\Tr_{0,1}')^*, $$
where the notation is the same as in Section \ref{star}. Moreover, the terms $ \Tr_{0,1}^* $ and $ (\Tr_{0,1}')^* $ can be explicitly computed using Proposition \ref{Computation}, and these two sums are indeed the first two terms in the theorem.



It remains to compute the correction term $ \Tr_q^0 + \Tr_{\mu_q, \mu_{q+1}} + (\Tr_{L-q}')^0 $, where $ \Tr_q^0 := \Tr_q - \Tr_q^* $, and $ (\Tr_{L-q}')^0 := \Tr_{L-q}' - (\Tr_{L-q}')^* $. It follows from the way we chose $q$, that $ \mu_q = \mu_{q+1} $ and that $ b_{q+1} = 2 $. With our explicit description of $ \Tr_q $ and $ \Tr_q^* $, it is then easy to compute that
$$ \Tr_q^0 = - \sum_{k=0}^{\mu_q - 1} \sum_{n=k}^{\mu_{q+1} - 1} y^{k P_q - n P_{q-1}} + 1. $$
Likewise, we can write
$$ \Tr_{\mu_q, \mu_{q+1}} = \sum_{k=0}^{\mu_q - 1} \sum_{n=0}^{\mu_{q+1} - 1} y^{k P_q - n P_{q-1}}. $$
An easy calculation now shows that
$$ \Tr_q^0 + \Tr_{\mu_q, \mu_{q+1}} = 1 + \sum_{k=1}^{\mu_q - 1} \sum_{n=0}^{k - 1} y^{k P_q - n P_{q-1}}. $$

After ``changing coordinates'' as in Section \ref{coordinate change} below, we get that
$$ (\Tr_{L-q}')^0 = - \sum_{s=0}^{\mu_{q+1} - 1} \sum_{t=s}^{\mu_q - 1} y^{t P_q - s P_{q-1}} + 1. $$
We would now like to calculate 
$$ \sum_{k=1}^{\mu_q - 1} \sum_{n=0}^{k - 1} y^{k P_q - n P_{q-1}} - \sum_{s=0}^{\mu_{q+1} - 1} \sum_{t=s}^{\mu_q - 1} y^{t P_q - s P_{q-1}}. $$
Before we do that, note that $ \mu_q = \mu_{q+1} = m $, and that the sum above may be written
$$ \sum_{k=1}^{m - 1} \sum_{l=0}^{k - 1} y^{k P_q - l P_{q-1}} - \sum_{l=0}^{m - 1} \sum_{k=l}^{m - 1} y^{k P_q - l P_{q-1}}. $$ 
However, it is easily seen that
$$ \sum_{k=1}^{m - 1} \sum_{l=0}^{k - 1} y^{k P_q - l P_{q-1}} = \sum_{l=0}^{m - 2} \sum_{k=l+1}^{m - 1} y^{k P_q - l P_{q-1}}, $$
and that
$$ \sum_{l=0}^{m - 2} \sum_{k=l+1}^{m - 1} y^{k P_q - l P_{q-1}} - \sum_{l=0}^{m - 1} \sum_{k=l}^{m - 1} y^{k P_q - l P_{q-1}} $$
$$ = - \sum_{l=0}^{m - 1} y^{l P_q - l P_{q-1}} = - \sum_{l=0}^{m - 1} (y^{P_q - P_{q-1}})^l, $$
so it follows that
$$ \Tr_q^0 + \Tr_{\mu_q, \mu_{q+1}} + (\Tr_{L-q}')^0 = 2 - \sum_{l=0}^{m - 1} (y^{P_q - P_{q-1}})^l. $$

Note now that $ y^{P_q - P_{q-1}} = \chi^{-r_q + r_{q+1}} $. It is easily seen that
$$ - \mu_l r_l + \mu_{l+1} r_{l-1} = - \mu_{l+1} r_{l+1} + \mu_{l+2} r_{l} $$
for all $l$, so by induction, we get that 
$$ - \mu_q r_q + \mu_{q+1} r_{q-1} = \mu_{L+1} r_{L-1}, $$ 
since $ r_L = 0 $. As $ \mu_q = \mu_{q+1} = m $, we get, after dividing by $m$, that
$$ -r_q + r_{q+1} = \mu_{L+1}/m = m_1/m, $$
remembering that $ r_{L-1} = 1 $. So it follows that 
$$ y^{P_q - P_{q-1}} = \chi^{\mu_{L+1}/m} = \xi^{\alpha_1 (m_1/m)}, $$ 
and the proof of the formula is complete in the case $m \geq 3$.

The two cases $ m = 1 $ and $ m = 2 $ remain. The proof in these cases is very similar to the one above, and is therefore omitted here. The main difference is that we have to cut the chain in three pieces: The part where the multiplicities descend, the part where the multiplicities are constant, and the part where the multiplicities ascend. On each of these parts, we have to compute certain sums of trace expressions, which can be done along the same lines as in the case where $ m \geq 3 $.

\end{pf}

\vspace{0.3cm}

We would like to end this section with the remark that now that we have obtained the explicit formula in Theorem \ref{Formula}, Theorem \ref{thm. 9.13} gives an effective formula for computing the trace $ \Tr_{\beta}(e(H^{\bullet}(g|_{\mathcal{Y}_k}))) $, where the notation is the same as in Section \ref{special filtrations}. We formulate this in Theorem \ref{improved formula} below. Recall the standard assumptions in Sections \ref{assumption on degree} and \ref{assumption on surface}.

\begin{thm}\label{improved formula}
Let $ g \in G $ be a group element corresponding to a root of unity $\xi \in \boldsymbol{\mu}_n$, where $ n \gg 0 $. Then we have that
$$ \Tr_{\beta}(e(H^{\bullet}(g|_{\mathcal{Y}_k}))) = \sum_{\upsilon \in \mathcal{V}} \Tr_{\upsilon}(\xi) + \sum_{\varepsilon \in \mathcal{E}} \Tr_{\sigma(\varepsilon)}(\xi), $$ 
where $ \sigma(\varepsilon) $ is the unique singularity of $ \mathcal{X}' $ corresponding to $ \varepsilon $. The contributions $ \Tr_{\upsilon}(\xi) $ are given by Proposition \ref{prop. 9.7}, and the contributions $ \Tr_{\sigma(\varepsilon)}(\xi) $ are given by Theorem \ref{Formula}.

Furthermore, this expression depends only on $ \Gamma(\mathcal{X}_k) $ and the functions $ \mathfrak{g} $ and $ \mathfrak{m} $ defined in Section \ref{9.1}. 
 \end{thm}

\subsection{Coordinate change}\label{coordinate change}
We consider now positive integers $ m_1 $, $m_2$ such that $ \Gcd(m_1,m_2) = 1 $. Let $ M = \Lcm(m_1,m_2) $, and let $ n $ be a positive integer such that $ \Gcd(n,M) = 1 $. Let $ \alpha_i $ denote the inverse to $ m_i $ modulo $n$. We denote by $ \mu_l $, $ b_l $ and $ r_l $ the numerical data associated to the singularity $ \sigma = (m_1, m_2, n) $ as usual.

Let us now define $ m_1' := m_2 $ and $ m_2' := m_1 $, and consider the singularity $ \sigma' = (m_1', m_2', n) $, which is the same singularity as $(m_1,m_2,n)$, but with reverse ordering of the branches. Let $ \mu'_j = \mu_{L+1-j} $, $c_j = b_{L+1-j} $. Then the numerical data for $ (m_1', m_2', n) $ consists of $ \mu'_j $, $ c_j $ and $ s_j $, where the $c_j$ and $s_j$ satisfy the equations $ s_{j-1} = c_{j+1} s_j - s_{j+1} $. 

The intersection points of the components in the exceptional locus are $ y_l = C_l \cap C_{l+1} $, so we let $ y'_j = y_{L-j} $, where $ l $ and $ j $ run from $ 0 $ to $L$. 

Let $ \alpha'_1 $ denote the inverse to $m_1' $ modulo $n$. Notice then that the equation $ m_1 + r_0 m_2 = n \mu_1 $ gives that $ \alpha'_1 \equiv_n - \alpha_1 r_0 $, and $ m_1' + s_0 m_2' = n \mu_1' $ gives $ \alpha_1 \equiv_n - \alpha'_1 s_0 $. 

\begin{lemma}\label{lemma 10.25}
Let $ \chi = \xi^{\alpha_1} $ and $ \chi' = \xi^{\alpha_1'} $. Then we have that $ \chi^{r_{L-1-j}} = \chi'^{-s_j} $ for all $ 0 \leq j \leq L-1 $.
\end{lemma}
\begin{pf}
As $ r_{L-1} = 1 $, we get that $ \chi^{r_{L-1}} = \xi^{\alpha_1} = \xi^{- \alpha'_1 s_0} = \chi'^{ - s_0 } $. Assume now that $ j > 0 $. By induction it then follows that 
$$ \chi^{r_{L-1-j}} = \xi^{\alpha_1 r_{L-1-j}} = \xi^{\alpha_1 b_{L+1-j} r_{L-j}} \xi^{ - \alpha_1 r_{L + 1 -j}} = \xi^{ - \alpha'_1 s_{j-1} c_j} \xi^{ \alpha'_1 s_{j-2}}$$
$$ = \xi^{ - \alpha'_1(c_j s_{j-1} - s_{j-2})} = \xi^{ - \alpha'_1 s_j } = \chi'^{-s_j}. $$
\end{pf}

Let us now assume that $ \xi \in \boldsymbol{\mu}_n $ is a primitive root. Then we can write 
$$ \Tr_{\sigma}(\xi) = \Tr_{y_0} + \ldots + \Tr_{y_l} + \ldots + \Tr_{y_L}, $$
where $ \Tr_{y_l} = \frac{1 - \xi}{(1 - \chi^{r_{l-1}})(1 - \chi^{- r_{l}})} $, $ \Tr_{y_0} = \frac{\mu_1}{1 - \chi^{- r_0}} $ and $ \Tr_{y_L} = \frac{\mu_L}{1 - \chi^{r_{L-1}}} $.

We similarly have that 
$$ \Tr_{\sigma'}(\xi) = \Tr_{y'_0} + \ldots + \Tr_{y'_j} + \ldots + \Tr_{y'_L}, $$
where $ \Tr_{y'_j} = \frac{1 - \xi}{(1 - \chi'^{s_{j-1}})(1 - \chi'^{- s_{j}})} $, $ \Tr_{y'_0} = \frac{\mu_1'}{1 - \chi'^{- s_0}} $ and $ \Tr_{y'_L} = \frac{\mu'_L}{1 - \chi'^{s_{L-1}}} $.

Using Lemma \ref{lemma 10.25}, we can draw the following conclusion:
\begin{cor}\label{cor. 10.26}
We have that $ \Tr_{y'_j} = \Tr_{y_{L-j}} $ for all $ j $, where $ 0 \leq j \leq L $. 
\end{cor}

\vspace{0.5cm}

Let us also remark that we will use the notation $ Tr_{\mu_{j-1}',\mu_j'} $, $ Tr_{\mu_j'} $ etc., for the expressions defined as in Section \ref{trace formula}, and the notation $ Tr'_{j-1,j} $ and $ Tr'_j $ for the expressions defined as in Section \ref{lemmas}.

\section{Character computations and jumps}\label{computations and jumps}
Let $X/K$ be a smooth, projective and geometrically irreducible curve such that $X(K) \neq \emptyset $, and let $ \mathcal{X}/S $ be the minimal SNC-model of $X$. We have in previous sections studied properties of the action of $ \boldsymbol{\mu}_n $ on the cohomology groups $ H^i(\mathcal{Y}_k, \mathcal{O}_{\mathcal{Y}_k}) $, where $ \mathcal{Y} $ is the minimal desingularization of the pullback $ \mathcal{X}_{S'} $ for some tame extension $S'/S$ of degree $n$.   

Let $\mathcal{J}/S$ be the N\'eron model of the Jacobian of $X$. We will in this section apply our results to the study of the filtration $ \{ \mathcal{F}^a \mathcal{J}_k \} $, where $ a \in \mathbb{Z}_{(p)} \cap [0,1] $, that we defined in Section \ref{ratfil}. We will first prove some general properties for these filtrations, and then present some computations for curves of genus $g = 1$ and $g = 2$.

We would at this point like to remark that in order to make the $ \boldsymbol{\mu}_n $-action on $ H^1(\mathcal{Y}_k, \mathcal{O}_{\mathcal{Y}_k}) $ compatible with the action on $ T_{\mathcal{J}'_k, 0} $, we have to let $ \boldsymbol{\mu}_n $ act on $ R' $ by $ [\xi](\pi') = \xi^{-1} \pi' $, for any $ \xi \in \boldsymbol{\mu}_n $. We made the choice in previous sections, when working with surfaces, to let $ \boldsymbol{\mu}_n $ act by $ [\xi](\pi') = \xi \pi' $, in order to get simpler notation. This means that the irreducible characters for the representation on $ T_{\mathcal{J}'_k, 0} $ are the \emph{inverse} characters to those we compute when using our formulas for the representation on $ H^1(\mathcal{Y}_k, \mathcal{O}_{\mathcal{Y}_k}) $.  

\subsection{Filtrations for N\'eron models of Jacobians}
Theorem \ref{thm. 9.13} states that the Brauer trace of the automorphism induced by any group element $ \xi \in \boldsymbol{\mu}_n $ on the formal difference $ H^0(\mathcal{Y}_k, \mathcal{O}_{\mathcal{Y}_k}) - H^1(\mathcal{Y}_k, \mathcal{O}_{\mathcal{Y}_k}) $ only depends on the combinatorial structure of $ \mathcal{X}_k $. With the assumption that $X(K) \neq \emptyset$, we can actually improve this result, and get a similar result for the character of the representation of $ \boldsymbol{\mu}_n $ on $ H^1(\mathcal{Y}_k, \mathcal{O}_{\mathcal{Y}_k}) $:

\begin{thm}\label{main character theorem}
Let $ X/K $ be a smooth, projective and geometrically connected curve having genus $ g(X) > 0 $, and assume that $ X(K) \neq \emptyset $. Let $ \mathcal{X} $ be the minimal SNC-model of $ X $ over $S$. Furthermore, let $ S'/S $ be a tame extension of degree $n$, where $ n $ is relatively prime to the least common multiple of the multiplicities of the irreducible components of $ \mathcal{X}_k $, and let $ \mathcal{Y}/S' $ be the minimal desingularization of $ \mathcal{X}_{S'} $. 

Then the irreducible characters for the representation of $ \boldsymbol{\mu}_n $ on $ H^1(\mathcal{Y}_k, \mathcal{O}_{\mathcal{Y}_k}) $ only depend on the intersection graph $ \Gamma(\mathcal{X}_k) $, together with the functions $ \mathfrak{g} $ and $ \mathfrak{m} $.  
\end{thm}
\begin{pf}
For any $ g \in G $, corresponding to a root $ \xi \in \boldsymbol{\mu}_n $, we have that 
$$ \Tr_{\beta}(e(H^{\bullet}(g|_{\mathcal{Y}_k}))) = \sum_{\upsilon \in \mathcal{V}} \Tr_{\upsilon}(\xi) + \sum_{\varepsilon \in \mathcal{E}} \Tr_{\varepsilon}(\xi), $$
by Theorem \ref{thm. 9.13}. The contributions $ \Tr_{\upsilon}(\xi) $ can be computed using Proposition \ref{prop. 9.7}, and for the contributions $ \Tr_{\varepsilon}(\xi) $, we use Proposition \ref{lemma 9.8} and Proposition \ref{lemma 9.10}. In this way, we obtain a formula for the Brauer trace of the automorphism induced by any $ \xi \in \boldsymbol{\mu}_n $ on the formal difference
$$ H^0(\mathcal{Y}_k, \mathcal{O}_{\mathcal{Y}_k}) - H^1(\mathcal{Y}_k, \mathcal{O}_{\mathcal{Y}_k}). $$ 

Since $X$, and hence $X_{K'}$, has a rational point, it follows from \cite{Liubook}, Corollary 9.1.32 that at least one of the irreducible components of $ \mathcal{Y}_k $ has multiplicity $1$. We can therefore conclude that $ H^0(\mathcal{Y}_k, \mathcal{O}_{\mathcal{Y}_k}) = k $ (\cite{Arwin}, Lemma 2.6). Furthermore, the $ \boldsymbol{\mu}_n $-action on $ \mathcal{Y}_k $ is relative to the ground field $k$, so it follows that the eigenvalue for the automorphism induced by $ \xi $ on $ H^0(\mathcal{Y}_k, \mathcal{O}_{\mathcal{Y}_k}) $ equals $1$. 

We therefore obtain the formula
$$ \Tr_{\beta}(H^1(g|_{\mathcal{Y}_k})) = 1 - (\sum_{\upsilon \in \mathcal{V}} \Tr_{\upsilon}(\xi) + \sum_{\varepsilon \in \mathcal{E}} \Tr_{\varepsilon}(\xi)). $$
Since the expressions $ \Tr_{\upsilon}(\xi) $ and $ \Tr_{\varepsilon}(\xi) $ only depend on the combinatorial structure of $ \mathcal{X}_k $, the same is true for $ \Tr_{\beta}(H^1(g|_{\mathcal{Y}_k})) $. This completes the proof, since the Brauer character for the representation of $ \boldsymbol{\mu}_n $ on $ H^1(\mathcal{Y}_k, \mathcal{O}_{\mathcal{Y}_k}) $ is determined by the Brauer trace for the group elements $ \xi \in \boldsymbol{\mu}_n $.
\end{pf}

\vspace{0.5cm}

Let $ \mathcal{J}/S $ be the N\'eron model of the Jacobian of $X/K$, and let $ \{ \mathcal{F}^a \mathcal{J}_k \} $, where $ a \in \mathbb{Z}_{(p)} \cap [0,1] $, be the filtration of $ \mathcal{J}_k $ defined in Section \ref{ratfil}. Then Theorem \ref{main character theorem} has the following consequence:

\begin{cor}\label{main jump corollary}
The jumps in the filtration $ \{ \mathcal{F}^a \mathcal{J}_k \} $ with indices in $ \mathbb{Z}_{(p)} \cap [0,1] $ depend only on the intersection graph $ \Gamma(\mathcal{X}_k) $, together with the functions $ \mathfrak{g} $ and $ \mathfrak{m} $. In particular, they don't depend on $p$.
\end{cor}
\begin{pf}
Let $S'/S$ be a tame extension of degree $n$, where $n$ is prime to $l$, the least common multiple of the multiplicities of the irreducible components of $ \mathcal{X}_k $. Let $ \mathcal{J}'/S' $ be the N\'eron model of the Jacobian of $X_{K'}$. Recall from Section \ref{jacobiancase} that we could make the identification $ H^1(\mathcal{Y}_k, \mathcal{O}_{\mathcal{Y}_k}) \cong T_{\mathcal{J}'_k,0} $. 

The jumps in the filtration of $ \mathcal{J}_k $ induced by the extension $S'/S$ are determined by the irreducible characters for the representation of $ \boldsymbol{\mu}_n $ on $ T_{\mathcal{J}'_k,0} $. However, this representation is precisely the representation of $ \boldsymbol{\mu}_n $ on $ H^1(\mathcal{Y}_k, \mathcal{O}_{\mathcal{Y}_k}) $, if we let $ \boldsymbol{\mu}_n $ act on $R'$ by $ [\xi](\pi') = \xi^{-1} \pi' $, for every $\xi$. By Theorem \ref{main character theorem}, the character for this representation only depends on $ \Gamma(\mathcal{X}_k) $, $ \mathfrak{g} $ and $ \mathfrak{m} $.  

Since $ \mathbb{Z}_{(lp)} \cap [0,1] $ is \emph{dense} in $ \mathbb{Z}_{(p)} \cap [0,1] $, we conclude that the jumps of the filtration $ \{ \mathcal{F}^a \mathcal{J}_k \} $ with indices in $ \mathbb{Z}_{(p)} \cap [0,1] $ only depend on $ \Gamma(\mathcal{X}_k) $, $ \mathfrak{g} $ and $ \mathfrak{m} $. 
\end{pf}

\vspace{0.5cm}

With the two results above at hand, we can draw some conclusions about \emph{where} the jumps occur in the case of Jacobians. Let us first recall the following terminology from \cite{Tame}: An irreducible component $C$ of $ \mathcal{X}_k $ is called \emph{principal} if either $ P_a(C) > 0 $, or if $ C $ is smooth and rational and meets the rest of the components of $ \mathcal{X}_k $ in at least three points.  

\begin{cor}\label{specific jump corollary}
Let $ \tilde{n} $ be the least common multiple of the multiplicities of the principal components of $ \mathcal{X}_k $. Then the jumps in the filtration $ \{ \mathcal{F}^a \mathcal{J}_k \} $ occur at indices of the form $ i/\tilde{n} $, where $ 0 \leq i < \tilde{n} $. 
\end{cor}
\begin{pf}
Let us first note that if $X$ obtains semi-stable reduction over a tame extension $K'/K$, then the Jacobian of $X$ obtains semi-abelian reduction over the same extension. Furthermore, the minimal extension that gives semi-abelian reduction is the tame extension $ \widetilde{K}/K $ of degree $ \tilde{n} $ (\cite{Tame}, Theorem 7.1). So in this case, the statement follows from Proposition \ref{tamejumpprop}.

Let us now assume that $X$ needs a wildly ramified extension to obtain semi-stable reduction. Consider the combinatorial data $ (\Gamma(\mathcal{X}_k), \mathfrak{g}, \mathfrak{m}) $. It follows from \cite{Winters}, Corollary 4.3, that we can find an SNC-model $ \mathcal{Z}/\Spec(\mathbb{C}[[t]]) $, where the generic fiber of $ \mathcal{Z} $ is smooth, projective and geometrically connected, and where the special fiber of $ \mathcal{Z} $ has the \emph{same} combinatorial data as $ \mathcal{X}_k $. 

Let $ \mathcal{J}_{\mathcal{Z}} $ be the N\'eron model of the Jacobian of the generic fiber of $ \mathcal{Z} $. Then the jumps of the filtration $ \{ \mathcal{F}^a \mathcal{J}_{\mathcal{Z},\mathbb{C}} \} $  occur at indices of the form $ i/\tilde{n} $, where $ 0 \leq i < \tilde{n} $. The result follows now from Corollary \ref{main jump corollary}.
\end{pf}

\subsection{}
Let $X/K$ be a smooth, projective and geometrically connected curve, and assume that $ X(K) \neq \emptyset $. Let $ \mathcal{X}/S $ be the minimal SNC-model of $X/K$. It is known that for a fixed genus $ g \geq 2 $, there are only finitely many possibilities for the combinatorial structure of the special fiber of $ \mathcal{X}/S$, modulo chains of $(-2)$-curves (\cite{Arwin},  Theorem 1.6). The same statement is, as we shall see below, also true for elliptic curves. 

Let $ \mathcal{J}/S $ be the N\'eron model of the Jacobian of $X$. Since, by Corollary \ref{main jump corollary},  the jumps of the filtration $ \{ \mathcal{F}^a \mathcal{J}_k \} $ only depend on the combinatorial structure of $ \mathcal{X}_k $, one can, for each $g > 0$, classify these jumps. In the next sections, we will give the jumps for every fiber type of genus $1$ and $2$. 

\begin{rmk}
It is easy to see that chains of $(-2)$-curves do not affect the jumps.
\end{rmk}

\subsection{Computations of jumps for  $ g =1 $} 
Let $X/K$ be an \emph{elliptic} curve, and let $ \mathcal{E} $ be the minimal regular model of $ X $. It is a well known fact that there are only finitely many possibilities for the combinatorial structure of the special fiber $ \mathcal{E}_k $, modulo chains of $(-2)$-curves. The various possibilities were first classified in \cite{Kod}, and this is commonly referred to as the \emph{Kodaira classification}. For another treatment of this theory, we refer to \cite{Liubook}, Chapter 10.2. If now $ \mathcal{X}/S $ denotes the minimal SNC-model of $X$, it follows that there are only finitely many possibilities for the combinatorial structure of $ \mathcal{X}_k $, each one derived from the Kodaira classification. The symbols $ I, II, \ldots $ appearing in Table \ref{table 1} below are known as the \emph{Kodaira symbols} and refer to the fiber types in the Kodaira classification.

Let $ \mathcal{J}/S $ be the N\'eron model of $ J(X) = X $. It follows from Corollary \ref{main jump corollary} and Corollary \ref{specific jump corollary} that the (unique) jump in the filtration $ \{ \mathcal{F}^a \mathcal{J}_k \} $ only depends on the fiber type of $ \mathcal{X}/S $, and can only occur at finitely many \emph{rational} numbers. In Table \ref{table 1} below, we list the jumps for the various Kodaira types. Note that we obtain the same list as the one computed in \cite{Edix} by R. ~Schoof.

We would like to say a few words about how these computations are done. For each fiber type, we consider an infinite sequence $ (n_j)_{j \in \mathbb{N}} $, depending on the fiber type, where $ n_j \rightarrow \infty $ as $ j \rightarrow \infty $. For each $n_j$ in this sequence, let $ R_j/R $ be the tame extension of degree $n_j$, and let $ \pi_j $ be the uniformizing parameter of $ R_j $. Furthermore, let $ \boldsymbol{\mu}_{n_j} $ act on $ R_j $ by $ [\xi](\pi_j) = \xi \pi_j $. We can then use Theorem \ref{thm. 9.13} to compute the character for the induced representation of $ \boldsymbol{\mu}_{n_j} $ on $ H^1(\mathcal{Y}^j_k, \mathcal{O}_{\mathcal{Y}^j_k}) $, where $ \mathcal{Y}^j $ denotes the minimal desingularization of $ \mathcal{X}_{S_j} $, and where $ S_j = \Spec(R_j) $. This character is on the form $ \chi(\xi) = \xi^{i(j)} $. 

The character for the representation of $ \boldsymbol{\mu}_{n_j} $ on $ T_{\mathcal{J}^j_k, 0} $ is the inverse of this character, $ \chi^{-1}(\xi) = \xi^{- i(j)} $. The jump of $ \{ \mathcal{F}^a \mathcal{J}_k \} $ will then be given by the limit of the expression $ [- i(j)]_{n_j}/n_j $ as $ j \rightarrow \infty $, where $ [- i(j)]_{n_j} \equiv_{n_j} - i(j) $, and $ 0 \leq [- i(j)]_{n_j} < n_j $.

In Example \ref{example genus 1} below, we explain in detail how these computations are done for fiber type $IV$ in the Kodaira classification.

\begin{ex}\label{example genus 1}
Let $ \mathcal{X}/S $ have fibertype $IV$. In this case, the combinatorial data of $ \mathcal{X}_k $ consists of the set of vertices $ \mathcal{V} = \{ \upsilon_1, \ldots, \upsilon_4 \} $, where $ \mathfrak{m}(\upsilon_i) = 1 $ for $ i \in\{ 1,2,3 \}$, and $ \mathfrak{m}(\upsilon_4) = 3 $. Furthermore, we have that $ \mathfrak{g}(\upsilon_i) = 0 $ for all $i$. The set of edges consists of $ \mathcal{E} = \{ \varepsilon_1, \varepsilon_2, \varepsilon_3 \} $, where $ \varepsilon_i $ corresponds to the unique intersection point of the components $ \upsilon_i $ and $ \upsilon_4 $, for $ i = 1,2,3 $. Let us choose the ordering $ (\upsilon_i,\upsilon_4) $ for all $i$. 

Let now $ n \gg 0 $ be a positive integer relatively prime to $ p $ and to $ \Lcm( \{ \mathfrak{m}(\upsilon_i) \} ) = 3 $, and let $R'/R$ be a tame extension of degree $n$. Let $ \boldsymbol{\mu}_n $ act on $R'$ by $ [\xi](\pi') = \xi \pi' $ for any $ \xi \in \boldsymbol{\mu}_n$, where $ \pi' $ is a uniformizing parameter for $R'$. 

For any $ g \in G $, corresponding to a root of unity $ \xi \in \boldsymbol{\mu}_n $, Theorem \ref{thm. 9.13} states that
$$ \Tr_{\beta}(e(H^{\bullet}(g|_{ \mathcal{Y}_k}))) = \sum_{\upsilon \in \mathcal{V}} \Tr_{\upsilon}(\xi) + \sum_{\varepsilon \in \mathcal{E}} \Tr_{\varepsilon}(\xi). $$

Let $ \sigma $ be the singularity $ (1,3,n) $. Then we have that $ \Tr_{\varepsilon_i}(\xi) = \Tr_{\sigma}(\xi) $ for all $ i \in \{ 1,2,3 \} $. It suffices to consider the case where $ n \equiv_3 1 $. One computes easily that $ \mu_l = 1 $ for all $ l \in \{ 1, \ldots, L(\sigma) \} $. From Theorem \ref{Formula}, we immediately get that $ \Tr_{\varepsilon_i}(\xi) = 1 $, for all $i$. 

Proposition \ref{prop. 9.7} states that 
$$ \Tr_{\upsilon}(\xi) = \sum_{k=0}^{m_{\upsilon}-1} (\xi^{\alpha_{m_{\upsilon}}})^{k} ((m_{\upsilon} - k)C_{\upsilon}^2  + 1 - p_a(C_{\upsilon})), $$
for any $ \upsilon \in \mathcal{V} $, where $ \alpha_{m_{\upsilon}} m_{\upsilon} \equiv_n 1 $. As $ C_{\upsilon_i}^2 = - 1 $ for $ i \in \{ 1,2,3 \} $, we see that $ \Tr_{\upsilon_i}(\xi) = 0 $ for these vertices, and since $ C_{\upsilon_4}^2 = - 1 $, it follows that $ \Tr_{\upsilon_4}(\xi) = - 2 - \xi^{\alpha_3} $. In total, we get
$$  \Tr_{\beta}(e(H^{\bullet}(g|_{ \mathcal{Y}_k}))) = 3 + (- 2 - \xi^{\alpha_3}) = 1 - \xi^{\alpha_3}. $$
We can therefore conclude that the character for the representation of $ \boldsymbol{\mu}_n $ on $ H^1(\mathcal{Y}_k, \mathcal{O}_{\mathcal{Y}_k}) $ is $ \chi(\xi) = \xi^{\alpha_3} $.


In order to compute the jump of the filtration $ \{ \mathcal{F}^a \mathcal{J}_k \} $, where $ \mathcal{J} $ is the N\'eron model of $J(X) = X$, we have to use the \emph{inverse} character, which is $ \chi^{-1}(\xi) = \xi^{[- \alpha_3]_n} $, where $ [- \alpha_3]_n = - \alpha_3 $ modulo $n$, and $ 0 \leq [- \alpha_3]_n < n $. The jump will be given by the limit of the expression $ ([- \alpha_3]_n)/n $ as $n$ goes to infinity over integers $n$ that are equivalent to $1$ modulo $3$. 

Since $ n = 1 + 3 \cdot h $, for some integer $h$, we get that $ \alpha_3 = \frac{1 + 2 n}{3} $, where $ 0 < \alpha_3 < n $. Therefore, the jump occurs at the limit of $ ([- \alpha_3]_n)/n = \frac{n - 1}{3 n} $ which is $ 1/3 $.


\end{ex}
 
\begin{table}[htb]\caption{Genus $1$}\label{table 1}
\begin{tabular}{|c|c|c|c|c|c|c|c|c|c|c|} 
\hline 
Fibertype & $(I)$ & $(I)^*$ & $(I_n)$ & $(I_n)^*$ & $(II)$ & $(II)^*$ & $(III)$ & $(III)^*$ & $(IV)$ & $(IV)^*$\\
\hline
Jumps & $0$ & $1/2$ & $0$ & $1/2$ & $ 1/6 $ & $ 5/6 $ & $ 1/4 $ & $ 3/4 $ & $ 1/3 $ & $ 2/3 $ \\\hline
\end{tabular}
\end{table}

\subsection{Computations of jumps for $g=2$}\label{genus 2}
Let $ X/K $ be a curve having genus equal to $2$. Like in the case for elliptic curves, there are finitely many possibilities, modulo chains of $(-2)$-curves, for the combinatorial structure of the special fiber of the minimal regular model of $X$. Moreover, there exists a complete classification of the various possible fiber types. This classification is mainly due to A.P. Ogg (\cite{Ogg}), with the exception of a few missing cases which were filled in by Y. Namikawa and K. ~Ueno in \cite{Ueno}. We will refer to the list of possible fiber types as the \emph{Ogg-classification}, and we will use the indexing from \cite{Ogg}, with the addition of the types $41_a$, $41_b$ and $41_c$ from \cite{Ueno} that were missing in \cite{Ogg}.

Let $ \mathcal{X}/S $ be the minimal SNC-model of $X$. Then there are only finitely many possibilities for the combinatorial structure of $ \mathcal{X}_k $, each derived from the fiber types in the Ogg-classification. Let $ \mathcal{J}/S $ be the N\'eron model of the Jacobian of $X$. The jumps in the filtration $ \{ \mathcal{F}^a \mathcal{J}_k \} $ depend only on the combinatorial structure of $ \mathcal{X}_k $, and can occur only at a finite set of rational numbers. 

In order to compute the jumps for each fibertype, we proceed more or less in the same manner as we did in the case of elliptic curves. In Example \ref{example 13.1}, we explain in detail how this is done for fiber type $4$ in the Ogg-classification. 


Most of the results for genus $2$ curves are gathered in Table \ref{table 2} below. However, some cases are treated separately in Section \ref{special types}.

\begin{ex}\label{example 13.1}
We consider fiber type $4$ in the Ogg-classification. In this case, the set of vertices of $ \Gamma(\mathcal{X}_k) $ is $ \mathcal{V} = \{ \upsilon_1, \ldots, \upsilon_7 \} $, where $ \mathfrak{g}(\upsilon_i) = 0 $ for all $i$. Furthermore, we have that $ \mathfrak{m}(\upsilon_i) = 1 $ for $ i = 1,7 $, $ \mathfrak{m}(\upsilon_i) = 2 $ for $ i = 2,5,6 $, $ \mathfrak{m}(\upsilon_3) = 3 $ and $ \mathfrak{m}(\upsilon_4) = 4 $. The set of edges is $ \mathcal{E} = \{ \varepsilon_1, \varepsilon_2, \varepsilon_3, \varepsilon_4, \varepsilon_5, \varepsilon_6 \} $, where $ \varepsilon_1 = (\upsilon_1,\upsilon_2) $, $ \varepsilon_2 = (\upsilon_2,\upsilon_3) $, $ \varepsilon_3 = (\upsilon_3,\upsilon_4) $, $ \varepsilon_4 = (\upsilon_5,\upsilon_4) $, $ \varepsilon_5 = (\upsilon_6,\upsilon_4) $  and $ \varepsilon_6 = (\upsilon_7,\upsilon_4) $.

We have that $ \Lcm( \{\mathfrak{m}(\upsilon_i)\}) = 12 $. Let $ n \gg 0$ be any integer not divisible by $p$, and such that $ n \equiv_{12} 1 $. Let $ R'/R $ be the extension of degree $n$, and let $ \pi'$ be the uniformizing parameter of $R'$. Let $ \mathcal{Y} $ be the minimal desingularization of $ \mathcal{X}_{S'} $. We let $ \boldsymbol{\mu}_n $ act on $ R'$ by $ [\xi](\pi') = \xi \pi' $, for any $ \xi \in \boldsymbol{\mu}_n $.

Now, let $ \xi \in \boldsymbol{\mu}_n$ be a root of unity. For any $ \upsilon \in \mathcal{V} $, Proposition \ref{prop. 9.7} gives that
$$ \Tr_{\upsilon}(\xi) = \sum_{k=0}^{m_{\upsilon}-1} (\xi^{\alpha_{m_{\upsilon}}})^{k} ((m_{\upsilon} - k)C_{\upsilon}^2  + 1 - p_a(C_{\upsilon})). $$
As the computations are similar for all $ \upsilon \in \mathcal{V} $, we only do this explicitly for $ \upsilon_3 $. We have that $ p_a(C_{\upsilon_3}) = \mathfrak{g}(\upsilon_3) = 0 $, so it remains only to compute $ C_{\upsilon_3}^2 $. The edge $ \varepsilon_2 $ corresponds to the singularity $ \sigma_2 = (2,3,n) $ and the edge $ \varepsilon_3 $ corresponds to the singularity $ \sigma_3 = (3,4,n) $. Denote by $ C_l^{\sigma_2} $ the exceptional components in the resolution of $ \sigma_2 $, and by $ C_l^{\sigma_3} $ the components in the resolution of $ \sigma_3 $. Then $ C_1^{\sigma_2} $ and $ C_L^{\sigma_3} $ are the only two components of $ \mathcal{Y}_k $ that meet $ C_{\upsilon_3} $ (note the ordering of the formal branches in $ \sigma_2 $ and $ \sigma_3 $). It is easily computed that $ \mu_1^{\sigma_2} = 2 $ and that $ \mu_L^{\sigma_3} = 1 $. So it follows that $ C_{\upsilon_3}^2 = - 1 $, and therefore
$$ \Tr_{ \upsilon_3 }(\xi) = - 2 - \xi^{\alpha_3}. $$
For the other vertices, we compute that
$$ \Tr_{ \upsilon_1 }(\xi) = \Tr_{ \upsilon_7 }(\xi) = 0, $$
$$ \Tr_{ \upsilon_2 }(\xi) = \Tr_{ \upsilon_5 }(\xi) = \Tr_{ \upsilon_6 }(\xi) = - 1, $$
and
$$ \Tr_{ \upsilon_4 }(\xi) = - 7 - 5 \xi^{\alpha_4} - 3 (\xi^{\alpha_4})^2 - (\xi^{\alpha_4})^3. $$

Next, we must compute the contributions from the singularities. We will only write out the details for $ \varepsilon_3 = (\upsilon_3,\upsilon_4) $. In this case, we need to compute $ \Tr_{\sigma_3}(\xi) $. It is easily computed that $ \mu_1^{\sigma_3} = 3 $ and $ \mu_L^{\sigma_3} = 1 $. Theorem \ref{Formula} then gives that 
$$ \Tr_{\varepsilon_3}(\xi) = \Tr_{\sigma_3}(\xi) = 3 + 2 \xi^{\alpha_4} + (\xi^{\alpha_4})^2. $$
For the contributions from the other edges, we compute in a similar fashion that
$$ \Tr_{\varepsilon_1}(\xi) = \Tr_{\varepsilon_6}(\xi) = 1, $$
$$ \Tr_{\varepsilon_2}(\xi) = 2 + \xi^{\alpha_3}, $$
and
$$ \Tr_{\varepsilon_4}(\xi) = \Tr_{\varepsilon_5}(\xi) = 3 + \xi^{\alpha_4} + (\xi^{\alpha_4})^2. $$

Summing up, we get 
$$ \sum_{i=1}^7 \Tr_{\upsilon_i }(\xi) + \sum_{i=1}^6 \Tr_{\varepsilon_i}(\xi) = 1 - \xi^{\alpha_4} - (\xi^{\alpha_4})^3. $$
We can therefore conclude that the irreducible characters for the induced representation of $ \boldsymbol{\mu}_n $ on $ H^1(\mathcal{Y}_k, \mathcal{O}_{\mathcal{Y}_k}) $ are $ \chi_1(\xi) = \xi^{\alpha_4} $ and $ \chi_2(\xi) = \xi^{3 \alpha_4} $. 

The irreducible characters for the representation of $ \boldsymbol{\mu}_n $ on $ T_{\mathcal{J}'_k,0} $ induced by the action $ [\xi](\pi') = \xi^{-1} \pi' $ on $R'$ are the inverse characters of these, $ \chi_1^{-1}(\xi) = \xi^{ - \alpha_4} $ and $ \chi_2^{-1}(\xi) = \xi^{ - 3 \alpha_4} $. It is easily seen that $ [- \alpha_4]_n = (n-1)/4 $, and that $ [- 3 \alpha_4]_n = (3n-3)/4 $. Hence the jumps occur at the limits $1/4$ and $ 3/4 $ of these expressions $ (n-1)/4n $ and $ (3n-3)/4n $ as $ n $ goes to infinity.


\end{ex}

\begin{table}[htb]\caption{Genus $2$}\label{table 2}
\begin{tabular}{|c|c|c|c|c|c|c|} 
\hline
Fibertype & $3$ & $4$ & $5$ & $6$ & $7$ & $8$ \\
\hline
Jumps & $ 1/6 $, $ 3/6 $ & $ 1/4 $, $ 3/4 $ & $ 1/4 $, $ 3/4 $ & $ 3/12 $, $ 10/12 $ & $ 3/10 $, $ 9/10 $ & $ 1/5 $, $ 3/5 $ \\
\hline  
\end{tabular} \vspace{3 mm}

\begin{tabular}{|c|c|c|c|c|c|c|c|c|}
\hline
$9$ & $10$ & $11$ & $12$ & $ 13 $ & $ 15 $ & $ 16 $ & $ 17 $ & $18$ \\
\hline
$ 0 $, $ 1/2 $ & $ 0 $, $ 2/3 $ & $ 0 $, $ 3/4 $ & $ 0 $, $ 1/2 $ & $ 0 $, $ 0 $ & $ 0 $, $ 1/2 $ & $ 0 $, $ 2/3 $ & $ 1/6 $, $ 5/6 $ & $1/6$, $ 5/6 $ \\
\hline
\end{tabular} \vspace{3 mm}

\begin{tabular}{|c|c|c|c|c|c|c|}
\hline
$ 19 $ & $20$ & $21$ & $22$ & $ 23 $ &  $ 24a $ & $ 24 $\\
\hline
 $ 4/6 $, $ 5/6 $ & $ 7/10 $, $9/10$ & $ 3/5 $, $ 4/5 $ & $5/8$, $7/8$ & $ 2/4 $, $ 3/4 $ & $ 1/4 $, $ 3/4 $ & $ 1/4 $, $ 3/4 $ \\ 
\hline
\end{tabular} \vspace{3 mm}

\begin{tabular}{|c|c|c|c|c|c|c|c|c|}
\hline
$25$ & $26$ & $27$ & $28$ & $30$ & $31$ & $32$ & $33$ \\
\hline
$ 3/6 $, $ 5/6 $ & $ 2/6 $, $ 5/6 $ & $ 3/8 $, $ 7/8 $ & $ 5/12 $, $ 11/12 $ & $ 0 $, $ 3/4 $ & $ 1/6 $, $ 4/6 $ & $ 0 $, $ 1/2 $ & $ 1/2 $, $ 1/2 $ \\
\hline
\end{tabular} \vspace{3 mm}

\begin{tabular}{|c|c|c|c|c|c|c|c|}
\hline
$34$ & $35$ & $36$ & $37$ & $38$ & $40$ & $41$ & $41_a$ \\
\hline
$ 1/6 $, $ 2/6 $ & $ 0 $, $ 1/4 $ & $ 1/5 $, $ 2/5 $ & $ 0 $, $ 1/3 $ & $ 1/4 $, $ 2/4 $ & $ 0 $, $ 0 $ & $ 0 $, $ 1/2 $ & $ 0 $, $1/2 $\\
\hline
\end{tabular} \vspace{3 mm}

\begin{tabular}{|c|c|c|c|}
\hline
$41_b$ & $41_c$ & $43$ & $44$\\
\hline
$0$, $ 2/3 $ & $0$, $3/4$ & $ 1/3 $, $ 2/3 $ & $ 2/5 $, $ 4/5 $\\
\hline
\end{tabular}
\end{table}

\begin{table}[htb]\caption{Genus $2$, $1_{KOD}$ and $14_{KOD}$}\label{table 3}
\begin{tabular}{|c|c|c|c|c|c|c|c|}
\hline
$KOD$ & $(I_n)$ & $(I_n^*)$ & $(II)^*$ & $(III)$ & $(III)^*$ & $(IV)$ & $(IV)^*$ \\
\hline
Jumps & 0, 0 & 0, 1/2 & 0, 5/6 & 0, 1/4 & 0, 3/4 & 0, 1/3 & 0, 2/3 \\
\hline
\end{tabular}
\end{table}

\subsection{Special fibertypes for $g=2$}\label{special types}
In the Ogg-classification, there are some fiber types that are built up from $g=1$ fiber types. We treat these separately in this section.

Let $KOD$ denote any of the seven specific fiber types $ I_n $, $ I_n^* $, $ II^* $, $ III $, $ III^* $, $ IV $ and $ IV^* $ from the Kodaira-classification. The fiber types $1$ and $14$ in the Ogg-classification are constructed in a certain way from choosing a type $KOD$. The various fibers obtained in this way are denoted by $1_{KOD}$ and $ 14_{KOD} $, and the associated jumps are listed in Table \ref{table 3} (they are actually the same for $1_{KOD}$ and $ 14_{KOD} $).

The fiber types $2$ and $39$ in the Ogg-classification are built up in a certain way by choosing \emph{two} types $ KOD_a $ and $ KOD_b $ from the seven specific Kodaira types $ I_n $, $ I_n^* $, $ II^* $, $ III $, $ III^* $, $ IV $ and $ IV^* $. For each of the 28 possible unordered pairs $ (KOD_a,KOD_b) $ made from this list, we get a fiber type denoted here by $ 2_{KOD_a,KOD_b} $ and $ 39_{KOD_a,KOD_b} $, respectively. Let $j(KOD)$ denote the jump associated to fiber type $KOD$ in Table \ref{table 1}. Then we have that the jumps for the fiber type $ 2_{(KOD_a,KOD_b)} $ are $ j(KOD_a) $ and $ j(KOD_b) $, and likewise, we find that the jumps for type $ 39_{KOD_a,KOD_b} $ are $ j(KOD_a) $ and $ j(KOD_b) $.

The fiber types $ 29 $, $ 29_a $ and $ 42 $ also break up into several cases in a similar way as in the cases mentioned above. Let $KOD' $ denote any of the three specific fiber types $ II^* $, $ III^* $ and $ IV^* $. The jumps for $ 29_{KOD'_a,KOD'_b} $ are listed in Table \ref{table 5}, the jumps for $ 29a_{KOD'} $ are listed in Table \ref{table 6} and the jumps for $ 42_{KOD'} $ are listed in Table \ref{table 7}.

We refer to \cite{Ogg} for precise details regarding the constructions mentioned here.



\begin{table}[htb]\caption{Genus $2$, $29_{KOD'_a, KOD'_b}$}\label{table 5}
\begin{tabular}{|c|c|c|c|}
\hline
$KOD'_a, KOD'_b$ & $(II)^*,(II)^*$ & $(II)^*,(III)^*$ & $(II)^*,(IV)^*$ \\
\hline
Jumps & 5/6, 5/6 & 9/12, 10/12 & 4/6, 5/6 \\
\hline
\end{tabular}
\end{table}

\begin{table}[htb]
\begin{tabular}{|c|c|c|}
\hline
$ (III)^*, (III)^* $ & $ (III)^*, (IV)^*$ & $ (IV)^*, (IV)^* $\\ 
\hline
3/4, 3/4 & 8/12, 9/12 & 2/3, 2/3 \\
\hline
\end{tabular}
\end{table}

\begin{table}[htb]\caption{Genus $2$, $29a_{KOD'}$}\label{table 6}
\begin{tabular}{|c|c|c|c|} 
\hline
$KOD'$ & $(II)^*$ & $(III)^*$ & $(IV)^*$ \\
\hline
Jumps & $ 3/6 $, $ 5/6 $ & $ 2/4 $, $ 3/4 $ & $ 3/6 $, $ 4/6 $\\
\hline
\end{tabular}
\end{table}

\begin{table}[htb]\caption{Genus $2$, $42_{KOD'}$}\label{table 7}
\begin{tabular}{|c|c|c|c|} 
\hline
$KOD'$ & $(II)^*$ & $(III)^*$ & $(IV)^*$\\
\hline
Jumps & $ 2/6 $, $ 5/6 $ & $ 4/12 $, $ 9/12 $ & $ 1/3 $, $ 2/3 $\\
\hline
\end{tabular}
\end{table}

\subsection{Final remarks and comments}
It would be interesting to know, for a curve $X/K$, the significance of the numerators and the denominators of the jumps in the filtration $ \{ \mathcal{F}^a \} $ of $ \mathcal{J}_k $, where $ \mathcal{J} $ is the N\'eron model of the Jacobian of $X$. 

One could also try to obtain a closed formula for the irreducible characters of the representation of $ \boldsymbol{\mu}_n $ on $ H^1(\mathcal{Y}_k, \mathcal{O}_{\mathcal{Y}_k}) $, where $ \mathcal{Y} $ is the minimal desingularization of $ \mathcal{X}_{S'} $, and $ S'/S $ is tamely ramified of degree prime to the least common multiple of the multiplicities of the irreducible components of $ \mathcal{X} $. It seems clear that such a formula would reflect combinatorial properties of the intersection graph of $ \mathcal{X}_k $.

We do not know if our results remain true in the case where the minimal SNC-model $ \mathcal{X} $ of $X/K$ does not fulfill our initial assumptions, that is, if distinct components of $ \mathcal{X}_k $ with multiplicities divisible by $p$ intersect nontrivially. The main problem is the lack of a good description of the minimal desingularization of $ \mathcal{X}_{S'} $, where $ S'/S $ is a tame extension.  

Finally, we think it would be interesting to study these filtrations for N\'eron models of abelian varieties that are not Jacobians. In that case, it is not so clear what kind of data would suffice in order to determine the jumps.

\bibliographystyle{plain}
\bibliography{algbib}

\end{document}